\documentclass[times,final]{elsarticle}

\usepackage{amsmath} 
\usepackage{amssymb}

\usepackage{enumerate}

\usepackage{jcomp}
\usepackage{framed,multirow}

\usepackage{latexsym}

\usepackage{wasysym}
\usepackage{dashrule}
\usepackage{hyperref}

\usepackage{url}
\usepackage{xcolor}
\usepackage{ulem}

\usepackage[ruled,lined,shortend]{algorithm2e}

\definecolor{newcolor}{rgb}{.8,.349,.1}

\journal{Journal of Computational Physics}

\begin{document}

\verso{C. S. Huang \textit{et al.}}

\begin{frontmatter}

\title{A 3D sharp and conservative VOF method for modeling the contact line dynamics with hysteresis on complex boundaries}%

\author[1]{Chong-Sen \snm{Huang}}
\author[2,3]{Tian-Yang \snm{Han}}
 
\author[1]{Jie \snm{Zhang}\corref{cor1}}
\cortext[cor1]{Corresponding author:}
\ead{j_zhang@xjtu.edu.cn}

\author[2]{Ming-Jiu \snm{Ni}\corref{cor2}}
\cortext[cor2]{Corresponding author:}
\ead{mjni@ucas.edu.cn}

\address[1]{ State Key Laboratory for Strength and Vibration of Mechanical Structures, School of Aerospace, Xi'an Jiaotong University, Xi'an, Shaanxi 710049, China}
\address[2]{ School of Engineering Science, University of Chinese Academy of Sciences, Beijing 101408, China}
\address[3]{ Sorbonne Université, CNRS, Institut Jean Le Rond d'Alembert, F-75005 Paris, France}

\bibliographystyle{unsrt}

\begin{abstract}
We propose a sharp and conservative 3D numerical method for simulating moving contact lines on complex geometries, developed within a coupled geometric Volume-of-Fluid (VOF) and embedded boundary framework. The first major contribution is a modified VOF advection and reconstruction scheme specifically designed for mixed cells containing liquid, gas, and solid phases. This formulation ensures strict local mass conservation in the presence of arbitrarily shaped embedded boundaries. To overcome the severe time-step limitation caused by small cut cells, a redistribution advection strategy is introduced, which completely removes the CFL constraint while preserving both local and global volume conservation. 
The second key contribution is a novel 3D contact angle imposition technique built upon the height function framework. By incorporating a pre-fitting paraboloid procedure, the method achieves robust curvature estimation and accurate enforcement of contact angle conditions on irregular solid surfaces. In addition, contact angle hysteresis is modeled to capture more realistic wetting dynamics. 
A series of challenging benchmark tests have been conducted to demonstrate the accuracy, robustness, and superiority of the proposed method compared with existing sharp-interface approaches. This study, for the first time, establishes a fully geometric and conservative VOF-based scheme capable of accurately resolving contact line dynamics on arbitrarily complex 3D surfaces.
\end{abstract}


\end{frontmatter}



\section{Introduction}\label{sec:introduction}

\begin{figure}
    \includegraphics[height=0.36\textwidth]{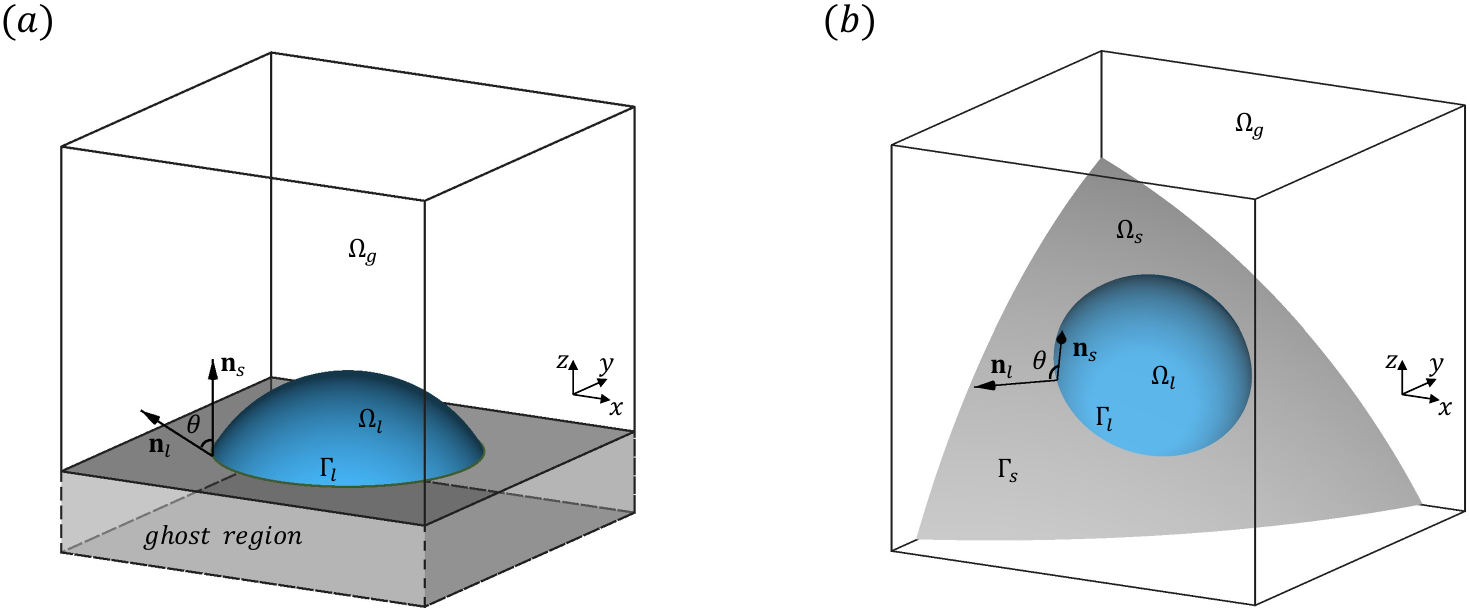}
    \centering
    \caption{
    Schematic illustration of the moving contact line problem on (a) a flat solid substrate (regular boundary) and (b) an irregular solid substrate (complex geometric boundary). $\Omega_g$, $\Omega_l$, and $\Omega_s$ denote the gas, liquid, and solid phases, respectively. $\Gamma_l$ and $\Gamma_s$ indicate the liquid/gas interface and the solid surface, with $\textbf{n}_l$ and $\textbf{n}_s$ representing their corresponding normal vectors. The symbol $\theta$ denotes the contact angle.
    }
    \label{fig:description}
\end{figure}

Contact line dynamics, referring to the motion of the liquid/gas/solid intersection, play a fundamental role in a wide range of engineering and scientific processes, such as water electrolysis, thin-film manipulation on microstructured substrates, and droplet deposition in inkjet printing. These phenomena involve complex multiphase interactions near solid boundaries and require accurate modeling of fluid interfaces subject to both geometric and physical constraints.

As illustrated in Fig.~\ref{fig:description}(a), previous studies have primarily examined droplet spreading on smooth, planar solid substrates. In these cases, the solid boundary is typically aligned with one of the Cartesian domain boundaries, hereafter referred to as a regular boundary. In contrast, practical applications increasingly demand the capability to simulate more general scenarios, such as droplet motion over curved or arbitrarily shaped substrates. These complex geometric boundaries, illustrated in Fig.~\ref{fig:description}(b), do not conform to the computational mesh and are irregularly embedded within the simulation domain.

The introduction of complex solid geometries presents significant numerical challenges, particularly in resolving the local fluid dynamics near solid boundaries. In this study, the term \textit{mixed cells} refers to computational cells intersected by the solid boundary $\Gamma_s$, which therefore contain partial volumes of both fluid and solid. Among these, special attention is given to \textit{contact line cells}, defined as mixed cells where the liquid/gas interface $\Gamma_l$ intersects the solid surface. These cells are crucial for enforcing the contact angle condition that characterizes the local geometry of the triple-phase junction. 
In the present context, three immiscible domains are considered: the liquid phase $\Omega_l$, the gas phase $\Omega_g$, and the solid phase $\Omega_s$. The interfaces separating these regions are the immobile solid surface $\Gamma_s = \Omega_s \cap (\Omega_l \cup \Omega_g)$ and the deformable liquid/gas interface $\Gamma_l = \Omega_g \cap \Omega_l$. The contact line is defined by the intersection $\Gamma_l \cap \Gamma_s$, where the contact angle $\theta$ is prescribed. Geometrically, this condition satisfies $\cos \theta = \mathbf{n}_l \cdot \mathbf{n}_s$, where $\mathbf{n}_l$ and $\mathbf{n}_s$ denote the unit normal vectors of $\Gamma_l$ and $\Gamma_s$ at the contact line, respectively. Accurate imposition of this angle within contact line cells—particularly when $\Gamma_s$ is arbitrarily oriented—constitutes a central technical challenge of the present work.

Over the past two decades, substantial progress has been achieved in simulating contact line dynamics on regular boundaries using both diffuse-interface and sharp-interface models~\cite{jacqmin2000contact,renardy2001numerical,spelt2005level,afkhami2009height}. These efforts encompass two- and three-dimensional implementations within various interface-capturing frameworks, including the Level Set, Phase Field, and Volume-of-Fluid (VOF) methods. However, extending these approaches to complex geometric boundaries introduces additional challenges, particularly in imposing accurate physical boundary conditions on non-aligned and curved surfaces. 
To overcome these difficulties, several embedded-boundary strategies have been proposed. Among them, the immersed boundary method (IBM)~\cite{peskin1972flow,mittal2005immersed} and the embedded boundary method (EBM)~\cite{johansen1998cartesian,schwartz2006cartesian} are two of the most widely adopted techniques. These frameworks discretize the governing equations on Cartesian grids by incorporating immersed geometries through modified fluxes or source terms near solid interfaces. When combined with front-tracking or VOF-type schemes, they have enabled simulations of contact line motion over irregular surfaces. Despite these advancements, two fundamental numerical challenges remain unresolved and constitute the central focus of the present study.

The first challenge concerns the accurate advection of the liquid/gas interface $\Gamma_l$ in the presence of complex solid boundaries $\Gamma_s$. Directly extending interface-advection techniques originally developed for regular cells to mixed cells is nontrivial, particularly when local and global mass conservation must be maintained. For example, Liu et al.~\cite{liu2015diffuse} coupled a phase-field method with the IBM, a formulation readily applicable to complex 3D geometries. The diffuse nature of the phase-field model introduces a finite interfacial thickness, which facilitates interface advection near curved boundaries; however, the loss of interface sharpness limits its applicability in scenarios requiring precise interface reconstruction. 
In contrast, O'Brien et al.~\cite{o2018volume} employed an algebraic VOF method combined with the IBM in 2D; however, their approach exhibited poor volume conservation even for regular boundaries.
More recent advances include the hybrid geometric VOF/IBM method proposed by Patel et al.~\cite{patel2017coupled} and the geometric VOF/EBM framework developed by Tavares et al.~\cite{tavares2024coupled}, both extended to 3D. Nevertheless, these studies reported noticeable fluid accumulation within solid regions over time, particularly in mixed cells intersected by both $\Gamma_l$ and $\Gamma_s$. This issue reflects the inherent difficulty of achieving strict mass conservation in sharp-interface formulations near complex geometric boundaries.

Our previous 2D study~\cite{huang20252d} addressed this limitation by developing a conservative geometric VOF method applicable to irregular solid boundaries. The method ensured local mass conservation within mixed cells while effectively preventing nonphysical fluid accumulation. Although it introduced a moderate computational overhead of approximately 32\%, it represented the first demonstration of accurate and conservative advection within a geometric VOF framework incorporating embedded geometries. A primary objective of the present work is to extend this 2D formulation to 3D, where both geometric representation and algorithmic implementation become considerably more complex. As demonstrated in this study, this extension is nontrivial and constitutes one of the core technical contributions of our work.

The second challenge involves the accurate enforcement of contact angle boundary conditions on complex solid surfaces, which directly affects interfacial curvature estimation and, consequently, the evaluation of surface tension forces near the contact line. Although considerable progress has been achieved in this area using both diffuse-interface~\cite{ding2007wetting} and sharp-interface formulations~\cite{renardy2001numerical,spelt2005level,afkhami2009height}, most existing methods have been developed for regular boundaries and rely on ghost-cell extrapolation to impose contact angle constraints. Among these, the height function (HF) method has become a widely adopted technique for computing interfacial curvature within the VOF framework, allowing accurate contact angle enforcement on flat, regular surfaces through geometric reconstruction~\cite{afkhami2009height,han2021consistent}. 
However, extending HF-based techniques to complex boundaries presents substantial challenges. For instance, Tavares et al.~\cite{tavares2024coupled} extended the liquid/gas interface into the solid region as a straight line in 2D or as a planar surface in 3D to generate ghost values, thereby approximating the contact angle condition. While this linear extension is conceptually simple and computationally efficient, it introduces significant errors for extreme contact angles (e.g., $\theta < 45^\circ$ or $\theta > 135^\circ$). In our recent 2D study~\cite{huang20252d}, we demonstrated that these inaccuracies stem from the sensitivity of the linear fitting to the relative position of the embedded boundary within mixed cells. This sensitivity often results in inconsistent contact line motion and inaccurate angle enforcement outside a limited range of $\theta$ values.

To overcome this limitation, we previously proposed a quadratic fitting strategy in 2D, which substantially enhanced both the accuracy and robustness of the imposed contact angle condition. The method reconstructs the liquid/gas interface using a parabolic curve that satisfies the geometric constraints defined by the surrounding interface configuration and the prescribed contact angle. Although this approach proved effective in 2D, its extension to 3D introduces significant algorithmic complexities, including the selection of appropriate interpolation stencils and the preservation of curvature consistency across mixed cells. These challenges have motivated the development of a new 3D algorithm capable of accurately and robustly enforcing contact angle conditions on arbitrarily complex solid boundaries.

In summary, this study presents a conservative, accurate, and robust numerical framework for simulating moving contact line dynamics on complex 3D geometries. The geometric volume-of-fluid (VOF) method is employed to track the liquid/gas interface $\Gamma_l$, while the height function (HF) method is utilized to compute the surface tension force. Embedded solid boundaries $\Gamma_s$ are represented using the embedded boundary method (EBM), enabling both $\Gamma_l$ and $\Gamma_s$ to be sharply resolved within a Cartesian grid. This unified treatment of sharp interfaces represents a central innovation of the present work. 
Broadly, this study makes two primary contributions:
\begin{enumerate}
    \item A mass-conservative 3D interface advection scheme is developed for cells intersected by $\Gamma_s$ (i.e., mixed cells). Building upon our previous 2D formulation~\cite{huang20252d}, the proposed scheme ensures accurate volume transport while removing the stringent time-step constraint imposed by the mixed cells with small fluid volume fractions. The method is fully geometric, computationally efficient, and robust across a wide range of interface configurations.
    \item A novel 3D HF-based contact angle formulation is introduced to accurately impose contact angle conditions on arbitrarily complex embedded boundaries. This approach extends the parabola-fitting strategy previously developed in 2D~\cite{huang20252d} and incorporates a consistent procedure for modeling contact angle hysteresis based on local volume conservation principles.
\end{enumerate}

The remainder of this paper is organized as follows. Section~\ref{sec:ge-bm} introduces the governing equations and discusses the numerical challenges associated with embedded boundaries, along with a brief review of our previous 2D algorithm to motivate its 3D extension. Section~\ref{sec:complexVOFgeometric} presents the 3D geometric VOF advection scheme for mixed cells and demonstrates its conservative properties. Section~\ref{sec:contactanglecondition} details the 3D HF-based approach for enforcing contact angle conditions on complex boundaries, including the treatment of hysteresis effects. Section~\ref{sec:results} provides a series of benchmark tests validating the accuracy and robustness of the proposed methods. Finally, concluding remarks are offered in Section~\ref{sec:conclusion}.

\section{Governing equations and basic methods}\label{sec:ge-bm}

The configuration considered in this study is illustrated in Fig.~\ref{fig:description}(b). The computational domain consists of three immiscible phases: gas $\Omega_g$, liquid $\Omega_l$, and solid $\Omega_s$, which are separated by two sharp interfaces—the deformable liquid/gas interface $\Gamma_l$ and the stationary solid boundary $\Gamma_s$. The liquid/gas interface $\Gamma_l$ is tracked using the Volume-of-Fluid (VOF) method, where a scalar field $c(\mathbf{x}, t) \in [0,1]$ denotes the liquid volume fraction of $\Omega_l$ within each computational cell. It is defined as
\[
c = \frac{|\Omega_l|}{|\Omega_g| + |\Omega_l| + |\Omega_s|}.
\]
This field evolves in time to capture the motion and deformation of $\Gamma_l$, and the interface is reconstructed using the Piecewise Linear Interface Calculation (PLIC) method, which provides sharp geometric resolution.
To represent complex solid geometries, the solid surface $\Gamma_s$ is embedded within the Cartesian grid and treated using a volume-fraction-based approach. In this framework, cells intersected by $\Gamma_s$ are referred to as \textit{mixed cells}, as they contain both solid and fluid phases. A separate solid volume-fraction field $c_s(\mathbf{x}) \in [0,1]$ is defined as
\[
c_s = \frac{|\Omega_g| + |\Omega_l|}{|\Omega_g| + |\Omega_l| + |\Omega_s|},
\]
which quantifies the fraction of each cell accessible to the fluid. Unlike $c$, which evolves in time, $c_s$ is initialized at the beginning of the simulation and remains constant, as the solid boundary is assumed to be stationary. Similar to $\Gamma_l$, the embedded solid boundary $\Gamma_s$ is reconstructed using planar segments within the mixed cells.

\subsection{Governing equations}\label{sec:governingequation}

We consider incompressible, Newtonian, and viscous two-phase flows featuring a moving liquid/gas interface $\Gamma_l$ and a stationary solid boundary $\Gamma_s$, as illustrated in Fig.~\ref{fig:description}(b). The governing equations for the conservation of mass and momentum are given by

\begin{equation}\label{eq:divu}
\nabla \cdot \textbf{u} = 0,
\end{equation}
\begin{equation}\label{eq:ns}
\rho \left(\frac{\partial \textbf{u}}{\partial t} + \textbf{u} \cdot \nabla \textbf{u}\right)
= -\nabla p + \nabla \cdot (2\mu \mathbb{D}) + \sigma \kappa \delta_s \textbf{n}_l + \rho \textbf{g},
\end{equation}
where $\textbf{u}$ and $p$ denote the velocity and pressure fields, respectively. Here, $t$, $\rho$, $\mu$, and $\textbf{g}$ represent time, density, dynamic viscosity, and gravitational acceleration. The deformation-rate tensor is defined as $\mathbb{D} = (\nabla \textbf{u} + \nabla \textbf{u}^\mathrm{T}) / 2$. The surface tension force is expressed as $\sigma \kappa \delta_s \textbf{n}_l$, where $\textbf{n}_l$ is the unit normal vector to $\Gamma_l$, $\kappa$ is the interface curvature, and $\delta_s$ is the Dirac delta function that localizes the force at the interface.

The evolution of the liquid/gas interface is governed by the advection equation for the volume fraction $c(\mathbf{x},t)$:
\begin{equation}\label{eq:advectiongoverned}
\frac{\partial c}{\partial t} + \textbf{u} \cdot \nabla c = 0,
\end{equation}
where $c=1$ indicates pure liquid, $c=0$ indicates pure gas, and $0 < c < 1$ defines interfacial cells. The volume fraction $c$ is then used to evaluate the effective material properties of the two-phase mixture:
\begin{equation}\label{eq:rho}
   \rho(c) = c \rho_l + (1 - c)\rho_g, 
\end{equation}
\begin{equation}\label{eq:mu} 
   \mu(c)  = c \mu_l + (1 - c)\mu_g, 
\end{equation}
where $\rho_l$ and $\rho_g$ ($\mu_l$ and $\mu_g$) are the densities (dynamic viscosities) of the liquid and gas, respectively.

The present method is implemented in the open-source flow solver \textit{Basilisk}, developed by S.~Popinet~\cite{popinet2009accurate}. In \textit{Basilisk}, all variables are discretized on a collocated grid, and time advancement is performed using an approximate projection method~\cite{Chorin1968Numerical,David2001Accurate}, which achieves second-order accuracy in both space and time. The surface tension term in Eq.~\eqref{eq:ns} is evaluated using the continuum surface force (CSF) model~\cite{brackbill1992continuum}, which expresses $\sigma \kappa \delta_s \textbf{n}_l$ as $\sigma \kappa \nabla c$, where $\nabla c$ is a representation of the interface normal vector, and $\kappa$ denotes the curvature computed by the height function (HF) method. The liquid/gas interface $\Gamma_l$ is geometrically reconstructed from $c$ using a piecewise linear interface construction (PLIC) scheme consistent with the geometric VOF framework.

In the original \textit{Basilisk}, contact angle boundary conditions on flat solid surfaces (regular boundaries) are imposed using a ghost-cell method~\cite{afkhami2009height}. In this approach, ghost HF values are assigned to cells adjacent to the solid boundary according to the prescribed contact angle $\theta$, which in turn influences the curvature and normal vectors computed. Although $\theta$ does not explicitly appear in the governing equations, it directly affects the surface tension force at the contact line. To accommodate irregular geometries, \textit{Basilisk} also incorporates embedded boundaries~\cite{johansen1998cartesian}, which represent complex solids through planar reconstructions within Cartesian grids. The no-slip condition is enforced by introducing additional fluxes across the reconstructed solid surface $\Gamma_s$ during the momentum discretization.

The next subsection reviews the original geometric VOF and HF methods for regular boundaries as implemented in \textit{Basilisk}, followed by a brief summary of our recent 2D study~\cite{huang20252d} addressing complex embedded boundaries. We then discuss the main challenges associated with extending this framework to 3D.

\subsection{VOF advection method}\label{sec:regularVOFgeometric}

\begin{figure}
    \includegraphics[width=1\textwidth]{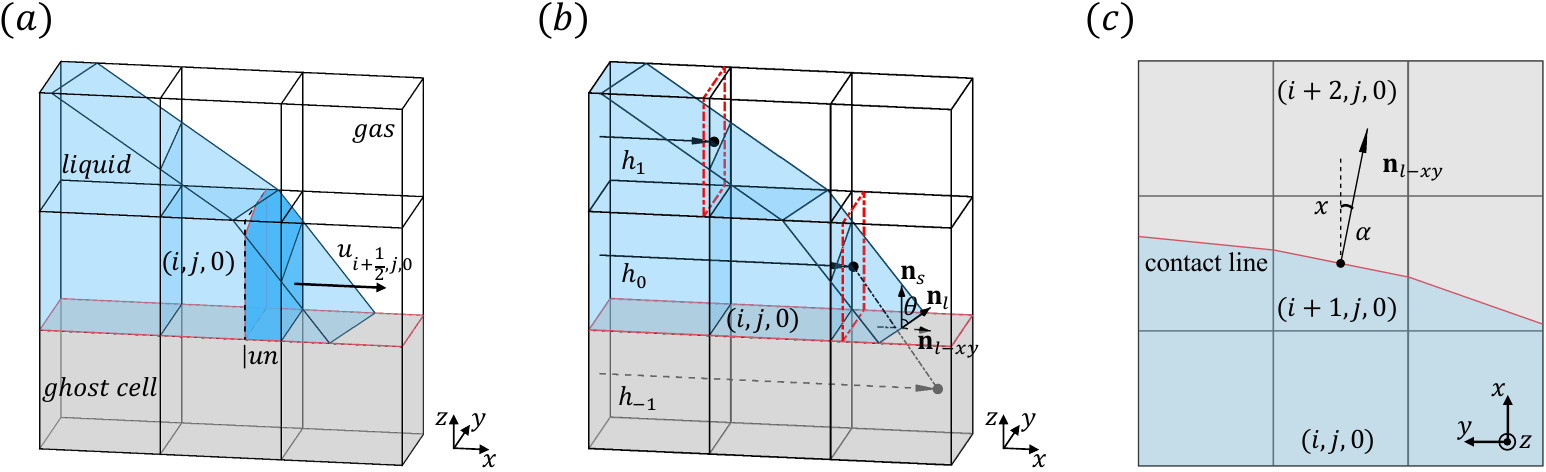}
    \centering
    \caption{
    (a) Schematic of the geometric VOF advection method on regular boundaries. (b) Definition of the height function in 3D on regular boundaries, including its evaluation in ghost cells. (c) Contact line contour on the solid boundary.
    }
    \label{fig:VOF/HF}
\end{figure}

The liquid volume fraction $c$ is updated by solving the advection equation~\eqref{eq:advectiongoverned} using a geometric, directionally split scheme~\cite{popinet2009accurate}, in which the interface is advected sequentially along the $x$-, $y$-, and $z$-directions. To account for the artificial compressibility introduced by directional splitting, Eq.~\eqref{eq:advectiongoverned} is reformulated as:

\begin{equation}\label{eq:advection1}
\frac{\partial c}{\partial t} + \nabla \cdot (c\textbf{u}) - c \nabla \cdot \textbf{u} = 0.
\end{equation}
Integrating Eq.~\eqref{eq:advection1} over a control volume $\Omega$ yields the volume-integral form:
\begin{equation}\label{eq:advection-int}
\frac{\partial}{\partial t} \int_{\Omega} c \, d\Omega
+ \oint_{\partial \Omega} c \, \textbf{u} \cdot \textbf{n} \, dS
- c \oint_{\partial \Omega} \textbf{u} \cdot \textbf{n} \, dS = 0,
\end{equation}
\begin{equation}\label{eq:advection discretization scheme1}
    c^* = c^{n - \frac{1}{2}} + \frac{dt}{\Delta} \left[
    \left(F_{i - \frac{1}{2},j,k}^{n - \frac{1}{2}} - F_{i + \frac{1}{2},j,k}^{n - \frac{1}{2}} \right)
    - c_c \left(u_{i - \frac{1}{2},j,k}^{n - \frac{1}{2}} - u_{i + \frac{1}{2},j,k}^{n - \frac{1}{2}} \right)
    \right],
\end{equation}
where $(i,j,k)$ is the index of the control volume, and $(i \pm 1/2,j,k)$ denote the left and right cell faces in the $x$-direction. The correction term $c_c$, proposed by Weymouth et al.~\cite{weymouth2010conservative}, accounts for the divergence of velocity and ensures local volume conservation, preventing artificial over- or under-filling of cells.

The accuracy of the scheme highly depends on evaluating the fluxes $F = c \, \textbf{u}_f \cdot \textbf{n}_f$ across each cell face, where subscript $f$ denotes a face-centered quantity. These fluxes are computed geometrically based on the reconstructed interface shape. For example, as illustrated in Fig.~\ref{fig:VOF/HF}(a), the flux through the right face of cell $(i,j,0)$ is calculated as:
\begin{equation}\label{eq:cf-normal}
    F_{i + \frac{1}{2}, j, 0} = \frac{V \cdot u_{i + \frac{1}{2}, j, 0}}{\Delta^2 \cdot un},
\end{equation}
where $V$ is the advected liquid volume (dark blue region), and $un = u_{i + \frac{1}{2}, j, 0} \cdot dt$ is the displacement over one time step.

The flux calculation requires the local interface normal vector $\textbf{n}_l$, which is typically computed using the height function (HF) method, as described in the next section. After the advection step in the $x$-direction, the interface $\Gamma_l$ is reconstructed with the PLIC method based on the updated intermediate volume fraction $c^*$. The same procedure is then applied sequentially in the $y$- and $z$-directions, resulting in the final volume fraction $c^{n + 1/2}$ after all substeps are completed.

\subsection{Height function method}\label{sec:surfacetensionterm}

A key challenge in solving the momentum equation~\eqref{eq:ns} lies in accurately evaluating the surface tension force $\sigma \kappa \delta_s \mathbf{n}_l$, which requires precise estimation of both the curvature $\kappa$ and the interface normal vector $\mathbf{n}_l$. This is particularly crucial in capillary-dominated flows, where surface tension governs the interfacial dynamics. Although the VOF method provides excellent mass conservation, it often yields inaccurate curvature and normal estimates due to the discontinuity of the volume fraction field. To overcome this limitation, the height function (HF) method is employed in \textit{Basilisk}, providing a continuous representation from which $\kappa$ can be computed with high accuracy.

As illustrated in Fig.~\ref{fig:VOF/HF}(b), the HF value $h_0$ at the first fluid layer adjacent to the bottom wall is obtained by summing the volume fractions $c_{i,j,0}$ along the $x$-direction:
\begin{equation}\label{eq:hsum}
    h_0 = \sum_i c_{i,j,0} \, \Delta,
\end{equation}
where the stencil typically extends from $(i-4,j,0)$ to $(i+4,j,0)$, unless truncated by the domain boundaries. The local interface curvature at cell $(i,j,0)$ is then evaluated using central finite differences applied to the HF field~\cite{popinet2009accurate}:
\begin{equation}\label{eq:h1-kappa}
    \kappa = \frac{h_{yy} + h_{zz} + h_{yy} h_z^2 + h_{zz} h_y^2 - 2 h_{yz} h_y h_z}{(1 + h_y^2 + h_z^2)^{3/2}},
\end{equation}
with derivatives computed as:
\begin{equation}\label{eq:h1-kappa-2}
\begin{aligned}
h_y &= \frac{h_{i,j+1,0} - h_{i,j-1,0}}{2\Delta}, \quad
h_{yy} = \frac{h_{i,j+1,0} - 2h_{i,j,0} + h_{i,j-1,0}}{\Delta^2}, \\
h_{yz} &= \frac{h_{i,j+1,1} + h_{i,j-1,-1} - h_{i,j+1,-1} - h_{i,j-1,1}}{4\Delta^2}.
\end{aligned}
\end{equation}

HF values in fluid cells (e.g., $h_{i,j,0}$) are computed directly, whereas ghost-cell values (e.g., $h_{i,j,-1}$) are extrapolated using the contact angle condition. As shown in Fig.~\ref{fig:VOF/HF}(b), the ghost value $h_{-1}$ is determined as
\begin{equation}\label{eq:standard-interpolation}
    h_{-1} = h_0 + \frac{\Delta}{\tan \theta \cdot \cos \alpha},
\end{equation}
where $\theta$ is the prescribed contact angle, and $\alpha$ denotes the inclination between the HF orientation and the projection of the contact line normal vector $\mathbf{n}_{l-xy}$ onto the $xy$-plane. In 2D, this expression simplifies to $h_{-1} = h_0 + \Delta / \tan \theta$~\cite{huang20252d}. Although the original method of Afkhami et al.~\cite{afkhami2008height} did not specify how to compute $\cos \alpha$, \textit{Basilisk} evaluates it at $z = \Delta/2$. A more accurate approach proposed by Han et al.~\cite{han2021consistent} fits a paraboloid at the contact line and evaluates $\cos \alpha$ directly at $z = 0$, thereby improving curvature estimation and achieving first-order convergence.

The interface normal vector $\mathbf{n}_l$ plays a crucial role in several components of the numerical framework. It appears first in the surface tension term of Eq.~\eqref{eq:ns}. In \textit{Basilisk}, this vector is approximated by the normalized gradient of the volume fraction:
\begin{equation}
    \mathbf{n}_{l,c} = \frac{\nabla c}{|\nabla c|},
\end{equation}
which is evaluated over a $3 \times 3 \times 3$ stencil. This formulation ensures consistency with the discretization of the pressure gradient and leads to a force-balanced CSF model~\cite{francois2006balanced}, effectively reducing spurious currents.

Second, $\mathbf{n}_l$ is also required for the geometric advection of $\Gamma_l$ (as described in Section~\ref{sec:regularVOFgeometric}) and for applying contact angle conditions at contact line cells. For these purposes, HF-based normals are preferred owing to their higher accuracy:
\begin{equation}\label{eq:hn}
    \mathbf{n}_{l,h} = (-1, h_y, h_z),
\end{equation}
which are normalized after construction. This computation additionally depends on ghost HF values to complete the finite-difference stencil.

In summary, both $\mathbf{n}_{l,c}$ (fraction gradient-based) and $\mathbf{n}_{l,h}$ (HF-based) are indispensable but serve distinct roles. The former is employed in the surface tension term of the momentum equation to maintain force balance and suppress spurious currents, whereas the latter is used for geometric interface reconstruction and the enforcement of contact angle conditions. This distinction becomes particularly important near complex solid boundaries, where inconsistencies between curvature estimation and contact angle enforcement may lead to significant numerical errors. Resolving this issue is a central objective in developing accurate and robust numerical methods for modeling contact line dynamics.

It is worth noting that Fig.~\ref{fig:VOF/HF}(b) illustrates the HF computed in the $x$-direction. In practice, the HF orientation is dynamically selected according to the local interface geometry and the value of $\alpha$. For moderate contact angles ($45^\circ < \theta < 135^\circ$), a ``horizontal HF'' aligned with the $x$- or $y$-direction is typically employed. In contrast, for extreme contact angles ($\theta < 45^\circ$ or $\theta > 135^\circ$), a ``vertical HF'' aligned with the $z$-direction becomes necessary. However, constructing vertical HFs in a stable and accurate manner is considerably more challenging, which is why the original \textit{Basilisk} implementation supports only horizontal HFs. To overcome this limitation, the present study introduces a consistent and accurate vertical HF formulation capable of handling extreme contact angle conditions, even over complex geometries. The validity and robustness of this extension will be demonstrated in \S~\ref{sec:results}.

In the following section, we briefly review our recent 2D method for modeling contact line dynamics on complex boundaries~\cite{huang20252d}, before discussing the new challenges and strategies involved in extending this approach to 3D.

\subsection{2D contact line model}\label{sec:review2D}

In our recent study~\cite{huang20252d}, we developed a two-dimensional (2D) numerical framework for simulating contact line dynamics on complex boundaries within the geometric VOF methodology. The scheme integrated a 2D flood-fill algorithm to reconstruct the liquid/gas interface $\Gamma_l$ in mixed cells partially intersected by embedded solids $\Gamma_s$, together with a modified advection velocity designed to compensate for flux imbalance in the discretization of Eq.~\eqref{eq:advection discretization scheme1}. To alleviate the time-step constraint associated with small cut cells, a threshold for $c_s$ was introduced, below which the corresponding mixed cells were forcibly removed at the start of the simulation. For curvature estimation, a new 2D height-function (HF) approach was proposed, employing parabolic fitting to reconstruct ghost-cell values and accurately impose contact angle conditions on arbitrarily shaped solid surfaces.

While several components of this framework naturally extend to three dimensions, others require significant reformulation. The key features and their corresponding 3D extensions are summarized below:

\begin{itemize}
    \item \textbf{Flood-fill algorithm.}  
    In 2D, the algorithm effectively reconstructed $\Gamma_l$ within mixed cells by identifying connected fluid regions. In 3D, this approach is generalized into a volume-filling flood algorithm~\cite{dyadechko2005moment}, where the primary challenge lies in determining the correct ordering of the vertices that define convex polyhedra. Furthermore, accurately computing the advected liquid volumes within unstructured 3D mixed cells introduces additional geometric complexity.

    \item \textbf{VOF advection.}  
    In 2D, applying the standard geometric VOF scheme within mixed cells often resulted in inaccurate flux evaluations and poor volume conservation. This limitation persists in 3D, where the coupling between $\Gamma_l$ and the embedded boundaries becomes even more complex. To address this challenge, we introduce a 3D geometric advection method specifically designed to ensure mass conservation in the presence of embedded solid surfaces.

    \item \textbf{Small cut cells.} 
    In 2D, small mixed cells were crudely removed to alleviate time-step constraints. However, this approach becomes impractical in 3D, where the fluid volume fraction scales as $c_s \sim {s_l}^N$ (with $s_l$ denoting the fluid fraction on the cell edge and $N$ the spatial dimension), resulting in a rapid increase in the number of small cells. To overcome this issue, we develop a flux redistribution scheme that completely eliminates the restrictive time-step limitation while enhancing computational efficiency and maintaining numerical accuracy.

    \item \textbf{Height-function method.}  
    In 2D, height-function (HF) values in unstructured fluid domains were obtained through parabolic fitting, allowing accurate evaluation of curvature and enforcement of contact angle conditions on complex surfaces. In 3D, we extend this concept using a paraboloid fitting approach. To further improve robustness and accuracy, a new two-step procedure is introduced: an initial pre-fitting stage identifies suitable interpolation points from an extended stencil, followed by a constrained fitting that explicitly enforces the prescribed contact angle.
\end{itemize}

In summary, the proposed 3D method is not a straightforward extension of our previous 2D framework. Instead, it integrates several newly developed strategies to address challenges that are intrinsic to 3D geometries, particularly those arising in the accurate modeling of interface dynamics near complex solid boundaries. The details of these developments are presented in the following sections.

\section{Geometric VOF advection in mixed cells}\label{sec:complexVOFgeometric}

This section presents a 3D flood-fill algorithm for reconstructing the liquid/gas interface $\Gamma_l$ in the presence of complex solid boundaries $\Gamma_s$. In addition, a flux redistribution scheme is introduced to enhance the robustness of volume-fraction advection in Eq.~\eqref{eq:advection1}, effectively eliminating the stringent time-step restriction imposed by small mixed cells.

\subsection{Reconstruction of $\Gamma_l$ in mixed cells}\label{sec:floodalgorithm}

\begin{figure}
    \includegraphics[width=1\textwidth]{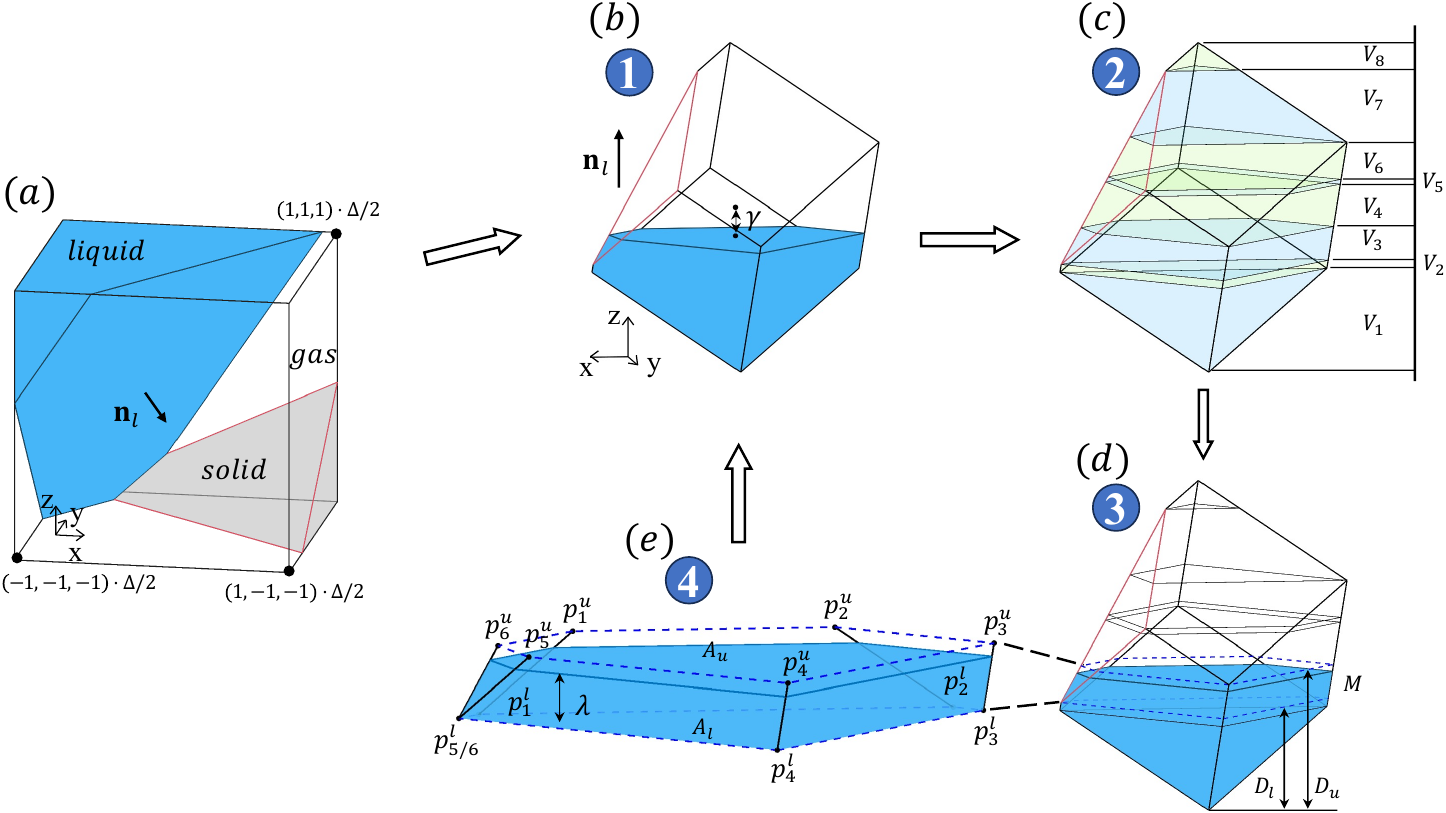}
    \centering
    \caption{
    Schematic of the liquid/gas interface reconstruction within a mixed cell. $\gamma$ denotes the signed distance from the cell center to the interface. $V_i$ represents the volume of each subdivided region. $D_l$ and $D_u$ are the heights from the lowest vertex to the lower and upper planes of the sub-region, respectively. The $M$-region contains the liquid/gas interface. $A_l$ and $A_u$ are the polygonal cross-sectional areas at the lower and upper planes, and $\lambda$ denotes the interface height within the sub-region.
    }
    \label{fig:floodalgorithm}
\end{figure}

In the geometric VOF method implemented in \textit{Basilisk}~\cite{scardovelli1999direct,popinet2009accurate}, the sharp liquid/gas interface $\Gamma_l$ is reconstructed within structured cubic cells using the PLIC scheme~\cite{popinet2009accurate}. The interface is defined by
\begin{equation}\label{eq:linearinterfaceeq}
    \textbf{n}_l \cdot \textbf{x} = \gamma,
\end{equation}
where $\textbf{n}_l$ is the unit normal vector to $\Gamma_l$, and $\textbf{x}$ is the position vector. Given $\textbf{n}_l$ and the volume fraction $c$ (obtained from the advection equation, Eq.~\eqref{eq:advection discretization scheme1}), the scalar $\gamma$ is uniquely determined by requiring that the volume of liquid beneath the reconstructed plane equals $c\Delta^3$. Geometrically, $\gamma$ represents the signed distance from the cell center to the interface. For regular cubic cells, $\gamma$ can be computed analytically. However, when the embedded solid surface $\Gamma_s$ intersects a cell, the fluid region becomes a general convex polyhedron (see Fig.~\ref{fig:floodalgorithm}(a)), rendering the standard PLIC algorithm inapplicable. In such cases, the interface must be reconstructed within the irregular sub-volume. To address this issue, we extend the flood-fill algorithm originally developed in our 2D work~\cite{huang20252d} to the 3D setting. This approach builds on prior work by Dai et al.~\cite{dai2019analytical}, who introduced an analytical VOF method for unstructured meshes. Here, we adapt the method to embedded-boundary Cartesian grids in 3D.

As illustrated in Fig.~\ref{fig:floodalgorithm}(b), the first step of the algorithm rotates the irregular fluid polyhedron so that the normal vector $\textbf{n}_l$ aligns with the positive $z$-axis. The rotation matrix $\mathbb{R}$ satisfies $\mathbb{R} \cdot \textbf{n}_l^\mathrm{T} = (0,0,1)^\mathrm{T}$ and is expressed as
\begin{equation}
    \mathbb{R} = 
    \left[
    \begin{array}{ccc}
    1 - \dfrac{n_{l,x}^2}{1 + n_{l,z}} & -\dfrac{n_{l,x} n_{l,y}}{1 + n_{l,z}} & -n_{l,x} \\
    -\dfrac{n_{l,x} n_{l,y}}{1 + n_{l,z}} & 1 - \dfrac{n_{l,y}^2}{1 + n_{l,z}} & -n_{l,y} \\
    n_{l,x} & n_{l,y} & n_{l,z}
    \end{array}
    \right],
\end{equation}
where $(n_{l,x}, n_{l,y}, n_{l,z})$ denote the components of the unit normal vector $\textbf{n}_l$.

In the second step (Fig.~\ref{fig:floodalgorithm}(c)), the rotated polyhedron is partitioned into sub-regions by inserting planes normal to the $z$-axis that pass through each vertex. Within each resulting slab, the liquid volume $V_i$ is computed using the cross-sectional method~\cite{halsted1904rational}:
\begin{equation}\label{eq:cross-section_area}
    V_i = \frac{D_u - D_l}{6}(A_u + A_l + 4A_{1/2}),
\end{equation}
where $D_u$ and $D_l$ are the heights of the upper and lower bounding planes of the sub-region. $A_u$ and $A_l$ denote the areas of the upper and lower polygonal cross-sections, respectively, and $A_{1/2}$ is the interpolated mid-plane area, evaluated as
\begin{equation}\label{eq:middle_area}
    A_{1/2} = \frac{1}{4}(A_u + A_l) + \frac{1}{8} \sum_{m=1}^{N} \left[
    x_m^u (y_{m+1}^l - y_{m-1}^l) + x_m^l (y_{m+1}^u - y_{m-1}^u)
    \right],
\end{equation}
where $N$ is the number of vertices on each face, and $(x_m, y_m)$ are the in-plane coordinates of the $m$-th vertex $p_m$. Superscripts $u$ and $l$ denote quantities associated with the upper and lower faces, respectively. A key technical step in this algorithm is to establish the vertex connectivity between the two bounding faces. As illustrated in Fig.~\ref{fig:floodalgorithm}(e), two types of connections may occur: (i) one-to-one connections (e.g., $p_1^u \leftrightarrow p_1^l$) and (ii) one-to-two connections (e.g., $p_5^u$ and $p_6^u$ both connect to $p_{5/6}^l$). Once these connections are determined, the edges of the polyhedron are uniquely defined, and the liquid volume within each sub-region can be accurately evaluated using Eqs.~\eqref{eq:cross-section_area} and \eqref{eq:middle_area}.

The third step is to identify the interval that contains the liquid/gas interface $\Gamma_l$. This is achieved by comparing the target liquid volume $c\Delta^3$ with the cumulative sum of the sub-volume integrals $\sum V_i$. The desired interval $M$ satisfies
\begin{equation}
    \sum_{i=1}^{M-1} V_i < c\Delta^3 < \sum_{i=1}^{M} V_i.
\end{equation}
Once $M$ is identified, the remaining liquid volume within this interval is given by
\begin{equation}
    V(\lambda) = c\Delta^3 - \sum_{i=1}^{M-1} V_i,
\end{equation}
where $\lambda$ denotes the unknown height of the interface within interval $M$, as illustrated by the blue region in Fig.~\ref{fig:floodalgorithm}(e). The volume $V(\lambda)$ is expressed as a cubic function of $\lambda$:
\begin{equation}\label{eq:cubic-f}
    V(\lambda) = a_1 \lambda^3 + a_2 \lambda^2 + a_3 \lambda,
\end{equation}
where the coefficients $a_1$, $a_2$, and $a_3$ are determined from the partial volumes evaluated at three reference heights:
\[
\lambda = \frac{1}{3}(D_u - D_l), \quad \frac{2}{3}(D_u - D_l), \quad \frac{3}{3}(D_u - D_l).
\]
The partial volumes evaluated at these reference heights are:
\begin{equation}
\begin{aligned}
    V\left(\tfrac{1}{3}(D_u - D_l)\right) &= \frac{D_u - D_l}{18}
    \left[
    \sum_{m=1}^{N} x_m^l \Delta y_m^l + \sum_{m=1}^{N} x_m^{\frac{1}{3}} \Delta y_m^{\frac{1}{3}} + 
    \frac{1}{2} \sum_{m=1}^{N} \left( x_m^l \Delta y_m^{\frac{1}{3}} + x_m^{\frac{1}{3}} \Delta y_m^l \right)
    \right], \\
    V\left(\tfrac{2}{3}(D_u - D_l)\right) &= \frac{D_u - D_l}{9}
    \left[
    \sum_{m=1}^{N} x_m^l \Delta y_m^l + \sum_{m=1}^{N} x_m^{\frac{2}{3}} \Delta y_m^{\frac{2}{3}} + 
    \frac{1}{2} \sum_{m=1}^{N} \left( x_m^l \Delta y_m^{\frac{2}{3}} + x_m^{\frac{2}{3}} \Delta y_m^l \right)
    \right], \\
    V\left(\tfrac{3}{3}(D_u - D_l)\right) &= \frac{D_u - D_l}{6}
    \left[
    \sum_{m=1}^{N} x_m^l \Delta y_m^l + \sum_{m=1}^{N} x_m^u \Delta y_m^u + 
    \frac{1}{2} \sum_{m=1}^{N} \left( x_m^l \Delta y_m^u + x_m^u \Delta y_m^l \right)
    \right].
\end{aligned}
\end{equation}
Here, superscripts $1/3$ and $2/3$ denote interpolated vertex positions at fractional heights $\lambda$ within the interval, obtained from
\[
x_m^{\lambda} = \lambda x_m^u + (1 - \lambda) x_m^l, \quad y_m^{\lambda} = \lambda y_m^u + (1 - \lambda) y_m^l.
\]
The coefficients in Eq.~\eqref{eq:cubic-f} are then computed as:
\begin{equation}
    \left[
    \begin{array}{c}
        a_1 \\
        a_2 \\
        a_3
    \end{array}
    \right]
    = \frac{D_u - D_l}{12}
    \left[
    \begin{array}{ccc}
        2 & 2 & -2 \\
        -6 & 0 & 3 \\
        6 & 0 & 0
    \end{array}
    \right]
    \left[
    \begin{array}{c}
        \sum x_m^l \Delta y_m^l \\
        \sum x_m^u \Delta y_m^u \\
        \sum (x_m^l \Delta y_m^u + x_m^u \Delta y_m^l)
    \end{array}
    \right],
\end{equation}
with $a_1$, $a_2$, and $a_3$ determined, $\lambda$ is obtained iteratively using the Newton–Raphson method applied to Eq.~\eqref{eq:cubic-f}. The parameter $\gamma$ in Eq.~\eqref{eq:hn}, representing the signed distance from the cell center to the interface, is then updated accordingly.

The complete 3D flood algorithm is summarized in Algorithm~\ref{alg:floodalgorithm}, which enables accurate reconstruction of $\Gamma_l$ in mixed cells. Compared with its 2D counterpart~\cite{huang20252d}, the primary challenges in 3D lie in evaluating sub-interval volumes and solving the nonlinear equation for $\lambda$.

\begin{algorithm}[H]\label{alg:floodalgorithm}
\caption{Reconstruction of the liquid/gas interface in a mixed cell (see Fig.~\ref{fig:floodalgorithm}).}
\begin{enumerate}
    \item Reconstruct the embedded solid surface $\Gamma_s$ and identify the fluid region as a convex polyhedron.
    \item Rotate the polyhedron so that the interface normal $\textbf{n}_l$ aligns with the $z$-axis.
    \item Subdivide the region by inserting horizontal planes at each vertex height; compute the volume of each sub-region.
    \item Locate the interval $M$ that contains $\Gamma_l$ using the cumulative volume test.
    \item Solve the cubic volume function $V(\lambda)$ to obtain the interface height $\lambda$ in interval $M$.
\end{enumerate}
\end{algorithm}

\begin{figure}
    \includegraphics[height=0.36\textwidth]{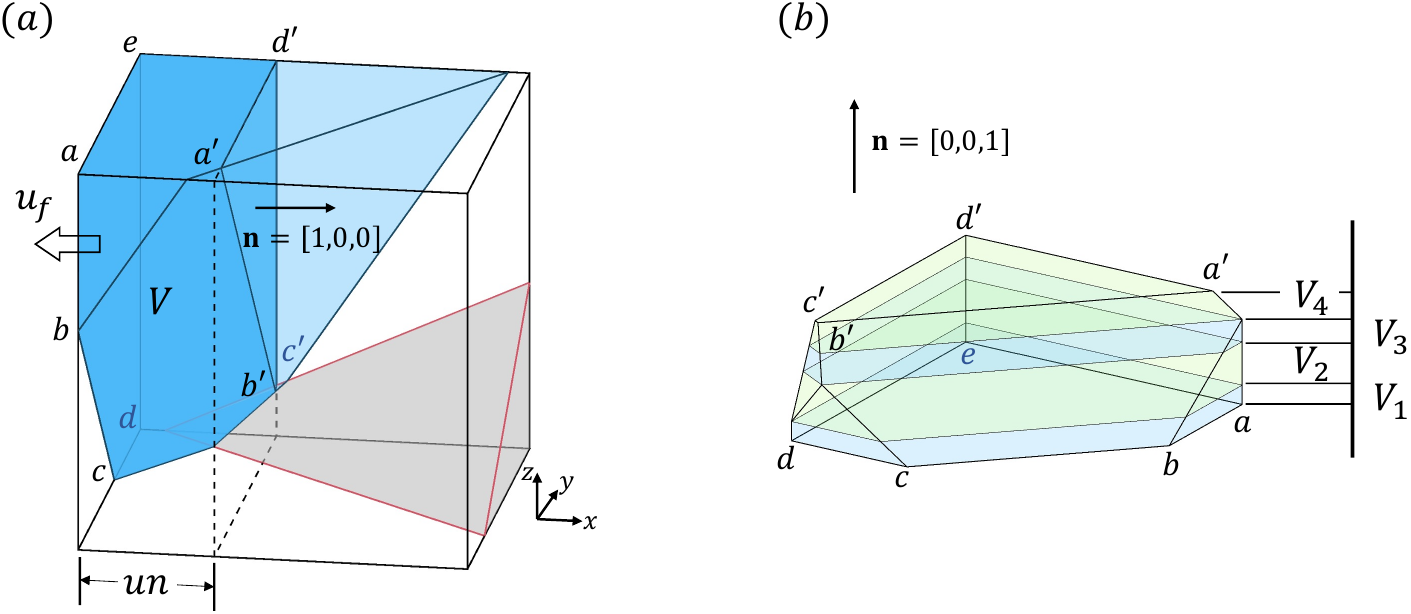}
    \centering
    \caption{
    Illustration of the advected liquid volume $V$ within a mixed cell. The advection interval has a width of $un = u_f \cdot dt$.
    }
    \label{fig:advectionvolume}
\end{figure}

A second challenge is the evaluation of the advected liquid volume $V$ in a mixed cell over a time step $dt$. As illustrated in Fig.~\ref{fig:advectionvolume}, this volume corresponds to the dark blue region transferred across face \textit{abcde} into a neighboring cell. The region is bounded by the cell faces, the sharp interfaces $\Gamma_s$ and $\Gamma_l$, and the advected surface \textit{a'b'c'd'}, which is displaced a distance $un = u_f \cdot dt$ from \textit{abcde}. Estimating this volume is effectively the inverse of the interface reconstruction problem: here, $\textbf{n}_l$ and $\gamma$ are known, while the transported volume $c = V/\Delta^3$ is unknown. The procedure follows the same principles as the flood algorithm:
\begin{enumerate}
    \item Rotate the cell so that the normal of face \textit{abcde} aligns with the $z$-axis.
    \item Subdivide the advected region using horizontal planes through each vertex.
    \item Evaluate each sub-volume using Eq.~\eqref{eq:cross-section_area}.
    \item Sum all contributions to obtain the total transported volume $V = \sum_i V_i$.
\end{enumerate}

This extended flood algorithm therefore enables both accurate planar reconstruction of $\Gamma_l$ within mixed cells and robust estimation of the liquid volume fluxed across cell faces.

\subsection{Conservative scheme for geometric VOF advection in mixed cells}\label{sec:conservationadvectionscheme}

\begin{figure}
    \includegraphics[height=0.36\textwidth]{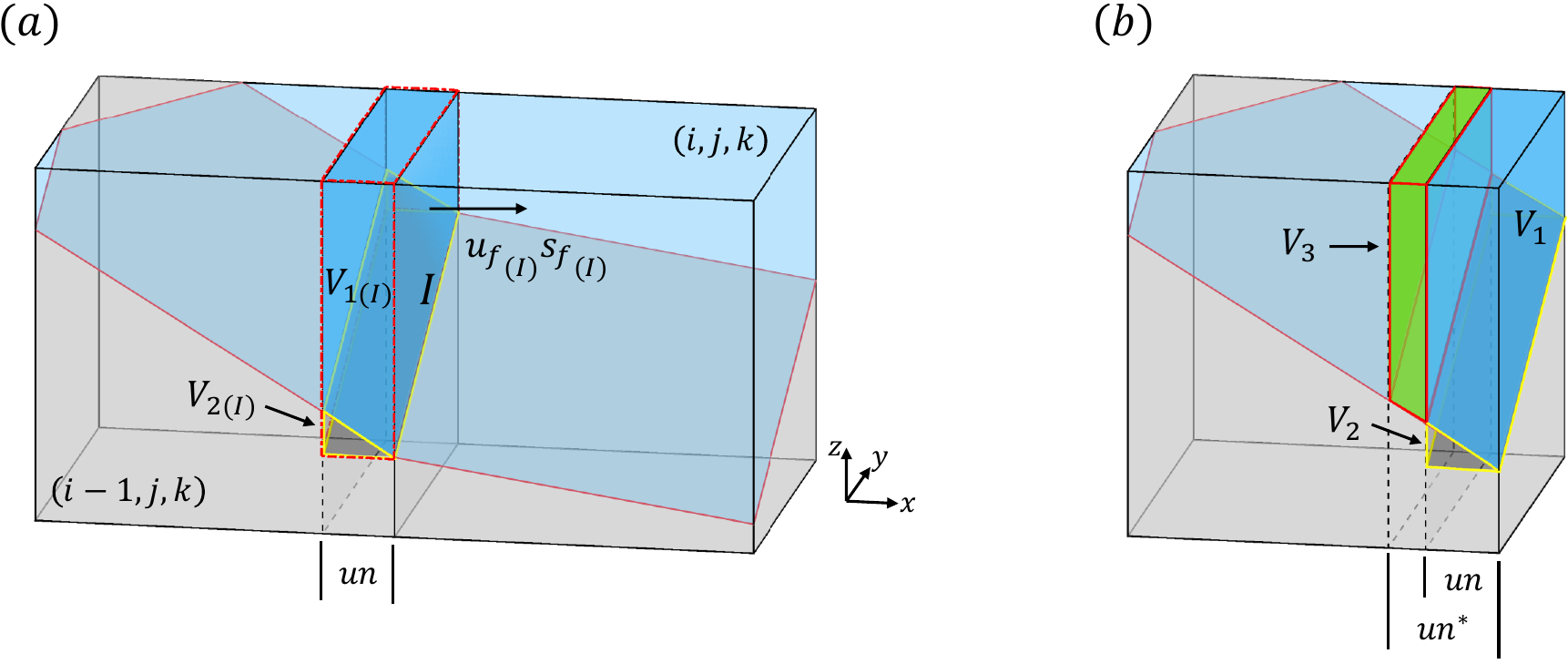}
    \centering
    \caption{
    Sketch of the geometric VOF advection scheme in mixed cells. (a) $s_f$ denotes the fluid fraction of the face, and $\Delta$ is the mesh size. $I$ indicates the face intersected by the embedded solid boundary. $V_1$ and $V_2$ represent the volumes of the dark blue and dark grey regions, respectively. (b) ${un}^*$ denotes the corrected width of the advection interval, and $V_3$ is the green region, numerically equal to $V_2$.
    }
    \label{fig:advectionmethod}
\end{figure}

As demonstrated in our companion 2D study~\cite{huang20252d}, the geometric advection scheme implemented in \textit{Basilisk}, i.e., Eq.~\eqref{eq:advection discretization scheme1} for regular cells, cannot be directly applied to mixed cells. The difficulty arises from the presence of $\Gamma_s$ within the advection interval: naively applying Eq.~\eqref{eq:advection discretization scheme1} results in systematic over- or under-estimation of fluxes, thereby violating local mass conservation. This issue persists in 3D. The 2D remedy—modifying the advection velocity to account for the volume displaced by $\Gamma_s$—remains effective in 3D, as detailed below.

As described in \S~\ref{sec:regularVOFgeometric}, geometric advection is performed in a directionally split manner along the $x$-, $y$-, and $z$-directions. Consider advection along the $x$-direction. For mixed cells, Eq.~\eqref{eq:advection discretization scheme1} is modified as follows:
\begin{equation}\label{eq:exampleoferrorinadvection}
    \frac{c^* - c^{n - \frac{1}{2}}}{dt} \Delta^3 = 
    \left(\widetilde{F}_{i - \frac{1}{2}}^{n - \frac{1}{2}} - \widetilde{F}_{i + \frac{1}{2}}^{n - \frac{1}{2}}\right)
    - c_c \left(u_{i - \frac{1}{2}}^{n - \frac{1}{2}} s_{i - \frac{1}{2}} - u_{i + \frac{1}{2}}^{n - \frac{1}{2}} s_{i + \frac{1}{2}}\right),
\end{equation}
where $s_{i \pm 1/2}$ denotes the fluid fraction on the corresponding face, and $\widetilde{F}_{i \pm 1/2}$ represents the corrected flux. Compared with Eq.~\eqref{eq:advection discretization scheme1}, three modifications are introduced: (i) the explicit inclusion of the face area fraction $s_f$, (ii) a corrected evaluation of fluxes $\widetilde{F}$ to ensure exact mass conservation, and (iii) according to Ref.\cite{huang20252d}, $c_c$ is defined as 
\begin{equation}\label{eq:ccmixedcell}
    c_{c}=\left\{
    \begin{array}{rcl}
        1,       &      & c>c_s/2\\
        0,       &      & c\leq c_s/2\\
    \end{array}\right..
\end{equation}.

To illustrate this, consider a liquid-filled mixed cell $(i,j,k)$ at time $t^{n - 1/2}$, surrounded by similarly filled neighbors (see Fig.~\ref{fig:advectionmethod}(a)). In this scenario, the flux depends solely on the geometry and velocity, as $c$ is constant. Since there is no net liquid inflow or outflow, $c$ should remain unchanged after advection:
\[
c^{n + 1/2} = c^{n - 1/2}.
\]
This condition enforces strict volume conservation: the total flux into the cell must be zero. Let us consider the left face $I = (i - 1/2, j, k)$, which is intersected by $\Gamma_s$. According to Eq.~\eqref{eq:cf-normal}, the advected liquid volume into cell $(i,j,k)$ through face $I$ is
\begin{equation}
    V_{(I)} = V_{1(I)} = \left( u_f dt \cdot s_f \cdot \Delta^2 \right)_I - V_{2(I)},
\end{equation}
where the first term represents the intended advection volume (red region in Fig.~\ref{fig:advectionmethod}(a)), and $V_{2(I)}$ denotes the volume blocked by the solid (gray region). Summing contributions from all $N$ faces of the cell, the net volume change becomes:
\begin{equation}\label{eq:errorinconservation}
    \left(c^{n + \frac{1}{2}} - c^{n - \frac{1}{2}}\right)\Delta^3 = (1 - c_c) \nabla \cdot (u_f s_f) \cdot dt \Delta^3 + \sum_{i=1}^{N} V_{2,i}.
\end{equation}
For incompressible flow, $\nabla \cdot (u_f s_f) = 0$. 
Hence, in a liquid-filled cell ($c_c = 1$), Eq.~\eqref{eq:errorinconservation} simplifies to:
\begin{equation}\label{eq:clocalerror}
    \left(c^{n + \frac{1}{2}} - c^{n - \frac{1}{2}}\right)\Delta^3 = \sum_{i=1}^{N} V_{2,i}.
\end{equation}
This residual arises purely from neglecting the geometry of $\Gamma_s$, and leads to local errors in $c$. If $\sum V_{2} > 0$, the cell becomes overfilled; if $\sum V_{2} < 0$, it becomes underfilled. To address this issue, we redefine the flux through face $I$ using a corrected advection width ${un}^*$ that compensates for the missing volume. As illustrated in Fig.~\ref{fig:advectionmethod}(b), a virtual extension is introduced to recover a volume $V_3$ satisfying $V_3 = V_2$. This is equivalent to reconstructing a fictitious liquid/gas interface with normal vector $\mathbf{n}'_l = [-1, 0, 0]$. Using the flood algorithm (Algorithm~\ref{alg:floodalgorithm}) with this normal and a specified volume fraction $c = u_f dt \cdot s_f \Delta^2/\Delta^3$, the distance $\gamma$ from the cell centroid to the interface is obtained:
\begin{equation}\label{eq:omega2}
    \gamma = \text{Algorithm~1}\left(c = \frac{u_f dt \cdot s_f \Delta^2}{\Delta^3}, \, \mathbf{n}'_l = [-1, 0, 0]\right).
\end{equation}
The corrected advection width is then given by
\begin{equation}
    {un}^* = \left( \frac{1}{2} + \gamma \right)\Delta.
\end{equation}
With this correction, the flux through face $I$ becomes:
\[
V = V_1 + V_3 = u_f dt \cdot s_f \cdot \Delta^2,
\]
recovering the full advection volume. Applying this correction across all faces and substituting back into Eq.~\eqref{eq:errorinconservation} yields:
\begin{equation}
    \left(c^{n + \frac{1}{2}} - c^{n - \frac{1}{2}}\right)\Delta^3 = (1 - c_c) \nabla \cdot (u_f s_f) \cdot dt \Delta^3.
\end{equation}
In a liquid-filled cell ($c_c = 1$), this confirms strict local conservation. Note, however, that this analysis assumes the absence of a gas phase. In more general three-phase cells containing both $\Gamma_s$ and $\Gamma_l$, the gas-phase contribution must also be accounted for, yielding
\[
\left(c^{n + \frac{1}{2}} - c^{n - \frac{1}{2}}\right)\Delta^3 = \sum (V_{2} + V_{\mathrm{gas}}),
\]
where $V_{\mathrm{gas}}$ denotes the gas-phase advection through face $I$. Since $c^{n - 1/2}$ and $V_{\mathrm{gas}}$ are independent, minor volume errors for each phase may still accumulate. Finally, to prevent over- or under-filling in mixed cells, the CFL condition must be satisfied:
\[
\frac{u_f s_f dt}{c_s \Delta} \leq \mathrm{CFL}.
\]
In 3D, however, cells intersected by $\Gamma_s$ may have very small $c_s$, resulting in $s_f / c_s \gg 1$ and thus excessively restrictive time steps. In the following section, we introduce a redistribution scheme that alleviates this constraint and enhances both robustness and efficiency in complex geometries.

\subsection{Small-cut-cell restriction}\label{sec:smallcutcell}

\begin{figure}
    \includegraphics[height=0.4\textwidth]{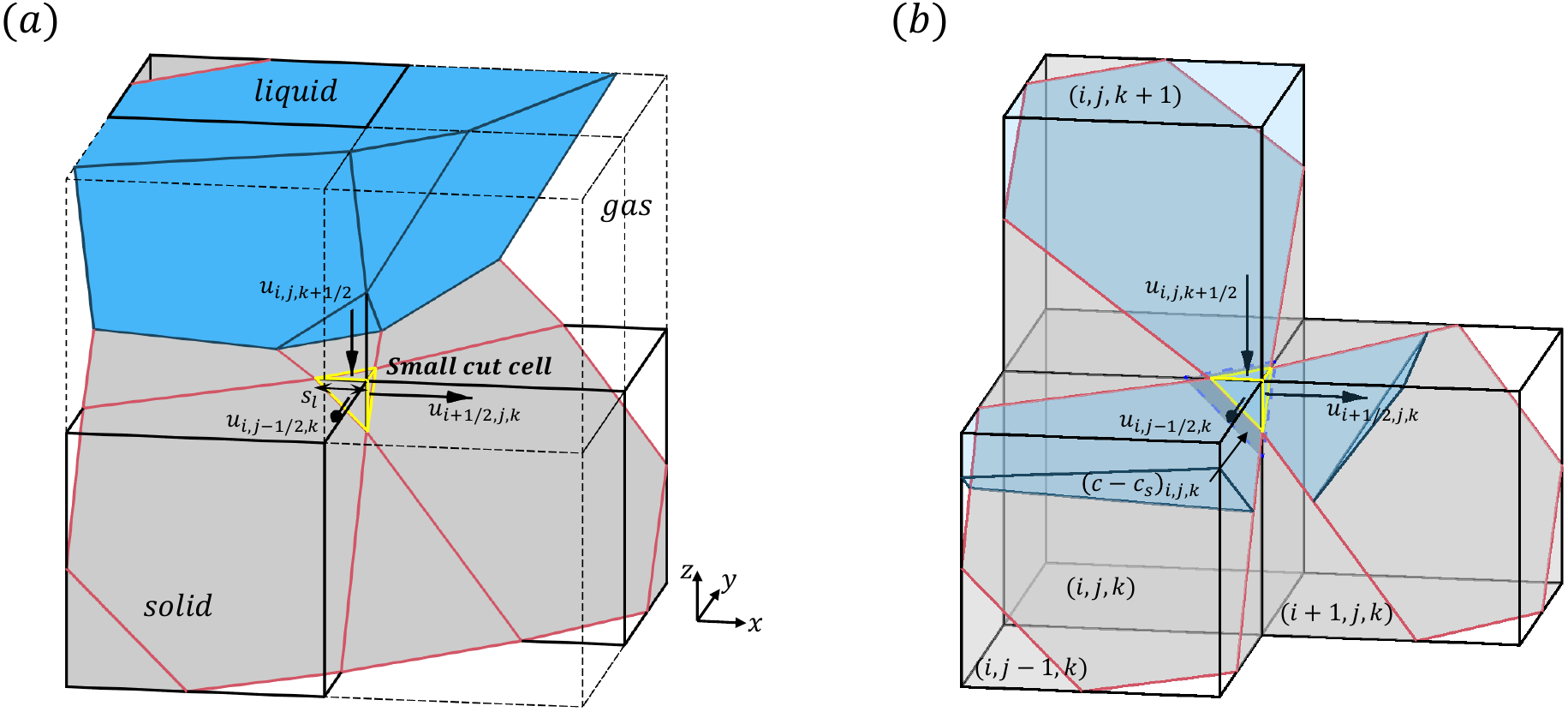}
    \centering
    \caption{
    Sketch of the redistribution advection method in mixed cells. The yellow solid cell represents a small cut cell. Velocities $u_{i+1/2,j,k}$, $u_{i,j-1/2,k}$, and $u_{i,j,k+1/2}$ are shown on their corresponding faces. (a) $s_l$ denotes the fluid line fraction along a grid edge. (b) The dark blue region within the small cut cell represents the overfilled volume $(c - c_s)_{i,j,k}$.
    }
    \label{fig:smallcutcell}
\end{figure}

It is well established that CFL constraints become particularly restrictive in Cartesian-grid methods with complex solid boundaries~\cite{kleefsman2005volume,colella2006cartesian,schneiders2016efficient}, where the presence of small cut cells with $s_f / c_s \gg 1$ is unavoidable. In our previous 2D study~\cite{huang20252d}, such cells were forcibly removed by discarding those with fluid fractions below a threshold ($c_s < 10^{-2}$), which slightly altered the solid boundary during initialization. However, this approach is no longer feasible in 3D. Removing cells with $c_s < 10^{-3}$ would significantly distort the geometry of embedded boundaries due to the stronger volumetric dependence (see \S~\ref{sec:resultsVOFadvection1}). As a result, 3D grids contain substantially more small cut cells than in 2D, making this removal strategy impractical.

Alternative strategies, such as cell merging~\cite{chung2006cartesian} and momentum redistribution~\cite{ghigo2021conservative}, have been employed to treat small cut cells in momentum equations. However, these methods have rarely been extended to volume-fraction advection. Accurate interface reconstruction in the presence of complex boundaries introduces additional challenges: cell merging complicates the geometric definition of the resulting fluid region, limiting the applicability of the ``flood algorithm'' described in \S~\ref{sec:floodalgorithm}. Motivated by the redistribution concept in~\cite{ghigo2021conservative}, we propose a new two-step advection scheme that redistributes fluxes in interfacial cells. This approach relaxes the timestep restrictions imposed by small cut cells while avoiding unnecessary advection in fully filled or empty cells. The scheme proceeds in two steps. First, the volume fraction is advected using a variant of the conservative scheme from \S~\ref{sec:conservationadvectionscheme}, with the timestep no longer constrained by $c_s$. Second, any excess or deficit volume in interfacial cells is redistributed to downstream neighbors.

In the first step, the dependence of the time step on $c_s$ and $s_f$ is eliminated. The updated CFL constraint becomes $dt < \Delta / u_f \cdot \text{CFL}$, making it independent of $s_f / c_s$. However, this relaxation of $dt$ may introduce inconsistencies in $c_c$, as defined by Weymouth et al.~\cite{weymouth2010conservative}, potentially leading to over-filled or under-filled cells. This issue is addressed in the subsequent step. In the present step, the advection method described in \S~\ref{sec:conservationadvectionscheme} is accordingly modified as follows:

\begin{enumerate}
    \item \textbf{Gas-filled cell} ($c = 0$): No correction is applied; the advected liquid volume is set to $V = 0$.
    \item \textbf{Liquid-filled cell} ($c = c_s$): No correction is applied; the advected liquid volume is $V = u_f \, dt \cdot s_f \Delta^2$.
    \item \textbf{Interfacial cell} ($0 < c < c_s$) with $c_s > u_f s_f \, dt \Delta^2 / \Delta^3$: the modified advection velocity from \S~\ref{sec:conservationadvectionscheme} is used, and the corrected width ${un}^*$ is computed. If the maximum line fraction ${s_l}_{\max}$ along the advection direction satisfies ${s_l}_{\max} > {un}^*$, the flood algorithm is applied to compute the advected volume geometrically. Otherwise, the entire fluid is advected as $V = c \Delta^3$.
    \item \textbf{Interfacial cell} ($0 < c < c_s$) with $c_s < u_f s_f \, dt \Delta^2 / \Delta^3$: because Eq.~\eqref{eq:omega2} is invalid, no velocity correction is applied. The original advection width $un$ is compared with ${s_l}_{\max}$, and the advected volume is computed following the same procedure as in case (3).
\end{enumerate}
Here, ${s_l}_{\max}$ denotes the maximum line fraction among all edges in the advection direction (see Fig.~\ref{fig:smallcutcell}(a)).

For full or empty cells away from the interface, this scheme preserves local conservation, as $c_c$ remains unchanged. In interfacial cells, however, volume errors may occur, resulting in over-filled ($c / c_s > 1$) or over-emptied ($c / c_s < 0$) cells. These inconsistencies arise from several factors: (i) the failure of interval correction when $c_s\Delta^3$ is smaller than the advected volume $u_f s_f dt \Delta^2$; (ii) minor divergence in the velocity field due to the approximate projection method; and (iii) the relaxation of the timestep, which undermines strict conservation of $c_c$ (Eq.~\eqref{eq:ccmixedcell}). To restore conservation, $c_c$ is redefined in its original form, $c / c_s$, and any excess or deficit volume is redistributed to downstream three-phase neighbors along the flow direction. The redistribution weights are determined by the local occupancy: for over-emptied cells ($c / c_s < 0$), the weights are $w = c / c_s$, and for over-filled cells ($c / c_s > 1$), the weights are $w = 1 - c / c_s$.

Figure~\ref{fig:smallcutcell} illustrates this redistribution. In panel (a), cell $(i,j,k)$ is a small cut cell prior to advection, with outgoing velocities indicated. In panel (b), after advection, the cell becomes over-filled ($c > c_s$), and the excess volume $(c - c_s)_{i,j,k}$ is shown in dark blue. This excess is then redistributed to downstream neighbors, for example,

\begin{equation}
(c - c_s)_{i,j,k} \cdot \frac{w_{i+1,j,k}}{w_{i+1,j,k} + w_{i,j-1,k}} \quad \text{to} \quad (i+1,j,k),
\end{equation}

\begin{equation}
(c - c_s)_{i,j,k} \cdot \frac{w_{i,j-1,k}}{w_{i+1,j,k} + w_{i,j-1,k}} \quad \text{to} \quad (i,j-1,k),
\end{equation}

where weights are defined as $w_{i+1,j,k} = 1 - c_{i+1,j,k} / c_{s,i+1,j,k}$ and $w_{i,j-1,k} = 1 - c_{i,j-1,k} / c_{s,i,j-1,k}$ for the over-filled case. In over-emptied scenarios, the corresponding weights are given by $w = c / c_s$. This redistribution strategy eliminates the timestep restriction imposed by small cut cells, while maintaining geometric accuracy and ensuring volume conservation.

\section{Contact angle condition}\label{sec:contactanglecondition}

In this section, we discuss the implementation of contact angle conditions in mixed cells, with a particular focus on cells containing all three phases: liquid, gas, and solid.

\subsection{Linear fitting --- for surface tension calculation}\label{sec:linearfitting}

\begin{figure}
    \includegraphics[height=0.36\textwidth]{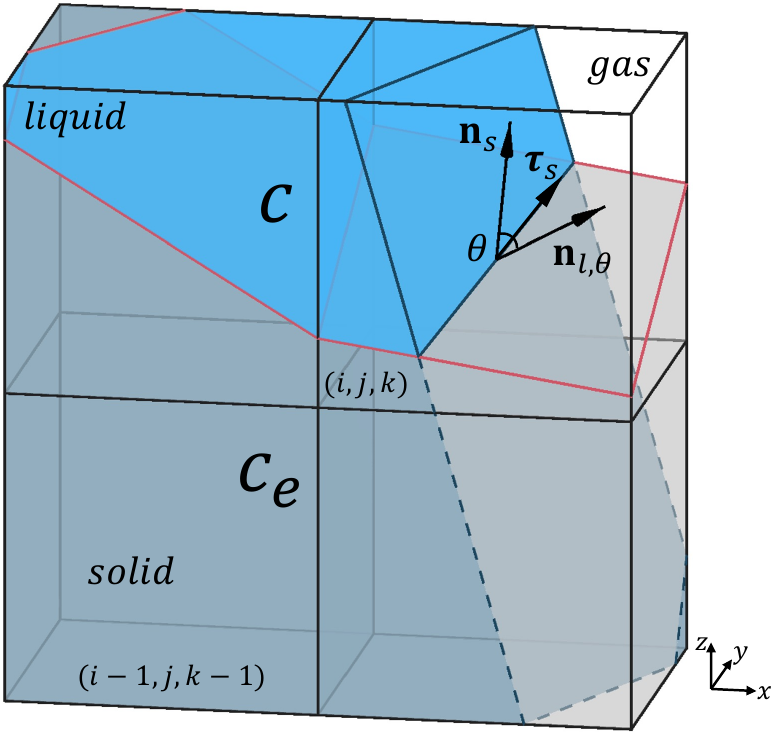}
    \centering
    \caption{
    Sketch illustrating the extension of the liquid/gas interface into the solid region. $\boldsymbol{\tau}_{s}$ denotes the unit tangential vector along the contact line. The dark blue region represents the ghost fluid volume, denoted by $c_e$.
    }
    \label{fig:extension}
\end{figure}

We begin by introducing the volume fraction extension method employed for computing surface tension forces. As discussed previously, an accurate evaluation of the volume fraction gradient $\nabla c$ is crucial for computing the surface tension term $\sigma \kappa \nabla c$ in the momentum equation. To achieve this near the contact line, a virtual volume fraction field must be constructed within the solid region to properly enforce the prescribed contact angle. This approach, originally proposed by Tavares et al.~\cite{tavares2024coupled}, is based on a linear extension of the interface, either as a straight line in 2D or as a planar surface in 3D. We refer to this method as the linear fitting scheme, and its 3D implementation is outlined below.

To extend the liquid/gas interface into the solid domain as a planar surface, we first define the interface normal vector in contact line cells. Denote this normal vector by $\textbf{n}_{l,\theta}$, where the subscript indicates its dependence on the contact angle $\theta$. The solid boundary normal vector $\textbf{n}_s$ is obtained directly from the geometry. The vector $\textbf{n}_{l,\theta}$ is then computed by rotating $\textbf{n}_s$ about the contact line direction $\boldsymbol{\tau}_s$ by the angle $\theta$, as illustrated in Fig.~\ref{fig:extension}. The rotation follows the right-hand rule and is expressed as:
\begin{equation}\label{eq:n-l-theta}
    \textbf{n}_{l,\theta} = \cos{\theta}\cdot \textbf{n}_s + \sin{\theta} \cdot (\boldsymbol{\tau}_s \times \textbf{n}_s),
\end{equation}

The contact line direction $\boldsymbol{\tau}_s$ is not known a priori and must be approximated. In this work, we employ the Mixed-Young-Centered (MYC) scheme~\cite{2007Interface} to estimate $\textbf{n}_{l,\mathrm{myc}}$, an approximate interface normal that avoids using ghost values of $c$ within the solid. The contact line direction is then computed via the cross product:
\[
\boldsymbol{\tau}_s = \textbf{n}_s \times \textbf{n}_{l,\mathrm{myc}}.
\]
Substituting this into Eq.~\eqref{eq:n-l-theta} yields the desired normal vector $\textbf{n}_{l,\theta}$, which is used to construct the extended volume fraction field, thereby ensuring that the contact angle condition is correctly reflected by the volume fraction gradient. Next, the ``flood algorithm'' described in \S~\ref{sec:floodalgorithm} is employed to reconstruct the 3D planar interface $\Gamma_l$ within the contact line cell, using the prescribed normal vector $\textbf{n}_{l,\theta}$. This reconstructed interface is then geometrically extended into the adjacent solid region to define ghost volume fractions, denoted by $c_e$, for solid cells.

The ghost value $c_e$ is assigned under the assumption that the portion of the solid lying below the extended interface belongs to the liquid phase. For instance, a solid cell such as $(i-1, j, k-1)$, initially assigned $c = \mathrm{NULL}$, receives $c_e = 1$ following this extension. These ghost values enable the construction of a complete $3 \times 3 \times 3$ stencil centered at the contact line cell $(i,j,k)$, which is necessary to compute $\nabla c$. In 3D configurations, a single solid cell may be adjacent to multiple contact line cells. To account for this, a weighted average is used to determine $c_e$, incorporating contributions from all neighboring contact line cells. The weights take the form $c_s (1 - c_s) c (1 - c)$, as proposed by Tavares et al.~\cite{tavares2024coupled}, favoring cells that are strongly interfacial.

With the extended ghost values $c_e$ available, the surface tension force $\sigma \kappa \delta_s \textbf{n}_l = \sigma \kappa \nabla c$ in Eq.~\eqref{eq:ns} can be computed consistently, even near embedded solid boundaries. This preserves the balanced-force condition between surface tension and pressure in the momentum equation. However, as emphasized in our previous 2D study~\cite{huang20252d}, the ghost values $c_e$ obtained from the linear extension should \emph{not} be used directly to enforce the contact angle condition. The linear extension captures only the local interface geometry from the contact line cell (or a nearby layer of mixed cells) and does not reflect the broader interface curvature required for accurate contact angle enforcement. In particular, precise estimation of the interface normal $\textbf{n}_l$ at the contact line depends on both the local and surrounding interface geometry, including contributions from interfacial cells farther from the solid surface. Therefore, enforcing the contact angle condition necessitates a more robust reconstruction method that integrates interface shape over a larger neighborhood.

In the 3D study by Tavares et al.~\cite{tavares2024coupled}, the original HF method in \textit{Basilisk} was employed, in which $\Gamma_l$ is linearly extended into the solid and ghost volume fractions are summed to define HFs over a single ghost layer to evaluate curvature and normal vectors. While straightforward to implement, our numerical tests (see \S~\ref{sec:surfacetensiondrivendropletspreading}) indicate that this approach exhibits noticeable errors in enforcing the contact angle condition.

These observations strongly motivate the development of a more rigorous HF-based method tailored for three-phase contact lines on complex 3D geometries. Although the need for such a method was established in our prior 2D study~\cite{huang20252d}, its extension to 3D introduces additional challenges in algorithmic implementation and geometric robustness. The details of the improved HF method are presented in the following subsections.

\subsection{Paraboloid fitting --- for the contact angle condition}\label{sec:paraboloidfittinghorizontal}

\begin{figure}
    \includegraphics[height=0.36\textwidth]{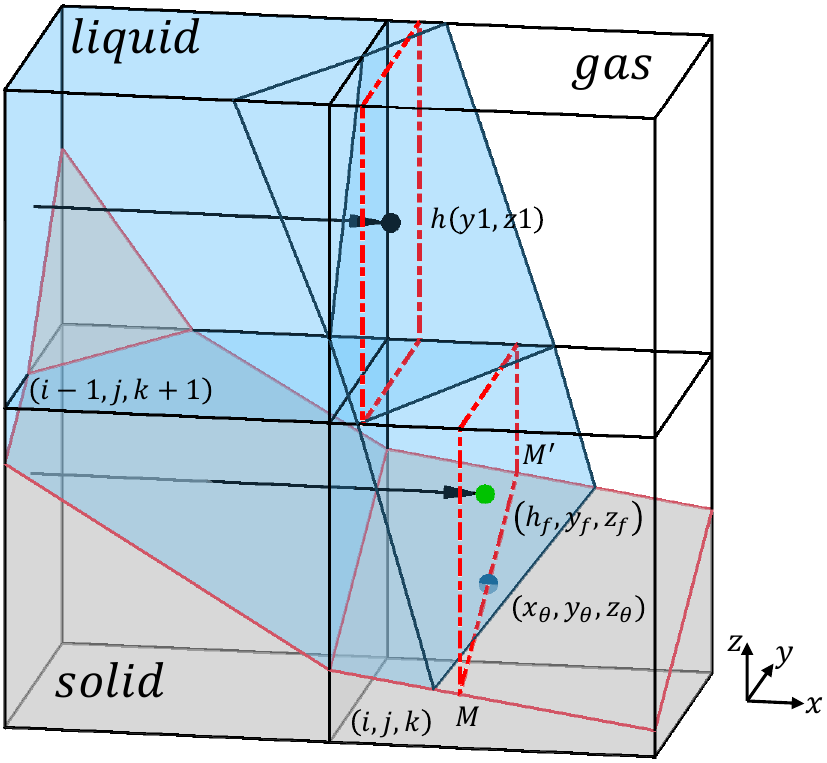}
    \centering
    \caption{
    Schematic diagram illustrating the definition of the height function within the fluid region. The red dashed plane represents the vertical interface used to determine the height function $(h_f, y_f, z_f)$ in the fluid, and the midpoint $(x_\theta, y_\theta, z_\theta)$ of segment $MM'$, where the contact angle condition is applied.
    }
    \label{fig:paraboloidHF}
\end{figure}

\begin{figure}
    \includegraphics[height=0.36\textwidth]{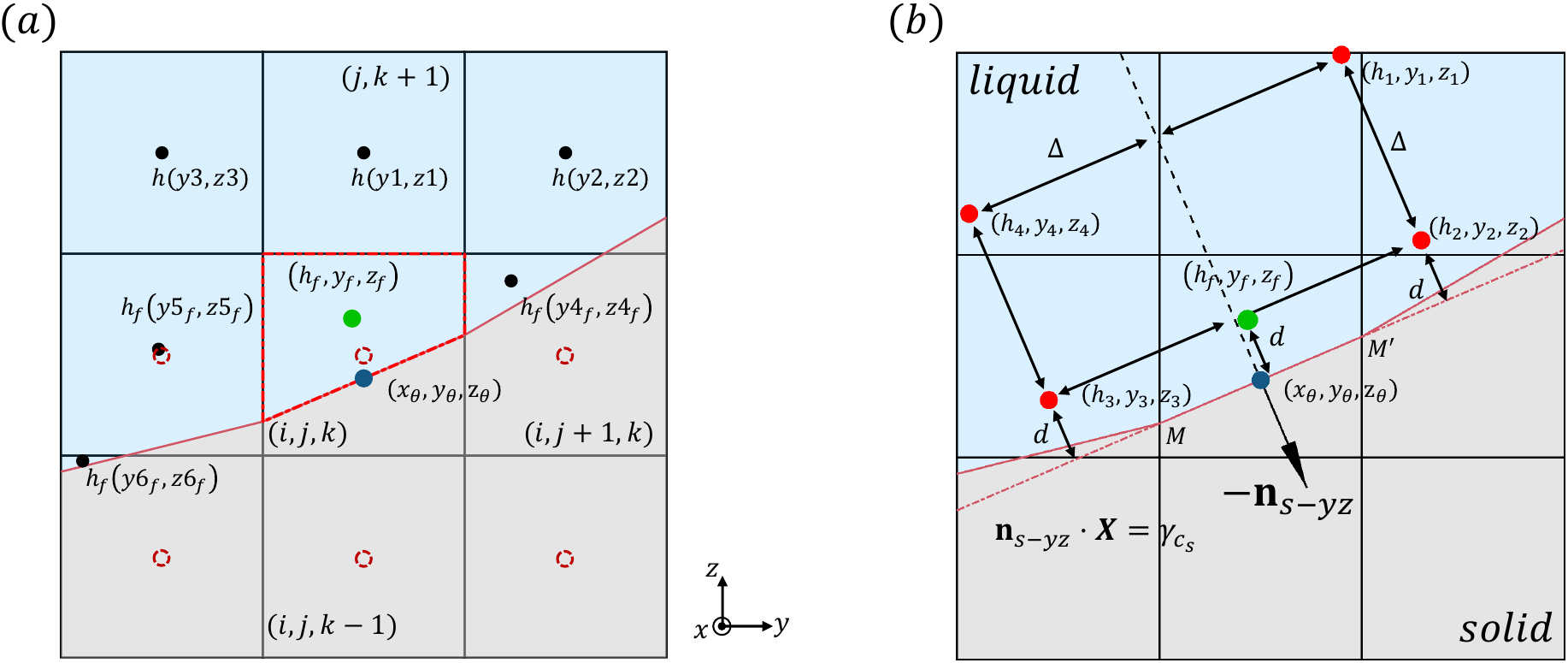}
    \centering
    \caption{
    (a) Schematic of the pre-fitting process used to identify the four red points (highlighted in panel b). Black and green dots represent height function points within the fluid, while the blue dot indicates the contact line location. Dashed red circles denote cell centers. (b) Final stencil showing the red points used to fit the paraboloid that incorporates the contact angle. The vector $\textbf{n}_{s\text{-}yz}$ denotes the projection of the solid normal $\textbf{n}_s$ onto the $yz$-plane. The dashed red line represents the 2D interface reconstruction, passing through $(x_\theta, y_\theta, z_\theta)$ and orthogonal to $\textbf{n}_{s\text{-}yz}$. The distance $d$ is the perpendicular gap between the height function point $(h_f, y_f, z_f)$ and the reconstructed interface.
    }
    \label{fig:paraboloidfit}
\end{figure}

In the presence of complex solid boundaries, standard interpolation methods (e.g., Eq.~\eqref{eq:standard-interpolation}), which are commonly used to enforce contact angles on flat surfaces, become inadequate. As demonstrated in our previous 2D study~\cite{huang20252d}, these methods can produce substantial errors when the solid surface is curved or inclined.

To address this limitation, we introduce a paraboloid fitting approach that incorporates the contact angle condition into the interface reconstruction near embedded boundaries.

\begin{algorithm}[H]
\caption{Definition of the $x$-direction height function in a mixed cell (see Fig.~\ref{fig:paraboloidHF}).}
\label{alg:hfdefinition}
\begin{enumerate}
    \item Sum the volume fractions along each $(j,k)$ row in the $x$-direction.
    \item Identify the interfacial cell in the row, i.e., the cell $(i,j,k)$ containing the interface.
    \item Reconstruct a planar interface perpendicular to the $x$-axis (red dashed plane).
    \item Compute the centroid of the interface surface, denoted by $(h_f, y_f, z_f)$, which defines the fluid height function. Also determine the midpoint of the intersection segment $MM'$, denoted by $(x_\theta, y_\theta, z_\theta)$, where the contact angle condition is enforced.
\end{enumerate}
\end{algorithm}

The first step of the method is to define the height function (HF) in each mixed cell near the solid boundary, which serves as a prerequisite for the subsequent paraboloid fitting. As illustrated in Fig.~\ref{fig:paraboloidHF}, we consider the $x$-direction as an example. The HF is constructed by summing the volume fractions along the $x$-axis for each $(j,k)$ row to determine the location of the liquid/gas interface $\Gamma_l$. Once identified, a vertical planar interface perpendicular to the $x$-direction is reconstructed using the ``flood algorithm'' introduced in \S~\ref{sec:floodalgorithm}. The reconstruction distinguishes between two cases:

\begin{itemize}
    \item \textbf{Case 1: Interface does not intersect the solid boundary} (e.g., row $k+1$). The reconstructed interface lies entirely within the fluid. The HF point is taken as the centroid of the reconstructed planar surface, denoted $h(y1, z1)$.
    
    \item \textbf{Case 2: Interface intersects the solid boundary} (e.g., row $k$). In this case, the interface cuts through the solid, forming an intersection segment $MM'$. The HF point is defined as the centroid of the reconstructed interface, $(h_f, y_f, z_f)$, while the midpoint of $MM'$—$(x_\theta, y_\theta, z_\theta)$—is identified as the precise location where the contact angle condition is enforced.
\end{itemize}

These procedures are summarized in Algorithm~\ref{alg:hfdefinition}. Once the HF points for cell $(i,j,k)$ and its neighbors are collected, a paraboloid surface is fitted to these points to approximate the interface geometry near the contact line. The fitted surface is expressed as

\begin{equation}\label{eq:targetfittingparaboloid}
    h(y,z) = b_1 y^2 + b_2 z^2 + b_3 yz + b_4 y + b_5 z + b_6,
\end{equation}
where the coefficients $b_1$ through $b_6$ are determined from the stencil of HF data.

This formulation follows the method proposed by Han et al.~\cite{han2021consistent}, which has demonstrated first-order convergence in curvature estimation for regular structured domains. Here, we extend it to complex solid geometries by incorporating a pre-fitting procedure that selects appropriate HF points from the fluid region and enforces the contact angle constraint via the point $(x_\theta, y_\theta, z_\theta)$.

In the present context, at least six independent data points are required to solve for the six coefficients in Eq.~\eqref{eq:targetfittingparaboloid}. These include one point from the fluid region $(h_f, y_f, z_f)$, one constraint derived from the contact angle $\theta$ at the contact point $(x_\theta, y_\theta, z_\theta)$, and at least four additional HF points from neighboring interfacial cells. However, identifying four suitable and geometrically consistent HF points is particularly challenging near complex solid boundaries. Our numerical experiments indicate that poor selection of neighboring points can lead to unstable or inaccurate interface reconstruction, thereby compromising the effectiveness of the height function method in solid-affected regions. To address this issue, we adopt a pre-fitting strategy that filters candidate HF points prior to the final paraboloid fitting. This strategy ensures that the selected points are geometrically consistent, lie within the fluid domain, and are sufficiently well-distributed to enable stable curvature reconstruction near embedded solids.

We emphasize that, for all cells adjacent to solid boundaries, the HF is evaluated either at the cell center (e.g., cell $(i,j,k+1)$) or at the partial fluid centroid (e.g., cell $(i,j,k)$), following the procedure illustrated in Fig.~\ref{fig:paraboloidHF}. Specifically, for a contact line cell such as $(i,j,k)$ in Fig.~\ref{fig:paraboloidfit}(a), all available HF values within a $3 \times 3$ or $5 \times 5$ stencil in the $(y,z)$ plane are collected. These points include values at full cell centers, e.g., $h(y1, z1)$, $h(y2, z2)$, $h(y3, z3)$, as well as partial fluid centroids, e.g., $h_f({y4}_{f}, {z4}_{f})$, $h_f({y5}_{f}, {z5}_{f})$, $h_f({y6}_{f}, {z6}_{f})$.

A critical constraint for selecting valid interpolation points is that all chosen points must lie on the same side of the solid boundary as the primary HF value $(y_f, z_f)$. This condition is enforced using a directional dot product check:
\begin{equation}
    \big[(y_f, z_f)-(y_\theta, z_\theta)\big] \cdot \big[(y_i, z_i)-(y_\theta, z_\theta)\big] > 0,
\end{equation}
which ensures that any candidate point $(y_i, z_i)$ lies on the same side of the projected solid surface as $(y_f, z_f)$ relative to the contact point $(y_\theta, z_\theta)$. This criterion guarantees geometric consistency and prevents mixing data across the solid boundary.

As the first step in the paraboloid fitting strategy, all valid HF data points satisfying the above constraint are used to fit a preliminary paraboloid surface of the form
\begin{equation}\label{eq:prefittingparaboloid}
    h_p(y, z) = p_1 y^2 + p_2 z^2 + p_3 y z + p_4 y + p_5 z + p_6,
\end{equation}
where $p_1, \dots, p_6$ are the coefficients of the pre-fitted surface $h_p(y,z)$. Since more than six HF values are typically available, thus unlike the accurate fitting of Eq.~\eqref{eq:targetfittingparaboloid}, Eq.~\eqref{eq:prefittingparaboloid} is solved using a least-squares approximation.
Importantly, this pre-fitting step employs only HF values within the fluid region (e.g., the green point in Fig.~\ref{fig:paraboloidfit}(a)) and does \emph{not} incorporate contact angle information (the blue point at $(x_\theta, y_\theta, z_\theta)$). Its purpose is to construct a stable geometric approximation of the local interface profile based solely on available volume fraction data within the admissible fluid domain.
Following the pre-fitting, four auxiliary points
\[
(h_1, y_1, z_1),\quad (h_2, y_2, z_2),\quad (h_3, y_3, z_3),\quad (h_4, y_4, z_4)
\]
are selected to complete the final paraboloid fitting required for enforcing the contact angle, as illustrated in Fig.~\ref{fig:paraboloidfit}(b). These points are chosen based on both the geometry of the pre-fitted paraboloid (Eq.~\eqref{eq:prefittingparaboloid}) and the known contact point $(x_\theta, y_\theta, z_\theta)$. The procedure consists of the following steps:

\begin{enumerate}
    \item Identify the interface intersection segment $MM'$ (red dash-dotted line in Fig.~\ref{fig:paraboloidfit}) and construct a perpendicular reference line (black dash-dotted) through the contact point in the $yz$-plane.
    \item Position the four interpolation points symmetrically around the black reference line, with their perpendicular distances to the red dash-dotted interface set as $d$ and $d + \Delta$, and their offsets from the black reference line set as $\Delta$.
\end{enumerate}
Here, $d$ denotes the shortest perpendicular distance from the HF point $(y_f, z_f)$ to the red dash-dotted interface (i.e., the projected contact line), and $\Delta$ is typically taken as one grid spacing. By placing the points in this controlled and symmetric manner, the final paraboloid accurately captures the local interface shape while enforcing the contact angle at the physically correct location $(x_\theta, y_\theta, z_\theta)$. The resulting interface normal $\textbf{n}_{l,h}$ is then evaluated from the gradient of the fitted surface and used to impose the contact angle boundary condition.

With these four height function values $h_{1}, h_{2}, h_{3}, h_{4}$, together with the values at the primary HF point $(h_f, y_f, z_f)$ and the contact point $(x_\theta, y_\theta, z_\theta)$ where the contact angle condition is prescribed, we now have sufficient information to fit the paraboloid $h(y, z)$ defined in Eq.~\eqref{eq:targetfittingparaboloid}. The coefficients $b_1$ through $b_6$ are determined by solving the following nonlinear system:

\begin{equation}\label{eq:systemequationsoffitting}
    \begin{aligned}
         &b_1y_1^2 + b_2z_1^2 + b_3y_1z_1 + b_4y_1 + b_5z_1 + b_6 - h_1 = 0,\\
         &b_1y_2^2 + b_2z_2^2 + b_3y_2z_2 + b_4y_2 + b_5z_2 + b_6 - h_2 = 0,\\
         &b_1y_3^2 + b_2z_3^2 + b_3y_3z_3 + b_4y_3 + b_5z_3 + b_6 - h_3 = 0,\\
         &b_1y_4^2 + b_2z_4^2 + b_3y_4z_4 + b_4y_4 + b_5z_4 + b_6 - h_4 = 0,\\
         &b_1y_f^2 + b_2z_f^2 + b_3y_fz_f + b_4y_f + b_5z_f + b_6 - h_f = 0,\\
         &\frac{\textbf{n}_{l,h} \cdot \textbf{n}_s}{|\textbf{n}_{l,h}||\textbf{n}_s|} = \cos\theta,
    \end{aligned}
\end{equation}
where the last equation enforces the contact angle constraint by requiring that the angle between the reconstructed interface normal $\textbf{n}_{l,h}$ and the solid surface normal $\textbf{n}_s$ matches the prescribed contact angle $\theta$. The interface normal vector derived from the fitted paraboloid is given by

\begin{equation}
\textbf{n}_{l,h} = \left[-1, \left.\frac{\partial h}{\partial y}\right|_{(y_\theta, z_\theta)}, \left.\frac{\partial h}{\partial z}\right|_{(y_\theta, z_\theta)} \right],
\end{equation}
with the partial derivatives evaluated at the contact point $(y_\theta, z_\theta)$. This formulation ensures that the reconstructed interface geometry not only interpolates the fluid height function data but also satisfies the prescribed local contact angle condition at the solid boundary.

Since Eq.~\eqref{eq:systemequationsoffitting} is nonlinear due to the contact angle constraint involving $\textbf{n}_{l,h}$, we solve it using Newton’s method. A contact line cell in which this nonlinear system is successfully solved is herein referred to as an \emph{interpolation cell}. Once the coefficients $b_1$ through $b_6$ are obtained, the height function at the center of the interpolation cell—represented by red dashed circles in Fig.~\ref{fig:paraboloidfit}(a)—can be evaluated directly from the fitted surface:
\begin{equation}
h_{\text{interp}}(y, z) = b_1 y^2 + b_2 z^2 + b_3 y z + b_4 y + b_5 z + b_6.
\end{equation}
This interpolated height function is used for curvature estimation in place of the original height function $(h_f, y_f, z_f)$ evaluated within the partial fluid region. This substitution serves two critical purposes. First, the discrete curvature formulas introduced (Eqs.~\eqref{eq:h1-kappa} and \eqref{eq:h1-kappa-2}) are constructed based on HF values located at full cell centers, whereas $(h_f, y_f, z_f)$ lies at the partial fluid centroid and may not be geometrically aligned with the uniform grid. Second, to construct a globally consistent height function field across neighboring interpolation cells, the interpolated height in each cell must be extended along the height direction (e.g., the $x$-direction in Fig.~\ref{fig:paraboloidfit}(a)) to its adjacent cells. For example, the height in cell $(i,j,k)$ must be extended to the neighboring cell $(i-1,j,k)$. This extension is performed under the following geometric constraint, referred to as the \emph{HF consistency condition}:
\begin{equation}
    |h_{i,j,k} - h_{i-1,j,k}| = \Delta,
\end{equation}
which ensures that the vertical interface reconstructed by the HF method remains consistent across adjacent cells. Violation of HF consistency can lead to discontinuities or artificial kinks in the reconstructed interface, and in severe cases, may cause numerical fragmentation of the liquid domain. This issue has also been reported by Han et al.~\cite{han2021consistent} for curvature estimation near regular (planar) solid boundaries. In the present study, enforcing HF consistency proves essential for maintaining geometric robustness in curvature computation near complex 3D solid surfaces.

After determining the coefficients of the fitted paraboloid in Eq.~\eqref{eq:targetfittingparaboloid} for all interpolation cells, the next step is to compute the height function values at the centers of neighboring solid cells. This process must be treated on a case-by-case basis, depending on the adjacency relationships between solid and mixed cells. For example, consider a solid cell adjacent to a single mixed cell, such as $(i, j{+}1, k{-}1)$ in Fig.~\ref{fig:paraboloidfit}(a). If a paraboloid has been successfully fitted in any interpolation cell along the same row $(j{+}1, k)$, it can be directly used to evaluate the HF value at $(i, j{+}1, k{-}1)$. 
In contrast, if a solid cell is adjacent to multiple mixed cells—e.g., $(i, j, k{-}1)$, which neighbors both $(i, j{-}1, k{-}1)$ and $(i, j, k)$—the paraboloids fitted in both rows $(j{-}1, k{-}1)$ and $(j, k)$ can each provide a candidate HF value for $(i, j, k{-}1)$. In such cases, the final HF value is computed by averaging the predictions from all applicable paraboloids. 
Finally, for solid cells that are not directly adjacent to any mixed cells, we search among their diagonal neighbors to identify mixed cells with valid fitted paraboloids in the corresponding rows. These neighboring paraboloids are then used to define the HF value following the same averaging procedure described above.

The complete procedure for the proposed paraboloid fitting method is summarized below.

\begin{algorithm}[H]\label{alg:hfcontactangle}
  \caption{Paraboloid fitting procedure for enforcing the contact angle condition}
  \begin{enumerate}
  
  \item Apply the standard height function method to evaluate the interface height in cells far from solid boundaries. For mixed cells near the solid surface, define the local height function using Algorithm~\ref{alg:hfdefinition};
  
  \item In cells where the liquid/gas interface intersects the solid boundary, construct a preliminary (pre-fitted) paraboloid surface using Eq.~\eqref{eq:prefittingparaboloid} based on surrounding HF values (see Fig.~\ref{fig:paraboloidfit}(a));
  
  \item Evaluate the HF values at the four red auxiliary points shown in Fig.~\ref{fig:paraboloidfit}(b) using the pre-fitted paraboloid;
  
  \item Combine the four auxiliary points, the fluid-region HF value $(h_f, y_f, z_f)$, and the contact angle constraint at the contact point $(x_\theta, y_\theta, z_\theta)$ to construct the final paraboloid surface by solving the nonlinear system in Eq.~\eqref{eq:systemequationsoffitting};
  
  \item Use the fitted paraboloid to evaluate the HF value at the center of the interpolation cell and at neighboring solid cells as needed.
  
  \end{enumerate}
\end{algorithm}

If the nonlinear system in Eq.~\eqref{eq:systemequationsoffitting} cannot be constructed or fails to converge during the iterative solution (e.g., via Newton's method), this is typically due to one of the following reasons:

\begin{enumerate}
    \item The local curvature of the liquid/gas interface is excessively large, rendering the selected HF points insufficient to define a well-posed fitting problem.
    
    \item Even when the system is mathematically solvable, the geometry of the selected HF points may lead to an unphysical or unstable fitted surface.
    
    \item An insufficient number of valid HF points is available within the local stencil, preventing successful construction of the pre-fitted paraboloid via Eq.~\eqref{eq:prefittingparaboloid}.
    
    \item In degenerate configurations where the solid boundary is nearly parallel to the $xy$-plane and the contact angle $\theta$ is small, it may be impossible to define a valid HF in the $x$- or $y$-directions. Specifically, the required bounding cells (with $c = 0$ and $c = c_s$) may not simultaneously exist within the accumulation stencil $[i{-}4, i{+}4]$, as per Eq.~\eqref{eq:hsum}. Although a HF may still be definable in the $z$-direction, the reconstructed interface may no longer intersect the solid surface, preventing enforcement of the contact angle condition.
\end{enumerate}

To improve robustness in cases (1)--(3), the HF stencil can be expanded from a $3 \times 3$ to a $5 \times 5$ region in the $yz$-plane. If the fitting still fails after this enlargement, we revert to the simpler linear extension method described in \S~\ref{sec:linearfitting}, which provides a less accurate but stable estimate of the HF at the cell center. To address issue (4), we introduce a \emph{vertical height function} (VHF) method, enabling contact angle enforcement in geometries where standard horizontal HFs fail. This approach, not available in the original \textit{Basilisk} implementation, is presented in the next section. 

Once all height function values near the contact line are determined and assigned to their corresponding cell centers, the interface curvature and normal vector can be consistently evaluated using the discrete HF-based expressions introduced earlier (Eqs.~\eqref{eq:h1-kappa} and \eqref{eq:hn}).

\subsection{Paraboloid fitting — additional ``vertical'' height function method}\label{sec:paraboloidfittingvertical}

\begin{figure}
    \includegraphics[height=0.28\textwidth]{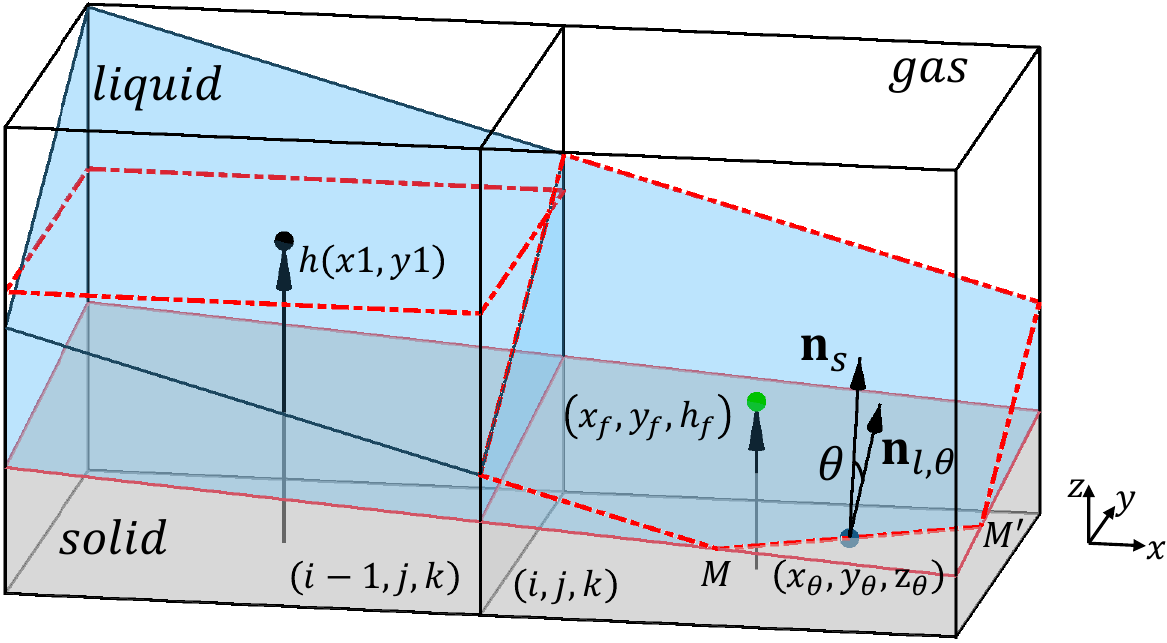}
    \centering
    \caption{
    Sketch of the ``vertical'' height function method.
    }
    \label{fig:paraboloidverticalHF}
\end{figure}

\begin{figure}
    \includegraphics[height=0.36\textwidth]{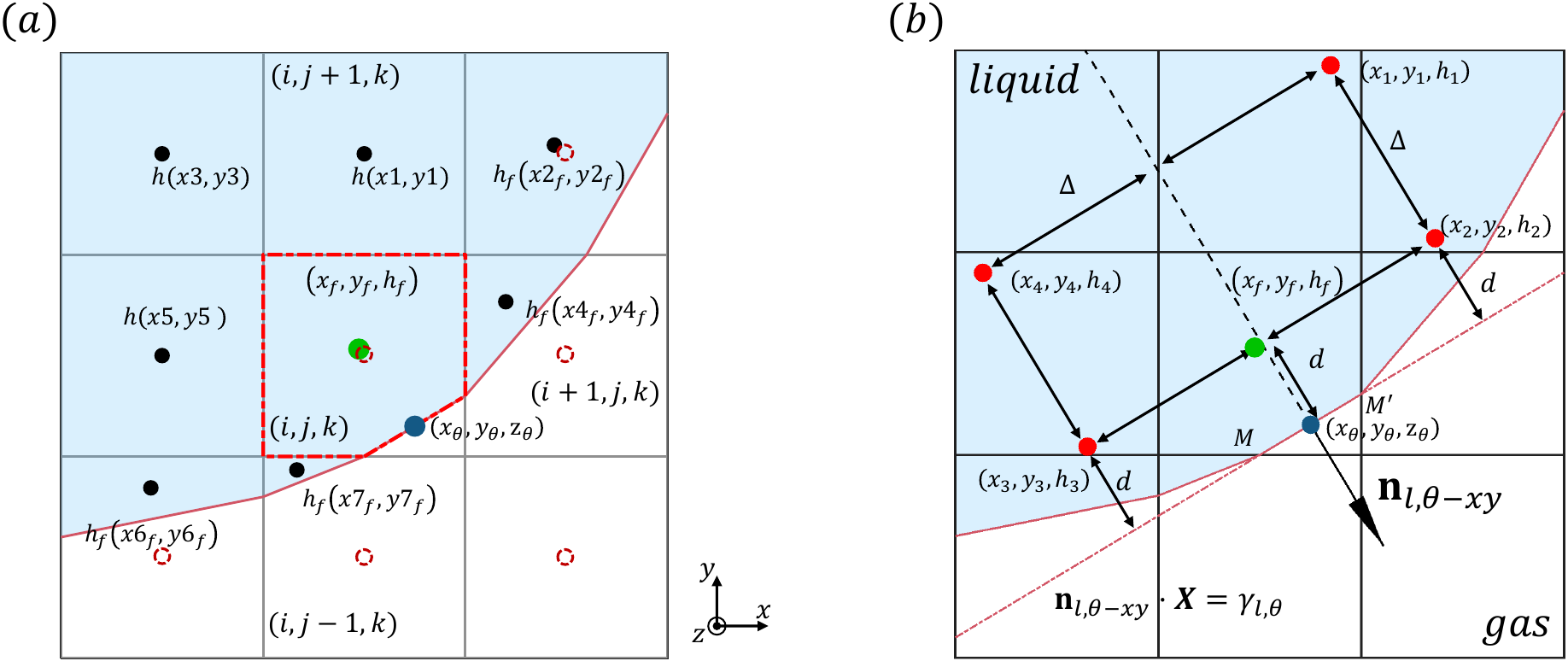}
    \centering
    \caption{
    Schematic diagram illustrating the pre-fitting of a paraboloid to determine the four red interpolation points (shown in panel b). $\textbf{n}_{l,\theta-xy}$ denotes the projection of the contact line normal $\textbf{n}_{l,\theta}$ onto the $xy$-plane.
    }
    \label{fig:paraboloidverticalfit}
\end{figure}

In \S~\ref{sec:paraboloidfittinghorizontal}, the height function in mixed cells was constructed based on an interface oriented perpendicular to the direction of height accumulation. This approach, analogous to the Simple Line Interface Calculation (SLIC) scheme, constitutes the ``horizontal'' height function (HHF) method, widely used for regular boundaries~\cite{afkhami2008height}. However, as discussed previously, the HHF method becomes invalid when the solid surface is nearly parallel to the $xy$-plane and the contact angle $\theta$ is small. In such cases, no valid interface can be reconstructed in any Cartesian direction, corresponding to issue (4), and the HHF strategy fails entirely.

To overcome this limitation, we extend the height function concept using a Piecewise Linear Interface Construction (PLIC) representation, thereby formulating a ``vertical'' height function (VHF) method. Originally proposed by Afkhami et al.~\cite{afkhami2008height} for regular boundaries, the VHF method allows contact angle enforcement even when $\theta < 45^\circ$ or $\theta > 135^\circ$. Despite its broader applicability, this method is computationally more involved and has not been included in existing implementations such as \textit{Basilisk}, which restricts the permissible contact angle range.

In contrast, the present work develops a generalized VHF scheme applicable to arbitrarily complex solid boundaries. The key innovation lies in relaxing the assumption that the reconstructed interface must be aligned with any specific coordinate direction. Instead, the liquid/gas interface $\Gamma_l$ is constructed to intersect the solid boundary $\Gamma_s$ at the prescribed contact angle $\theta$, regardless of orientation, as illustrated in Fig.~\ref{fig:paraboloidverticalHF}. This flexibility enables HF values to be defined even in geometrically degenerate scenarios, thereby restoring the ability to enforce accurate contact angle conditions.

To construct the vertical height function, we first reconstruct the interface in cell $(i, j, k)$ using the contact angle–dependent normal vector $\textbf{n}_{l,\theta}$ (red dash-dotted surface in Fig.~\ref{fig:paraboloidverticalHF}). The centroid of this reconstructed interface (green point) defines the fluid-region HF value $h_f$, while the midpoint of the intersection segment $MM'$ (blue point) identifies the location where the contact angle constraint is enforced. Aside from this change in interface orientation (step 3 in Algorithm~\ref{alg:hfdefinition}), all remaining steps of the HHF method are preserved and reused in the VHF method. These include:
\begin{itemize}
    \item Pre-fitting a paraboloid using Eq.~\eqref{eq:prefittingparaboloid} based solely on HF points within the fluid region (see Fig.~\ref{fig:paraboloidverticalfit}(a));
    \item Computing the four auxiliary HF points positioned relative to the contact line (Fig.~\ref{fig:paraboloidverticalfit}(b));
    \item Solving the nonlinear system in Eq.~\eqref{eq:systemequationsoffitting} to fit the final paraboloid incorporating the contact angle condition.
\end{itemize}

An important distinction between the HHF and VHF methods arises in the final step of Algorithm~\ref{alg:hfcontactangle}, namely the extension of height function values from the interpolation cell to neighboring cells. In the HHF framework, the HF is extended from the interpolation cell (typically a mixed cell) into adjacent solid cells. In contrast, the VHF method extends both the HF value and the prescribed contact angle from the liquid side of the contact line toward the gas side. This modification is essential for accurately capturing interface curvature in degenerate contact line geometries. Specifically, when computing the HF values at the four red interpolation points located at distances $d$ and $d{+}\Delta$ from the contact point (see Fig.~\ref{fig:paraboloidverticalfit}(b)), the contact line normal $\textbf{n}_{l,\theta}$ is projected onto the $xy$-plane, yielding the in-plane direction $\textbf{n}_{l,\theta-xy}$. This projection guides the placement of the auxiliary points, in contrast to the HHF method (Fig.~\ref{fig:paraboloidfit}(b)), where the projection of the solid normal $\textbf{n}_s$ onto the $yz$-plane, $\textbf{n}_{s-yz}$, was used instead.

Overall, the vertical height function method developed here extends the robustness of contact angle enforcement to a broader range of interface geometries, including cases involving steep or nearly horizontal solid boundaries. The integration of this VHF capability into the paraboloid fitting framework ensures accurate curvature and normal estimation even in challenging three-phase configurations.

\subsection{Hysteresis}\label{sec:hysteresis}

\begin{figure}
    \includegraphics[height=0.36\textwidth]{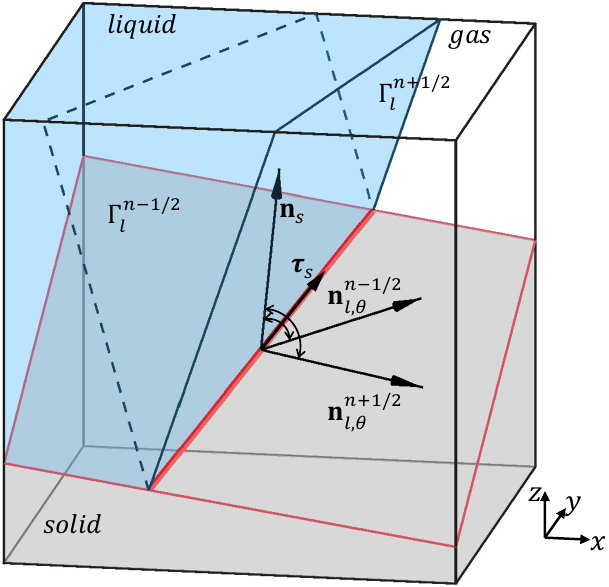}
    \centering
    \caption{
    Sketch illustrating the computation of $\theta^{n+1/2}$ within the hysteresis window. The dashed line represents the liquid/gas interface at the previous time step, $\Gamma_l^{n-1/2}$, while the solid line represents the updated interface at $n+1/2$, $\Gamma_l^{n+1/2}$.
    }
    \label{fig:hysteresis}
\end{figure}

In the preceding sections, the contact angle $\theta$ has been accurately enforced at the contact line using paraboloid-fitting methods. However, in many practical applications, microscale surface roughness or chemical heterogeneities give rise to contact angle hysteresis, wherein the contact angle varies within a bounded range while the contact line remains pinned. This subsection introduces a numerical strategy to model such hysteresis effects.

To simulate contact angle hysteresis, we employ an iterative procedure following the update of the volume fraction through advection. Specifically, in each contact line cell, the contact angle is determined via a bisection method within a prescribed hysteresis window $[\theta_{\text{rec}}, \theta_{\text{adv}}]$, where $\theta_{\text{rec}}$ and $\theta_{\text{adv}}$ denote the receding and advancing contact angles, respectively. The objective is to identify a value of $\theta$ such that the reconstructed liquid/gas interface, defined by its normal vector $\textbf{n}_{l,\theta}$, intersects the solid boundary at the same contact line position as in the previous time step. During this iteration, the volume fraction remains fixed at its updated value. Figure~\ref{fig:hysteresis} schematically illustrates this procedure: the interface at the previous time step, $\Gamma_l^{n-1/2}$ (dashed line), and the updated interface, $\Gamma_l^{n+1/2}$ (solid line), intersect the solid surface at an identical contact line location. This condition ensures that the contact line remains pinned, capturing the essential physical behavior associated with hysteresis.

In 3D, it is important to note that the reconstructed interface normal $\textbf{n}_{l,\theta}$, computed via Eq.~\eqref{eq:n-l-theta}, depends not only on the contact angle $\theta$ but also on the orientation of the contact line, represented by the tangential vector $\boldsymbol{\tau}_s$. During the hysteresis iteration, $\boldsymbol{\tau}_s$ is held fixed to preserve the contact line orientation and maintain geometric consistency across time steps. Once the appropriate value of $\theta$ is determined through this iterative search, the paraboloid fitting procedure described in Algorithm~\ref{alg:hfcontactangle} is applied to enforce the contact angle condition, as in the standard non-hysteresis case. This strategy enables accurate and stable modeling of contact line dynamics under hysteresis, while preserving the geometric fidelity of interface reconstruction near complex solid surfaces.

\section{Results}\label{sec:results}

This section presents a series of numerical tests designed to evaluate the accuracy, robustness, and conservation properties of the proposed numerical schemes. The objectives are twofold: to validate the volume-of-fluid (VOF) advection method, particularly in the presence of complex solid boundaries, and to assess the accuracy of the contact angle implementation, including hysteresis modeling.

In \S~\ref{sec:resultsVOFadvection}, we first examine the volume conservation performance of the conservative advection scheme described in \S~\ref{sec:conservationadvectionscheme} using 3D test cases, followed by an evaluation of the redistribution-based advection scheme introduced in \S~\ref{sec:smallcutcell}. These tests illustrate how the redistribution scheme alleviates the time step constraints imposed by small cut cells while maintaining interface integrity. Next, \S~\ref{sec:surfacetensiondrivendropletspreading} investigates the classical surface-tension-driven droplet spreading problem, verifying both the accuracy of the proposed contact angle enforcement method and the symmetry of contact line motion. In \S~\ref{sec:hysteresis-sheardriven}, we simulate droplet dynamics under shear flow with contact angle hysteresis, demonstrating the effectiveness of the iterative 3D hysteresis model introduced in \S~\ref{sec:hysteresis}. Finally, in the remaining examples,
we examine droplet behavior in the presence of complex 3D solid geometries and compare our results with existing numerical and experimental benchmarks to highlight the generality and practical applicability of the method.

Unless otherwise specified, all simulations in this section employ the redistribution advection scheme described in \S~\ref{sec:smallcutcell}, which avoids CFL restrictions associated with small cut cells. All solid surfaces are treated as no-slip boundaries, and grid resolution is expressed in terms of the number of cells per initial droplet radius, abbreviated as \textit{cpr}.

\subsection{Geometric VOF advection scheme}\label{sec:resultsVOFadvection}

\subsubsection{Concentric sphere model with rotating liquid shell}\label{sec:resultsVOFadvection1}

\begin{figure}
    \includegraphics[height=0.36\textwidth]{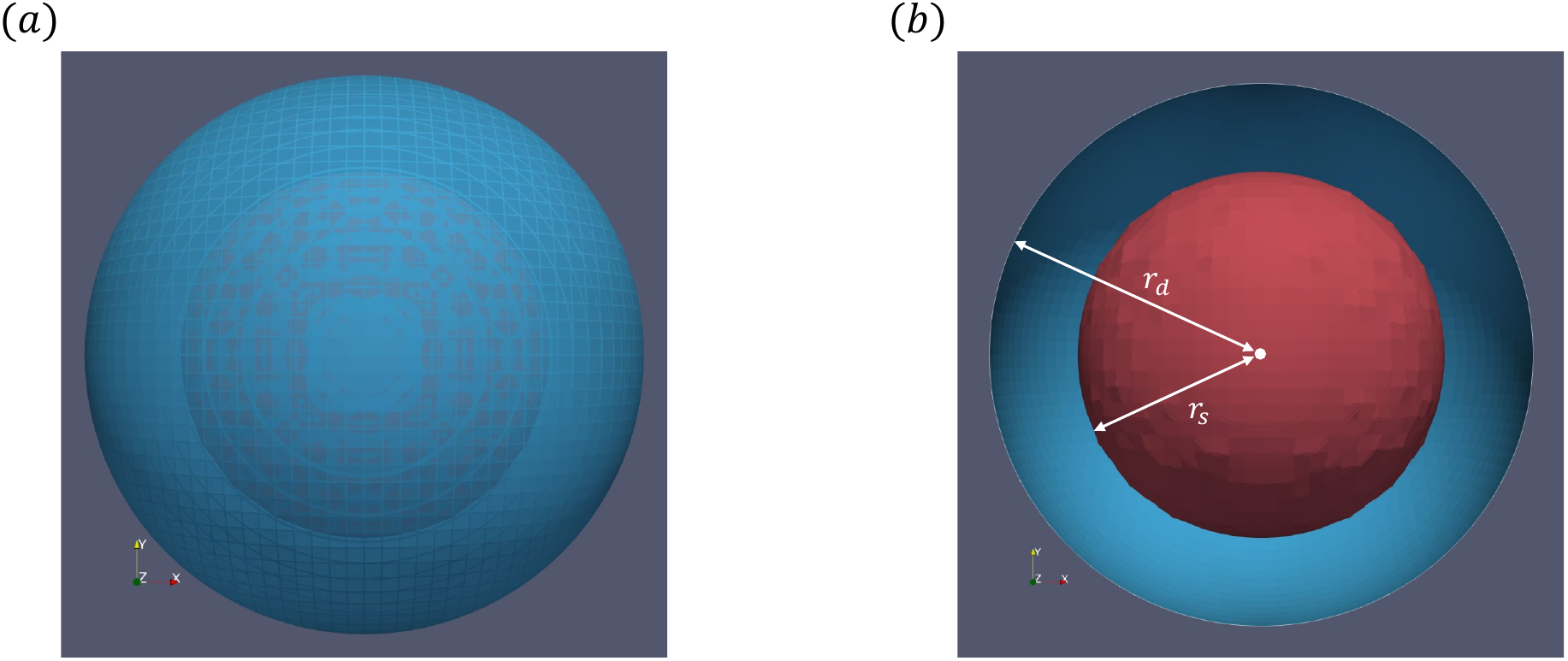}
    \centering
    \caption{
    Sketch of a concentric sphere model with a rotating liquid shell. $r_s$ and $r_d$ denote the radii of the solid sphere and the outer boundary of the liquid shell, respectively.
    }
    \label{fig:liquidring}
\end{figure}

\begin{figure}
    \includegraphics[height=0.36\textwidth]{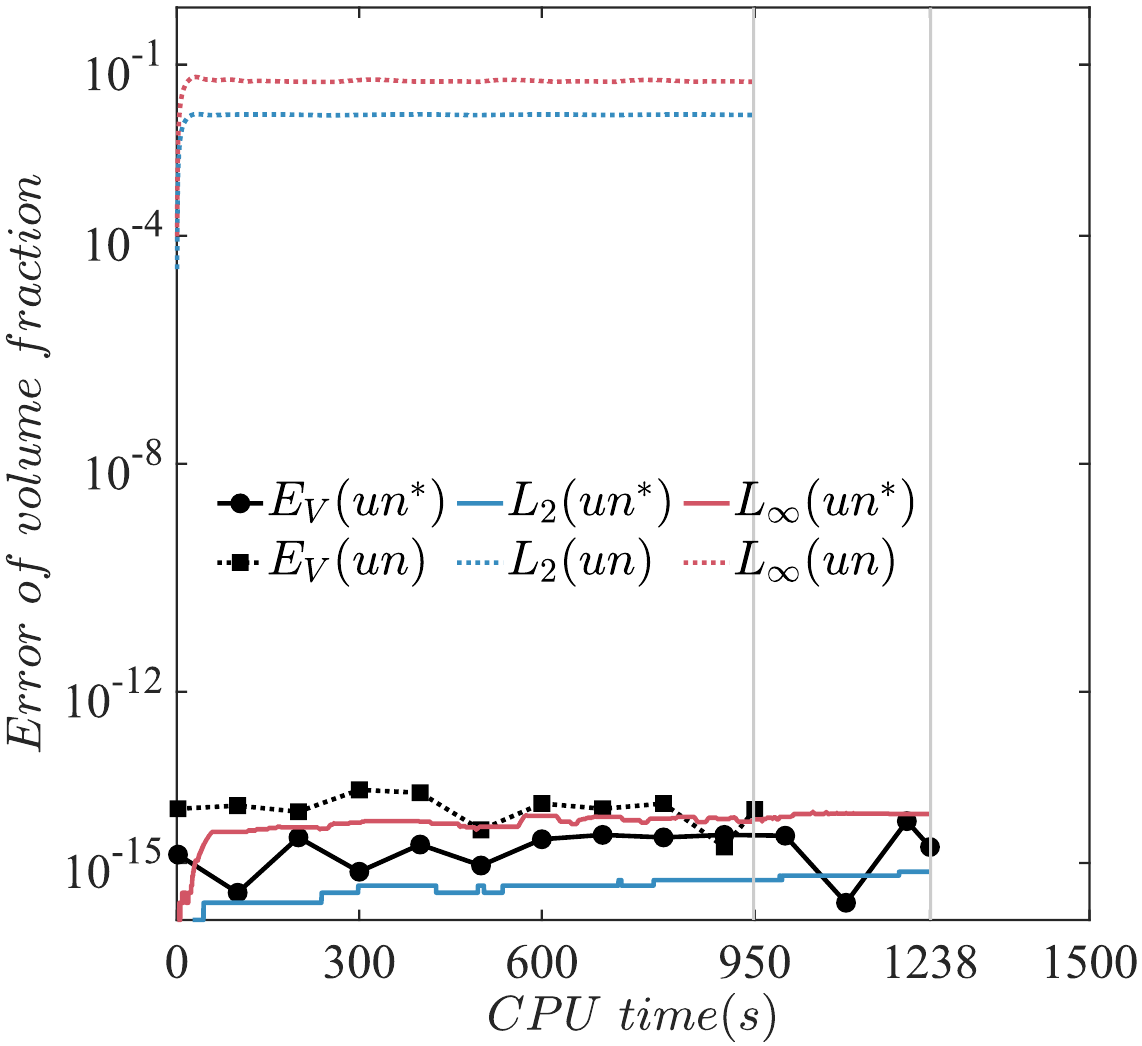}
    \centering
    \caption{
    Volume conservation error during one full rotation of the liquid shell. Blue and red lines indicate the $L_2$ and $L_\infty$ norms of local volume errors in mixed cells near the solid surface. Dashed lines correspond to the traditional advection using $un$, while solid lines correspond to the modified advection using ${un}^*$ (see \S~\ref{sec:conservationadvectionscheme}). The black line represents the global volume loss, $E_V = |V - V_0| / V_0$.
    }
    \label{figr:liquidring}
\end{figure}

This test highlights the inability of traditional geometric VOF advection schemes to preserve local volume near embedded boundaries in 3D, and demonstrates the effectiveness of the correction method proposed in \S~\ref{sec:conservationadvectionscheme}.

We consider a concentric spherical configuration (Fig.~\ref{fig:liquidring}) consisting of a stationary solid sphere of radius $r_s = 0.4\,\mathrm{m}$, surrounded by a liquid shell with outer radius $r_d = 0.6\,\mathrm{m}$. The liquid shell undergoes rigid-body rotation about the $z$-axis, driven by the prescribed velocity field:
\begin{equation}\label{eq:givingu1}
        \begin{aligned}
                &u_x = 2\pi (y - y_c), \\
                &u_y = 2\pi (x_c - x), \\
                &u_z = 0,
        \end{aligned}
\end{equation}
where $(x_c, y_c)$ denotes the center of the solid sphere. The simulation employs a uniform grid with a resolution of $19.2$ cells per radius (cpr). All cells within the liquid shell are initially full ($c = c_s$), so no interface reconstruction is required during advection. Consequently, the redistribution advection scheme described in \S~\ref{sec:smallcutcell}, which is specifically designed for interfacial cells, is not used in this test.

The presence of the embedded solid generates extremely small cut cells near the sphere surface, with $c_s / s_f \sim 2.3 \times 10^{-3}$, imposing severe time step constraints. To render the test feasible, cells with a line fraction $s_l < 0.15$ are discarded during initialization, slightly modifying the sphere geometry. This results in a minimum value of $c_s / s_f$ of approximately $7.3 \times 10^{-2}$. We emphasize that this pre-processing step is applied only in this test case.

To quantify volume conservation, we define two local error norms:
\begin{equation}
\begin{aligned}
    L_2 &= \sqrt{\frac{1}{N} \sum_{i=1}^{N} (c_{s,i} - c_i)^2}, \\
    L_\infty &= \max_{1 \leq i \leq N} \left| c_{s,i} - c_i \right|, 
\end{aligned}
\end{equation}
where the summation is taken over all mixed cells adjacent to the embedded solid boundary. Additionally, the global volume error is defined as
\[
E_V = \frac{|V - V_0|}{V_0},
\]
where $V$ denotes the total liquid volume after one full rotation.

Figure~\ref{figr:liquidring} compares the traditional geometric advection scheme with the proposed corrected scheme using the modified velocity ${un}^*$. After one rotation, the corrected method achieves local volume errors of $L_\infty = 7.9 \times 10^{-15}$ and $L_2 = 7.0 \times 10^{-16}$, effectively reaching machine precision. In contrast, the traditional scheme accumulates substantial local errors, with $L_\infty = 6.1 \times 10^{-2}$ and $L_2 = 1.4 \times 10^{-2}$, despite maintaining excellent global conservation ($E_V \sim 10^{-14}$ for both methods). These results confirm that the traditional advection scheme, which neglects the influence of the embedded solid when computing fluxes in mixed cells, introduces systematic local errors—a limitation previously identified in 2D~\cite{huang20252d} and now demonstrated to persist in 3D. The proposed correction, by modifying the effective advection width ${un}^*$, restores consistency and ensures local volume conservation.

To evaluate computational cost, we measure the CPU time required for one full rotation under different time step conditions. When both methods adopt the same CFL-constrained time step (5,434 steps), the corrected scheme incurs a 30\% increase in runtime (1,238\,s vs. 950\,s) but eliminates local conservation errors. However, if the timestep restriction imposed by small cut cells is disregarded, as is often done in practical simulations, the traditional method completes in only 407 steps (93\,s), over ten times faster. This underscores a key limitation of the conservative scheme: although it achieves machine-accurate volume conservation, its strict time step constraints can result in substantial computational cost. In practice, this motivates the development of the redistribution-based advection method introduced in \S~\ref{sec:smallcutcell}, which maintains conservation while removing the timestep restriction.

\subsubsection{Rotating liquid shell with a conical cutout}\label{sec:resultsVOFadvection2}

\begin{figure}
    \includegraphics[height=0.36\textwidth]{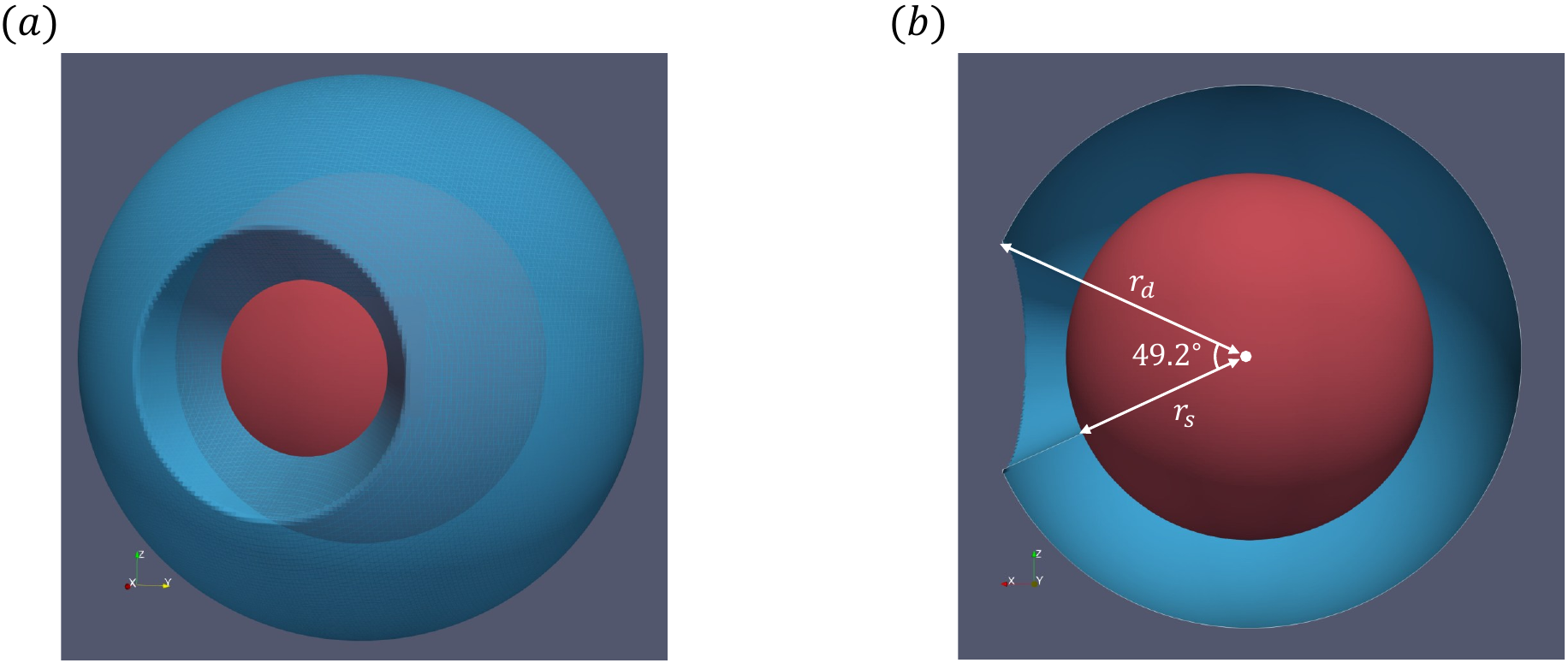}
    \centering
    \caption{
    Schematic of a rotating liquid shell with a conical cutout surrounding a solid sphere.
    }
    \label{fig:Zalesak}
\end{figure}

\begin{figure}
    \includegraphics[width=1\textwidth]{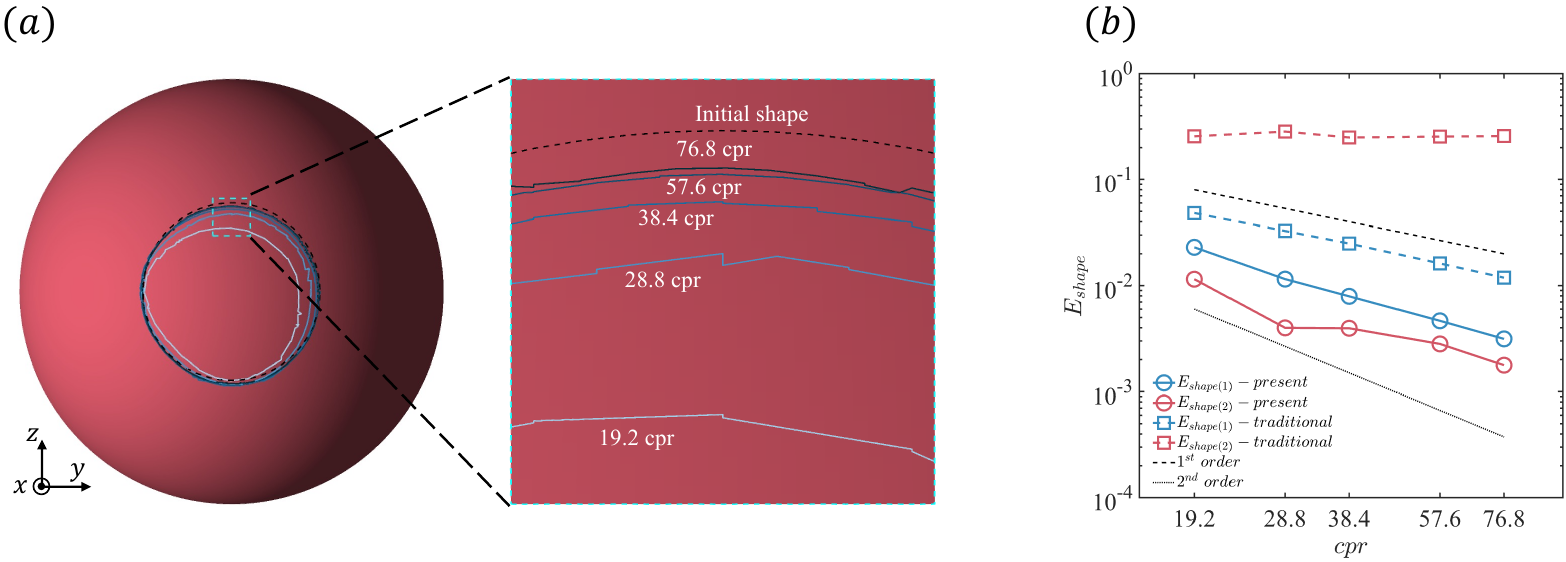}
    \centering
    \caption{
    (a) Contact line contours after 10 full rotations at various grid resolutions using the present redistribution advection method. The dashed black line represents the initial contact line, and increasingly darker blue lines correspond to higher resolutions. (b) Shape errors versus resolution: $E_{\text{shape(1)}}$ (global domain) and $E_{\text{shape(2)}}$ (mixed cells). Solid lines correspond to the present redistribution advection method, while dashed lines correspond to the traditional geometric VOF.
    }
    \label{figr:Zalesak}
\end{figure}

This test assesses the accuracy and robustness of the redistribution advection scheme described in \S~\ref{sec:smallcutcell} by simulating a rotating liquid shell with a conical cutout — a 3D analog of Zalesak disk. The configuration introduces a persistent contact line and complex interface deformations over time, providing a stringent benchmark for volume-preserving advection near embedded solid boundaries.

The liquid shell and solid sphere retain the same dimensions as in \S~\ref{sec:resultsVOFadvection1}, with the apex angle of the conical cutout set to $49.2^\circ$ (Fig.~\ref{fig:Zalesak}). The contact angle is fixed at $\theta = 90^\circ$. The rotational velocity is prescribed as
\begin{equation}\label{eq:givingu2}
        \begin{aligned}
                &{u}_x=2\pi (y - y_c) + 2\pi (z - z_c), \\
                &{u}_y=2\pi (x_c - x), \\
                &{u}_z=2\pi (x_c - x),
        \end{aligned}
\end{equation}
where $(x_c, y_c, z_c)$ denotes the center of the solid sphere. Five grid resolutions are considered: 19.2, 28.8, 38.4, 57.6, and 76.8 cpr. To quantify shape preservation during rotation, we define the relative shape error as
\begin{equation}
    E_{\text{shape(1/2)}} = \frac{\sum^{N_{(1/2)}}_i V_i\,|c_i(T) - c_i(0)|}{\sum^{N_{(1/2)}}_i V_i\, c_i(0)},
\end{equation}
where $N_{(1)}$ includes all grid cells, $N_{(2)}$ includes only mixed cells near the contact line, $V_i$ is the cell volume, and $c_i(t)$ is the volume fraction at time $t$.

After 10 full rotation cycles, the results are shown in Fig.~\ref{figr:Zalesak}. In panel (a), the contact line contours show that, with increasing resolution, the redistributed scheme (solid blue lines) better preserves the initial geometry, exhibiting minimal distortion even after prolonged rotation. Panel (b) presents the shape errors $E_{\text{shape(1/2)}}$ for both advection methods. The redistribution scheme achieves significantly lower errors, particularly in mixed cells, and demonstrates a clear convergence trend between first- and second-order accuracy. In contrast, the traditional geometric VOF method yields larger errors and does not show consistent convergence in mixed cells, underscoring its limitations near embedded boundaries.

We further compare performance at 38.4 cpr over 10 rotation cycles (Tab.~\ref{tab:liquidrotation}). The traditional method requires 5{,}959 steps and 24,600\,s of CPU time, with a global shape error of 0.0249 and a significantly larger error in mixed cells (0.2493). The redistribution advection scheme reduces these errors to 0.0079 and 0.0039, respectively—an order-of-magnitude improvement—while incurring only a 20\% increase in CPU time. By contrast, enforcing strict time step constraints due to small cut cells (e.g., $c_s / s_f \sim 5.1 \times 10^{-2}$) would require 11,248 steps for a single rotation cycle, rendering long-term simulations impractical. The proposed scheme circumvents this bottleneck, enabling accurate long-time advection near embedded boundaries without excessive computational cost.

\begin{table}
  \centering
  \caption{Results of the rotating liquid shell with a conical cutout at $t = 10T$ (38.4 cpr).}
  \vspace{5pt}
  \renewcommand\arraystretch{1.5}
  \begin{tabular}{l|cccc}
    \hline
    Scheme & Steps & CPU time (s) & $E_{\text{shape(1)}}$ & $E_{\text{shape(2)}}$ \\
    \hline
    Present (Redistribution) & 5{,}847 & 29{,}212 & 0.0079 & 0.0039 \\
    Traditional VOF & 5{,}959 & 24{,}600 & 0.0249 & 0.2493 \\
    \hline
  \end{tabular}
  \label{tab:liquidrotation}
\end{table}

\subsection{Surface tension-driven droplet spreading}\label{sec:surfacetensiondrivendropletspreading}

This test examines the static spreading of a droplet driven solely by surface tension, with gravity neglected. This classical benchmark problem~\cite{afkhami2009mesh,legendre2015comparison,liu2015diffuse} is widely used to evaluate the accuracy of contact angle implementation methods.

In the absence of external forces, the droplet evolves toward an equilibrium shape determined exclusively by the prescribed contact angle. Consequently, the final droplet configuration should be independent of the position and orientation of the embedded solid boundary. This property provides a stringent test of the geometric consistency of the numerical method, particularly when complex or inclined surfaces are involved. In 3D simulations, we also monitor the \textit{symmetry of the contact line} during spreading, which serves as a sensitive indicator of how accurately the contact angle condition is enforced throughout the evolution.

Unless otherwise specified, the physical parameters used in this section are:
\begin{itemize}
  \item Liquid and gas densities: $\rho_l = \rho_g = 1~\text{kg}/\text{m}^3$,
  \item Liquid and gas viscosities: $\mu_l = \mu_g = 7.5 \times 10^{-4}~\text{kg}/(\text{m} \cdot \text{s})$,
  \item Surface tension coefficient: $\sigma = 0.1~\text{N}/\text{m}$.
\end{itemize}
The grid resolution is fixed at 25.6 cpr (cells per initial droplet radius).

\subsubsection{Flat plate}\label{sec:flatplate}

\begin{figure}
    \includegraphics[height=0.36\textwidth]{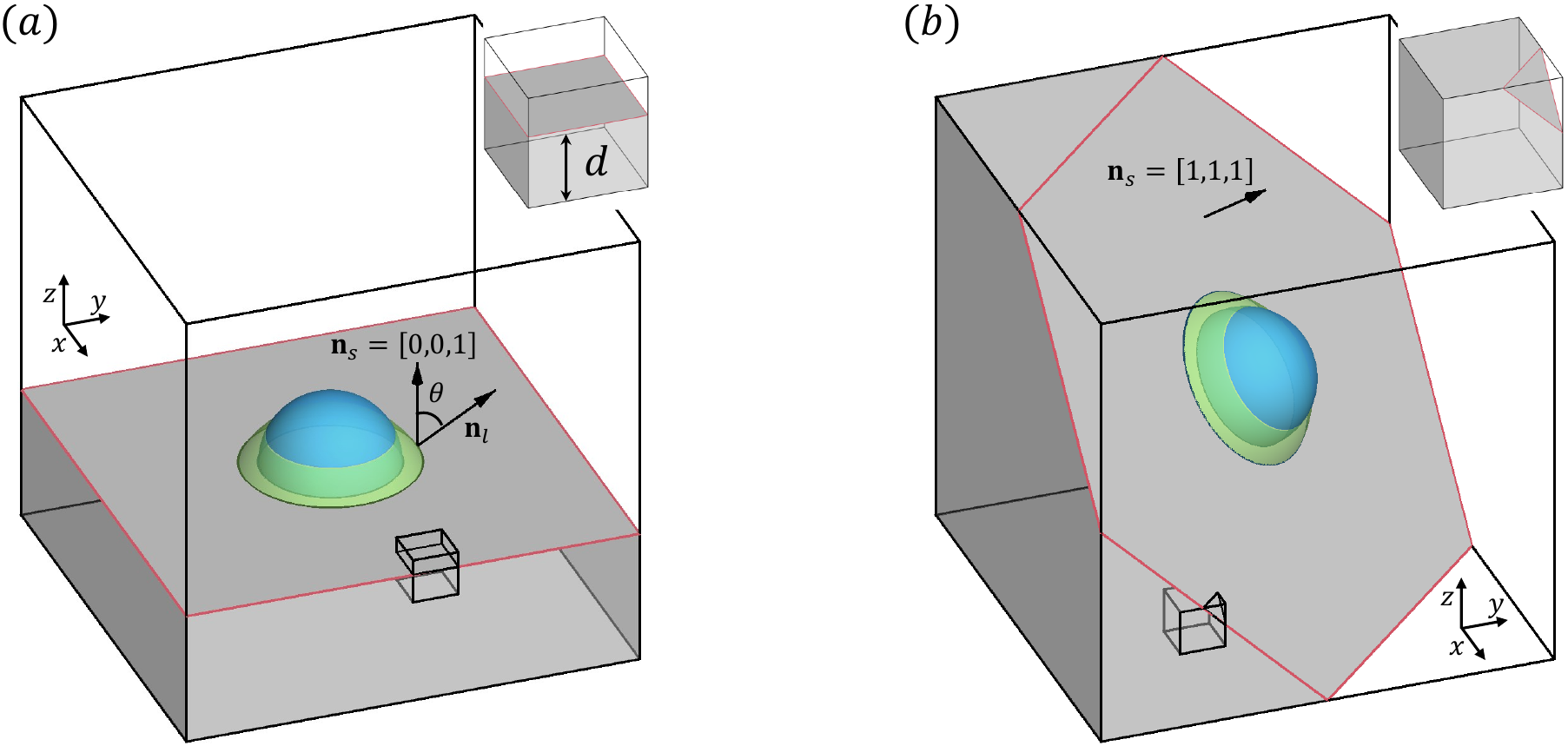}
    \centering
    \caption{
    Schematic of the initial (blue) and equilibrium (green) droplet shapes on a flat plate, with surface normal vectors (a) $\textbf{n}_s = [0, 0, 1]$ and (b) $\textbf{n}_s = [1, 1, 1]$.
    }
    \label{fig:planedroplet}
\end{figure}

We begin with a classical benchmark in which a droplet spreads over a flat solid plate. Two plate orientations are considered: horizontal ($\mathbf{n}_s = [0, 0, 1]$) and oblique ($\mathbf{n}_s = [1, 1, 1]$), as illustrated in Fig.~\ref{fig:planedroplet}. A hemispherical droplet with an initial contact angle of $\theta_0 = 90^\circ$ and radius $r_0 = 0.2$\,m is placed on the surface. The prescribed contact angle is then modified, causing the droplet to spread or retract under surface tension until reaching equilibrium.

Six target contact angles are tested: $15^\circ$, $40^\circ$, $70^\circ$, $110^\circ$, $140^\circ$, and $165^\circ$. To ensure stable adhesion during rapid retraction (particularly for large $\theta$), the liquid and gas viscosities are increased to $\mu_l = \mu_g = 5.0 \times 10^{-3}~\text{kg}/(\text{m}\cdot\text{s})$ for $\theta = 140^\circ$, and $\mu_l = \mu_g = 2.5 \times 10^{-2}~\text{kg}/(\text{m}\cdot\text{s})$ for $\theta = 165^\circ$.

The theoretical equilibrium contact line radius is given by~\cite{afkhami2009mesh}:
\begin{equation}\label{eq:rt}
    r_t = \sqrt[3]{\frac{3v}{\pi(2+\cos{\theta})(1-\cos{\theta})^2}}\sin{\theta},
\end{equation}
where $v$ denotes the droplet volume. To evaluate robustness against boundary embedding, four configurations are examined in which the flat surface intersects the computational cells at solid volume fractions of $d = 0.2$, $0.4$, $0.6$, and $0.8$. Two methods for enforcing the contact angle condition are compared: the linear fitting method (described in \S~\ref{sec:linearfitting}) and the proposed paraboloid fitting approach (\S~\ref{sec:paraboloidfittinghorizontal}–\ref{sec:paraboloidfittingvertical}). For reference, an additional simulation on a regular Cartesian grid without embedded solids is included using \textit{Basilisk}’s built-in ``contact.h'' module (labeled ``non-embedded'').

\paragraph{Accuracy of spreading radius and boundary insensitivity}

\begin{figure}
    \includegraphics[width=1\textwidth]{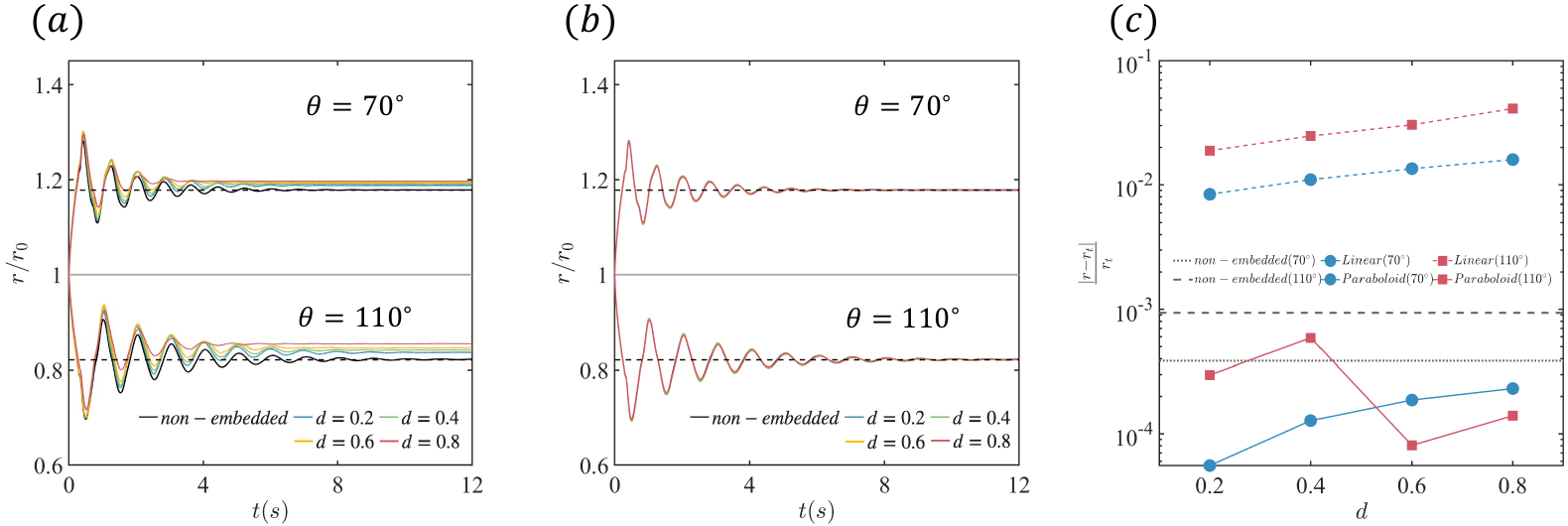}
    \centering
    \caption{
    Droplet spreading on a flat plate ($\textbf{n}_s = [0, 0, 1]$): (a) linear fitting, (b) paraboloid fitting, and (c) relative error in the equilibrium radius compared to theory. ``Non-embedded'' denotes the reference case without an embedded solid.
    }
    \label{figr:flatplanedroprt}
\end{figure}

Figure~\ref{figr:flatplanedroprt} compares the evolution of the droplet contact line radius across different methods and embedding depths for $\theta = 70^\circ$ and $110^\circ$. Panels (a) and (b) show results obtained using the linear and paraboloid fitting methods, respectively, while panel (c) presents the relative error with respect to Eq.~\eqref{eq:rt}. With the linear fitting approach, the error of the equilibrium radius increases with $d$, indicating a strong dependence on the boundary position. Larger values of $d$ lead to greater deviations from the non-embedded reference. In contrast, the paraboloid fitting method produces consistent results across all $d$ values, with small and nearly uniform relative errors. This demonstrates that the linear fitting scheme fails to accurately impose the contact angle at the true solid interface, particularly when the solid boundary lies deep within the mixed cell. Conversely, the proposed paraboloid fitting scheme enforces the contact angle directly at the geometric interface, yielding accurate and boundary-insensitive results.

\paragraph{Contact line symmetry during spreading}

\begin{figure}
    \includegraphics[width=1\textwidth]{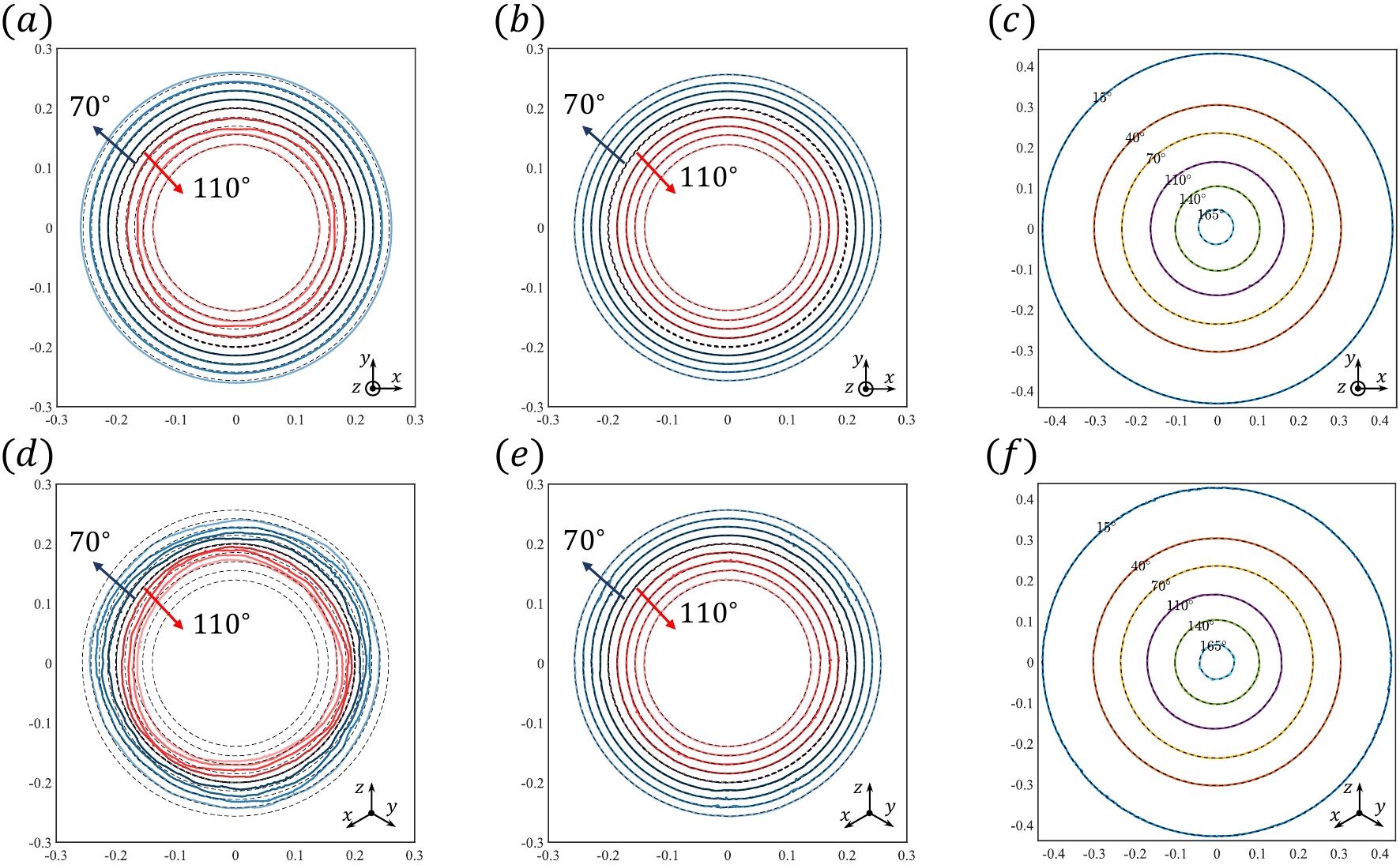}
    \centering
    \caption{
    Contact line evolution using linear (a,d) versus paraboloid fitting (b,e) on flat plates with $\textbf{n}_s = [0, 0, 1]$ (a–c) and $\textbf{n}_s = [1, 1, 1]$ (d–f). Dashed black lines indicate the non-embedded reference. Panels (c,f) show the final contact lines at various $\theta$ compared to the reference circles.
    }
    \label{figr:flatplanedropcontour}
\end{figure}

As shown in Fig.~\ref{figr:flatplanedropcontour}(a) and (b), the contact line contours during spreading reveal distinct differences in symmetry between the two methods. For $\theta = 70^\circ$ and $110^\circ$, the linear fitting method produces visibly asymmetric, ``rectangular'' spreading patterns, especially during retraction. This asymmetry becomes more pronounced at larger contact angles due to inaccurate curvature estimation. In contrast, the paraboloid fitting method preserves excellent rotational symmetry and maintains close agreement with the non-embedded reference case for all tested contact angles. Panel (c) further confirms that at equilibrium, the contact line remains circular under the paraboloid scheme across all angles and embedding depths.

\paragraph{Effect of surface orientation}

To further assess the robustness of the method, an inclined flat surface with $\mathbf{n}_s = [1, 1, 1]$ is considered (Fig.~\ref{fig:planedroplet}(b)). As shown in Fig.~\ref{figr:flatplanedropcontour}(d)–(e), the linear fitting method produces ``triangular'' contact lines that deviate significantly from the non-embedded reference, indicating a breakdown of geometric symmetry. In contrast, the paraboloid fitting approach maintains circular contact lines that closely match the reference case. Figure~\ref{figr:flatplanedropcontour}(f) further demonstrates that the contact line remains symmetric across all tested contact angles, even for inclined and embedded surfaces.

\paragraph{Quantitative summary}

\begin{figure}
    \includegraphics[height=0.36\textwidth]{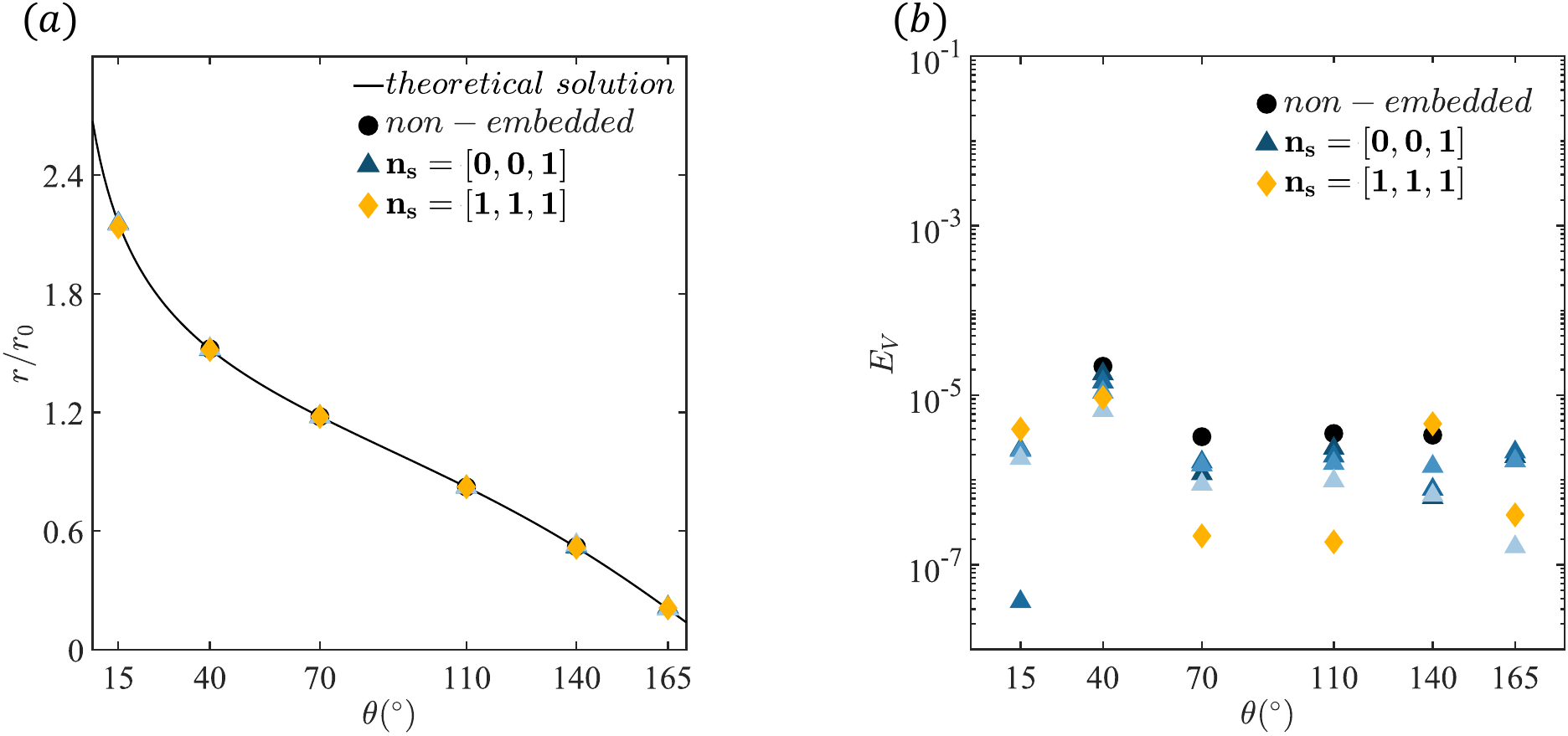}
    \centering
    \caption{
    (a) Equilibrium contact line radius compared to the theoretical prediction (Eq.~\eqref{eq:rt}); (b) global volume loss, $E_V$. Black line: theory; dots: non-embedded case; triangles: $\textbf{n}_s = [0, 0, 1]$ at various $d$; diamonds: $\textbf{n}_s = [1, 1, 1]$.
    }
    \label{figr:flatplanedroprv}
\end{figure}

Finally, Fig.~\ref{figr:flatplanedroprv}(a) summarizes the equilibrium contact line radius across all contact angles and boundary orientations. Results obtained with the proposed method show excellent agreement with theoretical predictions, regardless of the surface normal or embedding depth. Panel (b) presents the corresponding volume conservation performance, showing that the global volume variation remains within $2 \times 10^{-5}$ to $3 \times 10^{-8}$, comparable to the non-embedded baseline.

\subsubsection{Spherical surface}\label{sec:sphere}

\begin{figure}
    \includegraphics[height=0.36\textwidth]{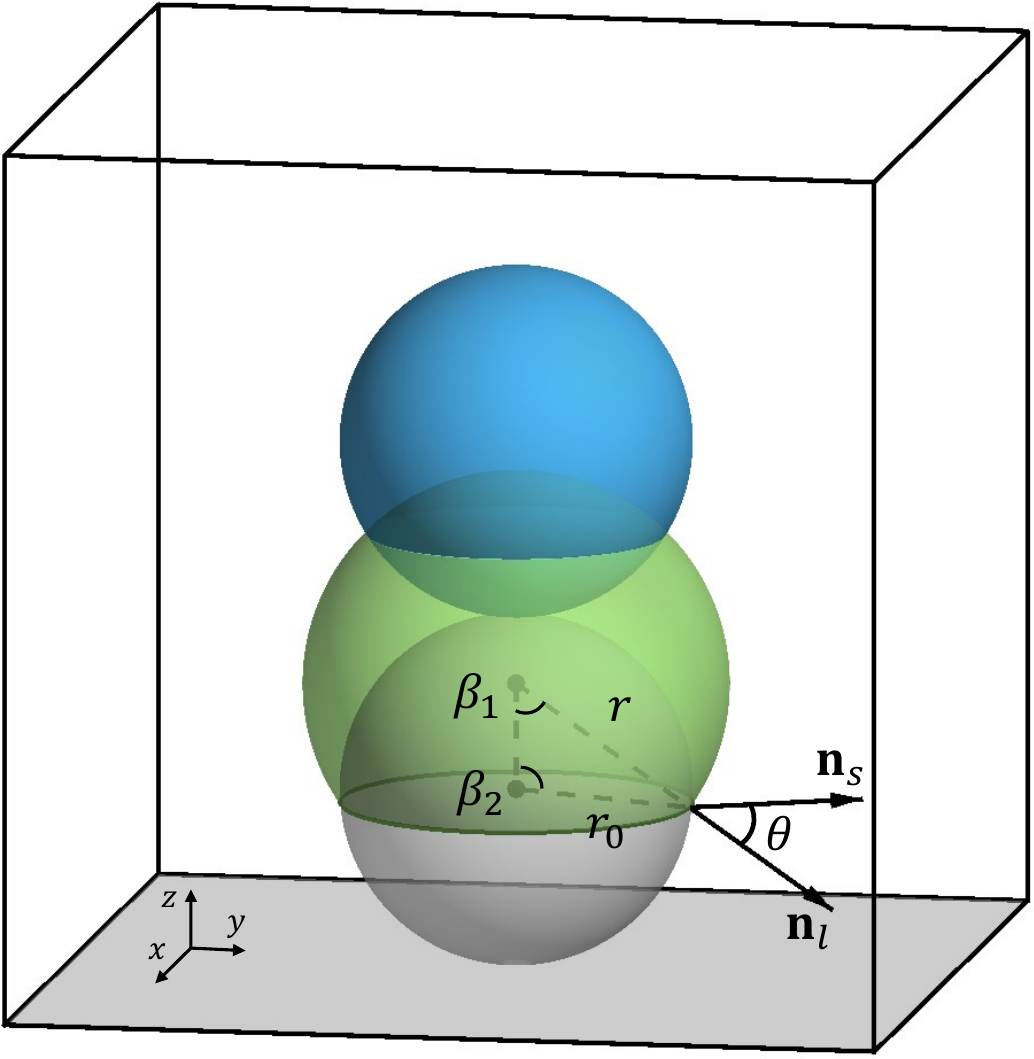}
    \centering
    \caption{
    Schematic of droplet spreading on a spherical surface. $r_0$ and $r$ denote the initial and equilibrium droplet radii, respectively. $\beta_1$ is the angle at the droplet center subtended by the contact point, while $\beta_2$ is the angle at the center of the solid sphere. $\theta$ denotes the prescribed contact angle.
    }
    \label{fig:spheredroplet}
\end{figure}

To further evaluate the accuracy of contact angle imposition on curved geometries, we simulate droplet spreading driven solely by surface tension on a \textbf{spherical surface}. This test examines the method’s capability to impose contact angles on arbitrarily curved boundaries while conserving volume and maintaining symmetry.

The physical parameters and grid resolution are identical to those used in the previous section. A droplet of radius $r_0 = 1$\,mm is initially placed in contact with a solid sphere of equal radius, as illustrated in Fig.~\ref{fig:spheredroplet}. The droplet then spreads under surface tension until it reaches an equilibrium configuration corresponding to a prescribed contact angle $\theta$. Five contact angles are tested: $\theta = 30^\circ$, $60^\circ$, $90^\circ$, $120^\circ$, and $150^\circ$.

The theoretical equilibrium droplet radius $r_t$ and droplet shape are determined from geometric relationships by solving the following nonlinear system:
\begin{equation}\label{eq:spheredroplet}
    \begin{aligned}
        &\beta_1 + \beta_2 + \theta - \pi = 0, \\
        &r_t\sin{\beta_1} = r_0\sin{\beta_2}, \\
        &8\pi r_0^3 = \pi r_t^3(1+\cos{\beta_1})^2(2-\cos{\beta_1}) + \pi r_0^3(1+\cos{\beta_2})^2(2-\cos{\beta_2}),
    \end{aligned}
\end{equation}
where $\beta_1$ and $\beta_2$ are the angles subtended at the centers of the droplet and the solid sphere, respectively (see Fig.~\ref{fig:spheredroplet}). The third equation enforces volume conservation by equating the initial droplet volume to the sum of the spherical-cap volumes associated with the droplet and the solid surface.

\begin{figure}
    \includegraphics[width=.8\textwidth]{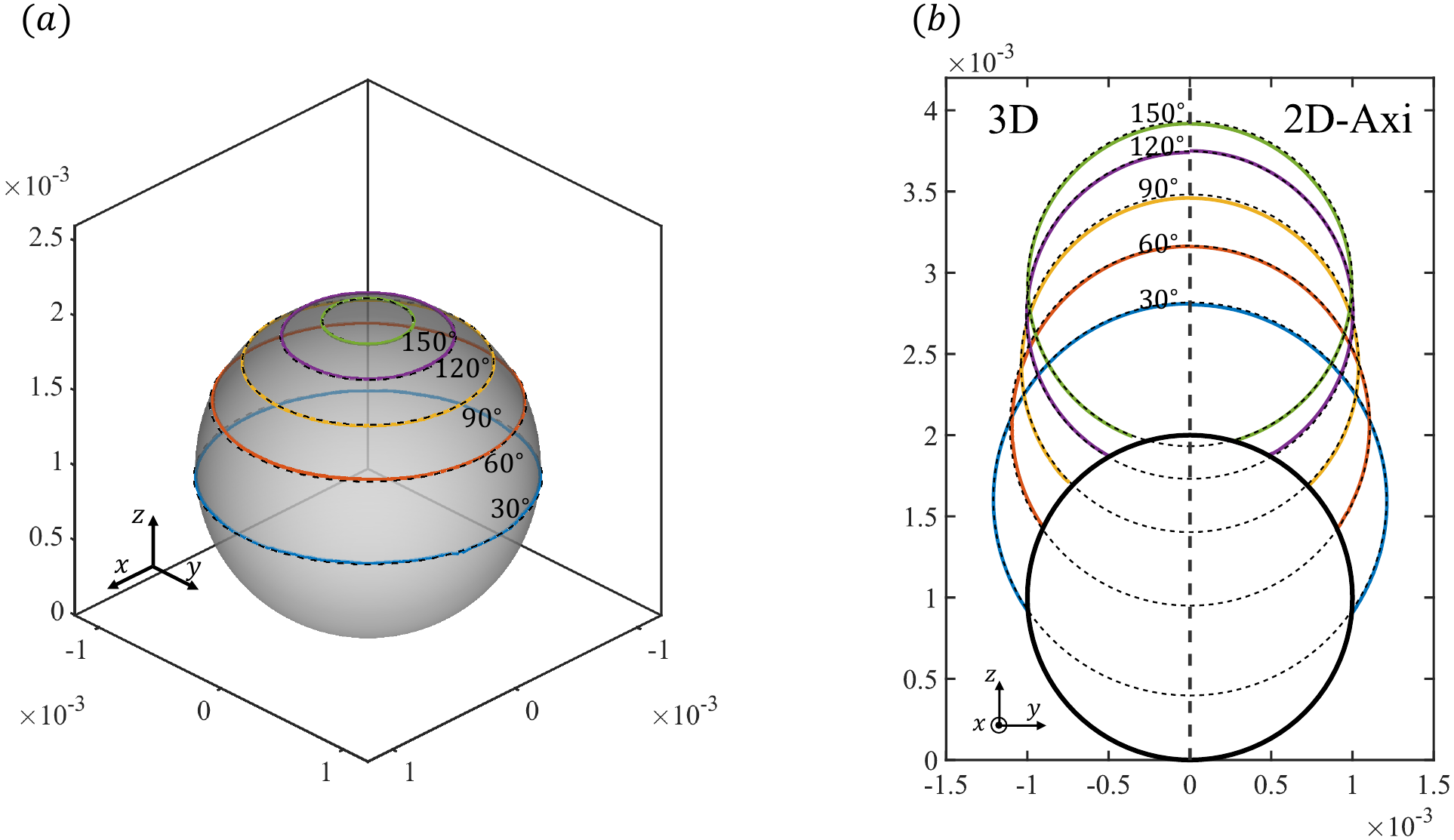}
    \centering
    \caption{
    Droplet spreading on a spherical surface: (a) contact line contours at equilibrium for different contact angles; (b) projection of droplet contours onto the $xz$-plane. Solid lines represent simulation results, while dashed lines correspond to theoretical predictions from Eq.~\eqref{eq:spheredroplet}.
    }
    \label{figr:spheredrop}
\end{figure}

Figure~\ref{figr:spheredrop}(a) shows the equilibrium contact line contours from the 3D simulations (solid lines) at various contact angles, together with the theoretical predictions from Eq.~\eqref{eq:spheredroplet} (dashed lines). Panel (b) presents the projection of these 3D results onto the $xz$-plane and the contours from a 2D axisymmetric simulation using the previously developed method in~\cite{huang20252d}, where the dashed lines denote the theoretical shapes. In all cases, the agreement is excellent, demonstrating both the geometric and physical consistency of the proposed 3D formulation.

Quantitative results are summarized in Table~\ref{tab:spheredroplet}, which reports the relative errors in the computed droplet radius and the global volume loss. The equilibrium radius is determined using a least-squares fit. Across all tested contact angles, the relative error in $r$ remains below $6 \times 10^{-3}$, while the volume loss $E_V$ is consistently smaller than $10^{-6}$. These results confirm the high accuracy of the contact angle enforcement on curved surfaces and the excellent volume conservation properties of the proposed method, even in geometrically complex configurations.

\begin{table}
  \centering
  \caption{
  Relative error in the equilibrium droplet radius and the global volume loss, $E_V$, for droplet spreading on a spherical surface.
  }
  \vspace{5pt}
  \renewcommand\arraystretch{1.5}
  \begin{tabular}{l|ccccc}
    \hline
    $\theta(^\circ)$ & 30 & 60 & 90 & 120 & 150 \\
    \hline
    $|r - r_t| / r_t$ & $1.0\times 10^{-3}$ & $2.6\times 10^{-3}$ & $3.5\times 10^{-3}$ & $5.4\times 10^{-3}$ & $5.5\times 10^{-3}$ \\
    $E_V$ & $9.4\times 10^{-7}$ & $1.0\times 10^{-6}$ & $1.4\times 10^{-7}$ & $5.8\times 10^{-8}$ & $4.0\times 10^{-7}$ \\
    \hline
  \end{tabular}
  \label{tab:spheredroplet}
\end{table}

\subsection{Contact angle hysteresis: shear-driven droplet motion on a flat substrate}\label{sec:hysteresis-sheardriven}

\begin{figure}
    \includegraphics[height=0.36\textwidth]{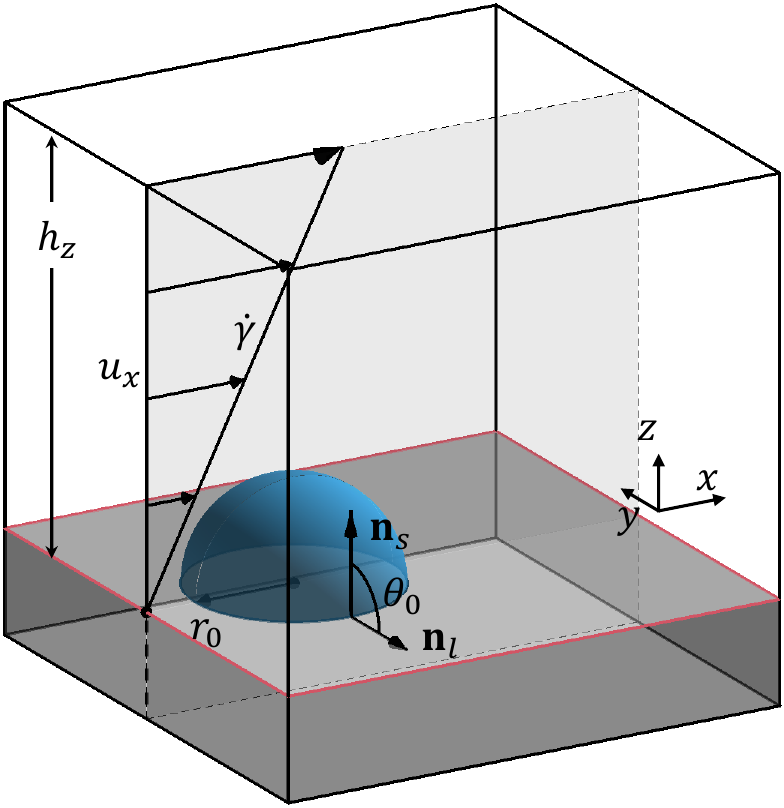}
    \centering
    \caption{
    Schematic of a droplet sliding under shear flow on a flat substrate.
    }
    \label{fig:hysteresis-shear}
\end{figure}

To assess the feasibility of the proposed 3D hysteresis modeling approach, we simulate a droplet sliding under steady shear flow on a flat substrate. The setup follows the benchmark configuration of Ding et al.~\cite{ding2008onset}. As illustrated in Fig.~\ref{fig:hysteresis-shear}, a droplet initially rests on a solid boundary with a prescribed contact angle $\theta_0$ and corresponding contact line radius $r_0$. A fully developed shear flow, defined by $u_x = \dot{\gamma} z$, is applied to drive the droplet motion, while gravity is neglected.

The solid boundary coincides with the $xy$-plane and is represented using the cut-cell method. Exploiting geometric symmetry, only half of the droplet is simulated. Boundary conditions are imposed as follows: along the $x$-direction, both domain faces enforce the shear profile $u_x = \dot{\gamma} z$, $u_y = u_z = 0$; along the $y$-direction, symmetry conditions are applied ($u_y = 0$, $\partial u_x/\partial y = \partial u_z/\partial y = 0$); along the $z$-direction, the top boundary maintains the same shear profile with $u_x = \dot{\gamma} h_z$, and the embedded solid surface at $z=0$ satisfies the no-slip condition.

Two test cases are considered, with the relevant physical parameters summarized in Table~\ref{tab:hysteresisshear-physics}. In both cases, the density ratio $\rho_l/\rho_g$, viscosity ratio $\mu_l/\mu_g$, and shear rate $\dot{\gamma}$ are set to unity. The Reynolds and Weber numbers are defined respectively as $Re = \rho_l \dot{\gamma} a^2 / \mu_l$ and $We = \rho_l \dot{\gamma} a^3 / \sigma$, where $a = \sqrt[3]{3v / 2\pi} = 0.2$\,m is the characteristic droplet length scale. The initial contact line radius $r_0$ is given by
\begin{equation}
  r_0 = a \cdot \sqrt[3]{\frac{2}{(1 - \cos \theta_0)^2 (2 + \cos \theta_0)}} \cdot \sin \theta_0.
\end{equation}
A grid resolution of 32 cells per $r_0$ is employed in all simulations.

\begin{table}
\centering
\caption{
Physical parameters used for shear-driven droplet hysteresis tests.
}
\vspace{5pt}
\renewcommand\arraystretch{1.5}
\begin{tabular}{l|ccccc}
\hline
Case & $Re$ & $We$ & $\theta_0$ & $\theta_a$ & $\theta_r$ \\
\hline
A & 0.684 & 0.095 & $40^\circ$ & $90^\circ$ & $40^\circ$ \\
B & 0.684 & 0.060 & $90^\circ$ & $90^\circ$ & $50^\circ$ \\
\hline
\end{tabular}
\label{tab:hysteresisshear-physics}
\end{table}

\begin{figure}
    \includegraphics[width=.8\textwidth]{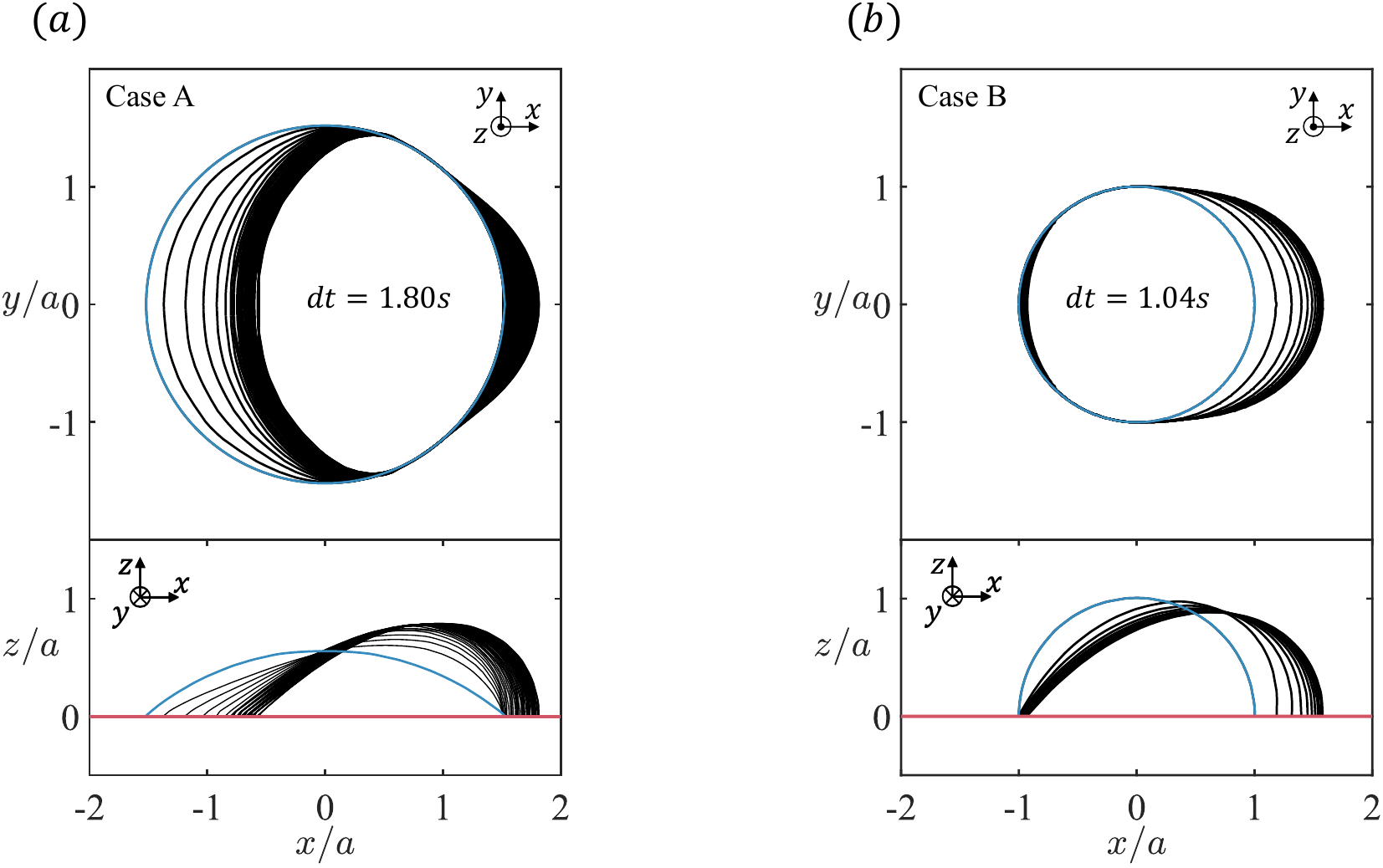}
    \centering
    \caption{
    Droplet evolution under shear flow for (a) Case A and (b) Case B. Top row: contact line contours in the $xy$-plane; bottom row: side views in the $xz$-plane. Time intervals between contours are 1.80\,s for Case A and 1.04\,s for Case B. Blue contours indicate the initial droplet shape.
    }
    \label{figr:hysteresisshear}
\end{figure}

Figure~\ref{figr:hysteresisshear} shows the time evolution of droplet contours during sliding for both test cases. In each panel, the upper row presents the contact line contours on the $xy$-plane, while the lower row shows side views projected onto the $xz$-plane. For Case~A, snapshots are recorded at intervals of 1.80\,s, and for Case~B, at intervals of 1.04\,s. The initial contour is indicated in blue.
In Case~A, the initial contact angle is equal to the receding angle $\theta_0 = \theta_r = 40^\circ$. Consequently, the rear contact line depins immediately and begins to slide, while the front contact line remains pinned until the local angle increases to the advancing value $\theta_a = 90^\circ$. This results in a distinct bulging of the front interface, a reduced droplet length $L$ in the direction of flow, and an increase in droplet height $H$. 
In Case~B, the initial contact angle equals the advancing angle $\theta_0 = \theta_a = 90^\circ$. Here, the front contact line begins to move immediately, while the rear contact line remains pinned until the local angle decreases to the receding angle $\theta_r = 50^\circ$. This leads to an elongated droplet in the streamwise direction and a visibly reduced height compared to the initial configuration.

The final droplet morphologies obtained in both cases are consistent with those reported in previous studies~\cite{ding2008onset,dimitrakopoulos1998displacement}. Table~\ref{tab:hysteresisshear} compares the normalized steady-state droplet—length ($L$), width ($W$), and height ($H$)—with the corresponding results from Ding et~al.~\cite{ding2008onset} and Dimitrakopoulos et~al.~\cite{dimitrakopoulos1998displacement}. The present results show excellent agreement, confirming that the proposed hysteresis scheme accurately reproduces the asymmetric pinning and depinning dynamics of droplets on complex 3D wetting surfaces.

\begin{table}
\centering
\caption{Comparison of normalized droplet dimensions ($L$, $W$, $H$) at steady state.}
\vspace{5pt}
\renewcommand\arraystretch{1.5}
\begin{tabular}{l|ccc|ccc}
\hline
\multirow{2}{*}{Reference} & \multicolumn{3}{c|}{Case A} & \multicolumn{3}{c}{Case B} \\
 & $L$ & $W$ & $H$ & $L$ & $W$ & $H$ \\
\hline
Dimitrakopoulos et al.~\cite{dimitrakopoulos1998displacement} & 2.23 & 3.15 & 0.75 & 2.49 & 2.00 & 0.86 \\
Ding et al.~\cite{ding2008onset} & 2.36 & 2.72 & 0.73 & 2.51 & 1.95 & 0.84 \\
Present study & 2.33 & 2.96 & 0.77 & 2.51 & 2.00 & 0.88 \\
\hline
\end{tabular}
\label{tab:hysteresisshear}
\end{table}

\subsection{Shear-driven droplet motion on a sine-wavy substrate}\label{sec:sinwave}

\begin{figure}
    \includegraphics[height=0.36\textwidth]{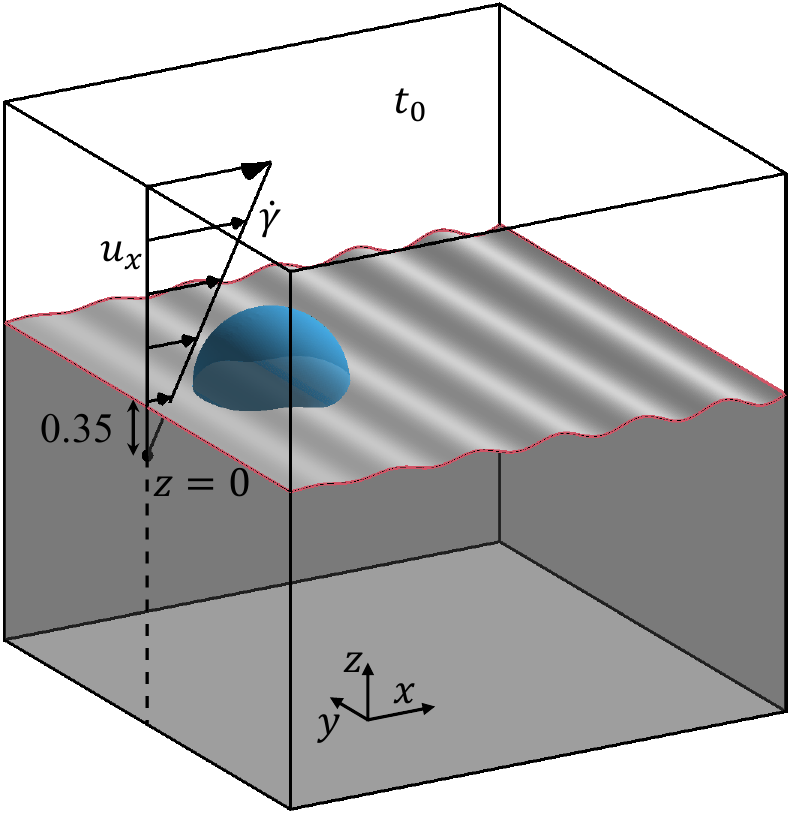}
    \centering
    \caption{
    Schematic of droplet motion under shear flow on a sine-wavy solid substrate.
    }
    \label{fig:sinwave}
\end{figure}

To further assess the robustness of the proposed method on geometrically complex surfaces, we simulate the motion of a droplet subjected to shear flow over a sinusoidal solid substrate. The test configuration follows that of Zhu et~al.~\cite{zhu2025diffuse}, where the solid boundary is defined as
\begin{equation}
    z = 0.035\sqrt{x/2} \sin{(3\pi x)} + 0.35.
\end{equation}
The droplet is initially spherical and tangent to the surface, described by
\begin{equation}
    (x - 1)^2 + y^2 + (z - 0.35)^2 = 0.55^2.
\end{equation}

A shear flow is imposed with velocity components $u_x = 0.2~\text{m/s}$, $u_y = 0$, and $u_z = 0$. The physical parameters are $\rho_l = \rho_g = 1~\text{kg}/\text{m}^3$, $\mu_l = \mu_g = 0.01~\text{kg}/(\text{m} \cdot \text{s})$, and $\sigma = 0.04~\text{N}/\text{m}$. Gravity is neglected in this simulation.

The computational domain extends over $[0, 4] \times [-2, 2] \times [-2, 2]$\,m, with a resolution of 35.2 cells per initial droplet radius (cpr). Boundary conditions are specified as follows: in the $x$-direction, periodic boundaries are applied with a velocity profile $u_x = \dot{\gamma} z$, $u_y = u_z = 0$; in the $z$-direction, shear flow is enforced at the upper boundary with $u_x = 0.2~\text{m/s}$, $u_y = u_z = 0$; and at the embedded solid surface, a fully slip condition is applied in the shear direction, with no penetration and zero normal velocity, satisfying $\partial u_x/\partial x = \partial u_z/\partial z = 0$, and $u_y = 0$. The prescribed static contact angle is $\theta = 60^\circ$.

\begin{figure}
    \includegraphics[width=0.8\textwidth]{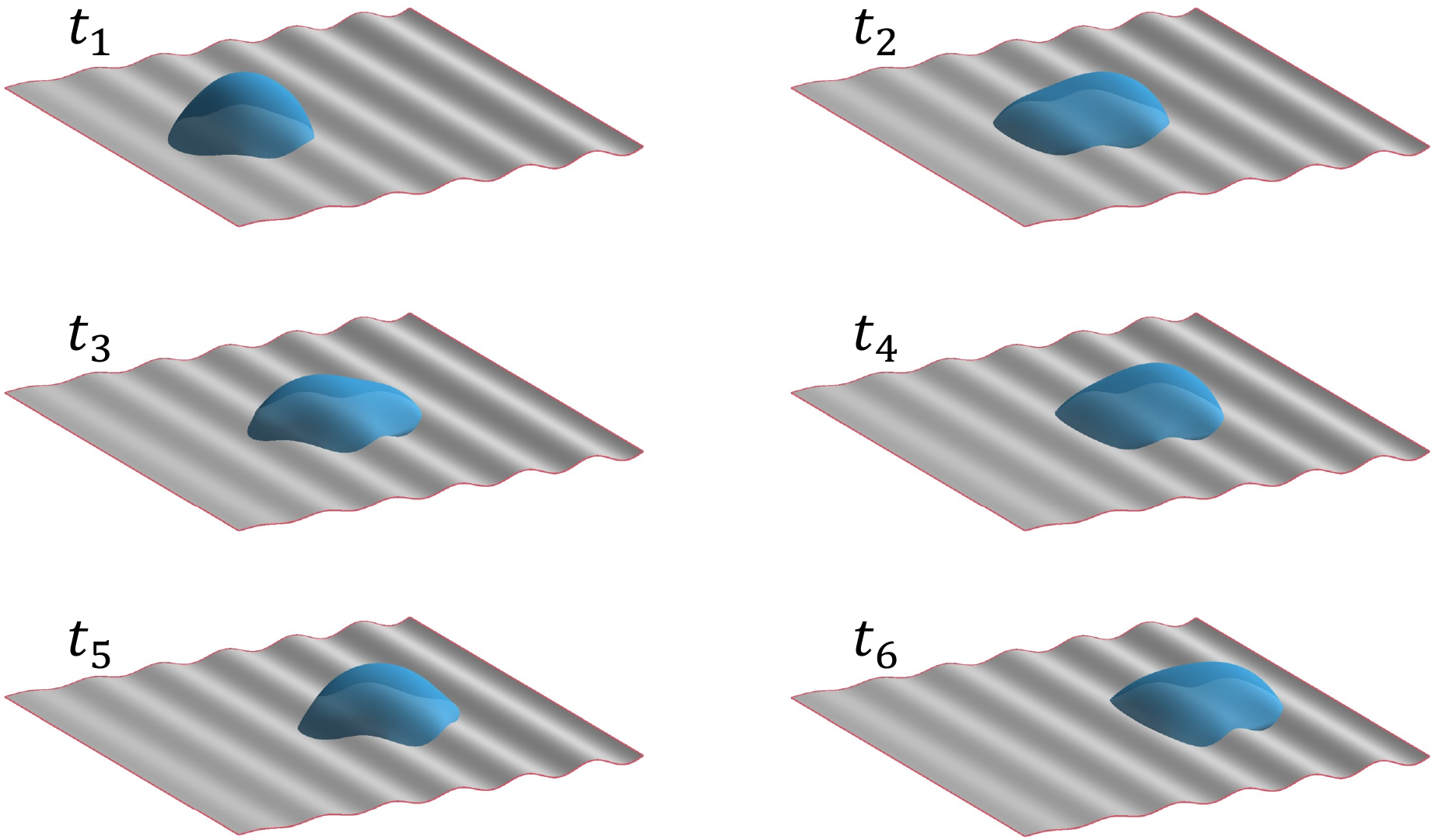}
    \centering
    \caption{
    Snapshots of droplet motion over a sine-wavy surface at time instances $t_i = 1.04,\ 4.35,\ 6.55,\ 8.76,\ 11.51,\ 14.82$\,s.
    }
    \label{figr:sinwave}
\end{figure}

The time evolution of the droplet is illustrated in Fig.~\ref{figr:sinwave}, showing snapshots at six time instances: $t_i = 1.04$, $4.35$, $6.55$, $8.76$, $11.51$, and $14.82$\,s. In the early stage ($t_0$–$t_1$), the droplet spreads symmetrically to relax from its approximate $90^\circ$ initial contact angle toward the imposed $60^\circ$ value. At this stage, motion due to shear is negligible, and the dynamics are dominated by capillary relaxation. Beyond $t_1$, the influence of shear flow becomes significant: the droplet begins to deform and gradually advects in the positive $x$-direction. The sinusoidal substrate topography, however, introduces resistance to motion. As the droplet encounters surface crests, its forward velocity decreases, while fluid accumulates in the troughs, forming localized bulges visible at $t_2$ and $t_4$. Once the advancing front surmounts a crest, the rear region is dragged forward, causing the tail to elongate and sharpen before the droplet clears the crest, as observed at $t_3$ and $t_5$.

The simulation captures the characteristic pinning and depinning behavior induced by surface roughness and the associated evolution of the contact line. The results qualitatively agree with those of Zhu et~al.~\cite{zhu2025diffuse}, confirming that the proposed method accurately and stably resolves complex wetting dynamics on geometrically intricate substrates.

\subsection{3D impact of a droplet on a flat plate with an orifice}\label{sec:orifice}

\begin{figure}
    \includegraphics[height=0.36\textwidth]{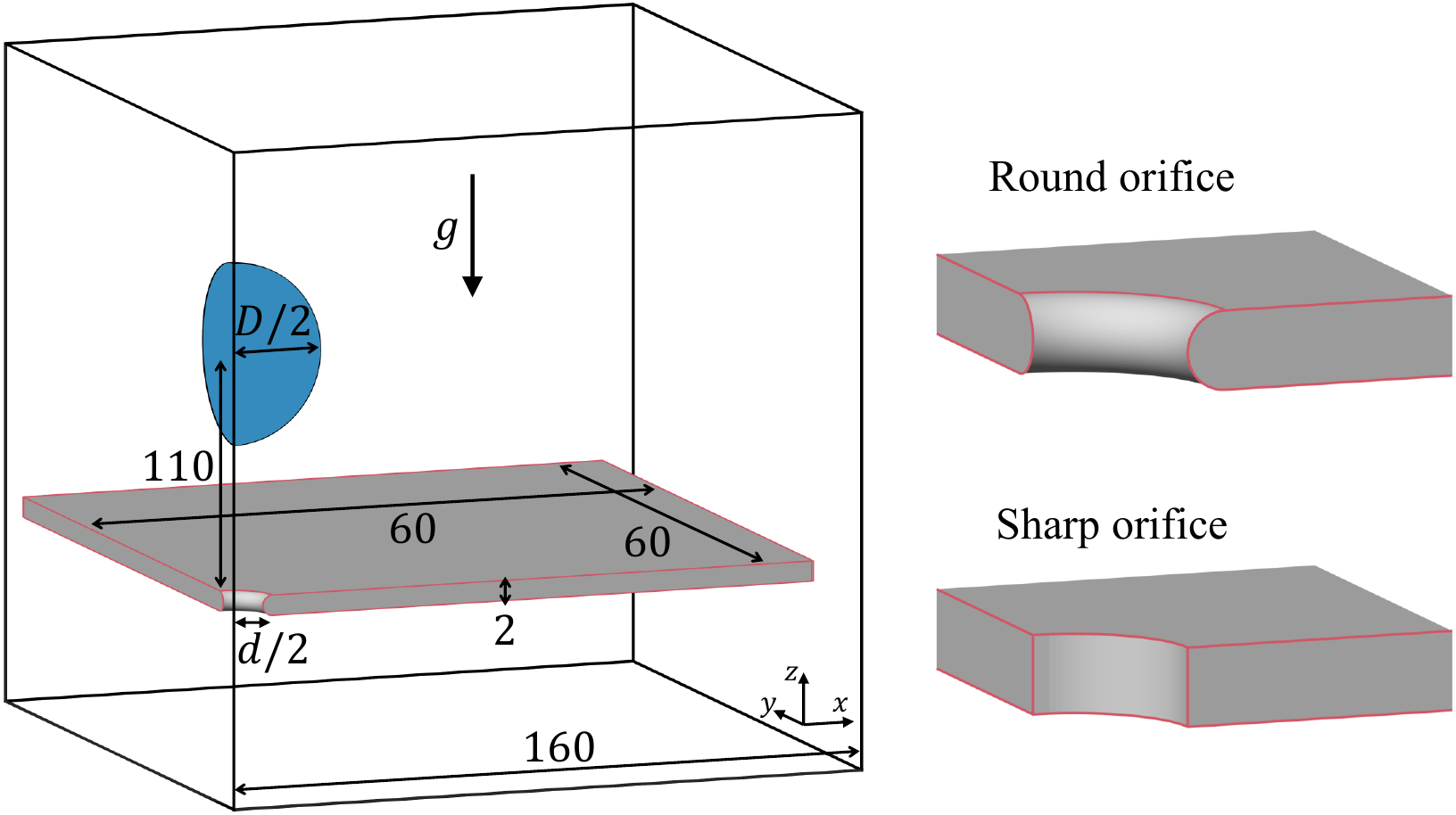}
    \centering
    \caption{
    Schematic representation of droplet impact on a flat plate with either a round or sharp orifice (units: mm).
    }
    \label{fig:orifice}
\end{figure}

To validate the applicability of the proposed method to dynamic contact line problems on complex geometries, we extend our previous 2D simulations~\cite{huang20252d} to a fully 3D configuration of a droplet impacting a plate with an orifice. As illustrated in Fig.~\ref{fig:orifice}, a droplet is released above a horizontal plate and falls freely under gravity, interacting with a centrally located orifice. To reduce computational cost, only one-quarter of the droplet is simulated by exploiting symmetry.

The physical properties are consistent with prior experimental studies~\cite{bordoloi2014drop}: the droplet density and viscosity are $\rho_d = 1130~\text{kg}/\text{m}^3$ and $\mu_d = 0.007~\text{kg}/(\text{m} \cdot \text{s})$, respectively, while the surrounding fluid has $\rho_s = 960~\text{kg}/\text{m}^3$ and $\mu_s = 0.048~\text{kg}/(\text{m} \cdot \text{s})$. Surface tension is $\sigma = 0.0295~\text{N}/\text{m}$, and gravitational acceleration is $g = 9.8~\text{m}/\text{s}^2$. The droplet diameter is denoted by $D$, and it is released from a height of 110\,mm above a circular orifice of diameter $d = 6$\,mm. The plate has a thickness of 2\,mm and half-length and half-width of 60\,mm. The computational domain spans $x \in [0,160]$\,mm, $y \in [0,160]$\,mm, and $z \in [-40,120]$\,mm. Symmetry boundary conditions are applied at $x = 0$ and $y = 0$, the top boundary ($z = 120$\,mm) is open, and all other boundaries are no-slip walls. The grid spacing is approximately $7.8\times 10^{-2}$\,mm, corresponding to a resolution of 64 cpr for a droplet of 10\,mm diameter.
Two types of orifices are considered, following the experimental configurations: (i) a round orifice, through which the droplet passes without wetting the plate surface, and (ii) a sharp orifice, where pinning and contact angle hysteresis occur at the edge. 

\begin{figure}
    \includegraphics[width=1\textwidth]{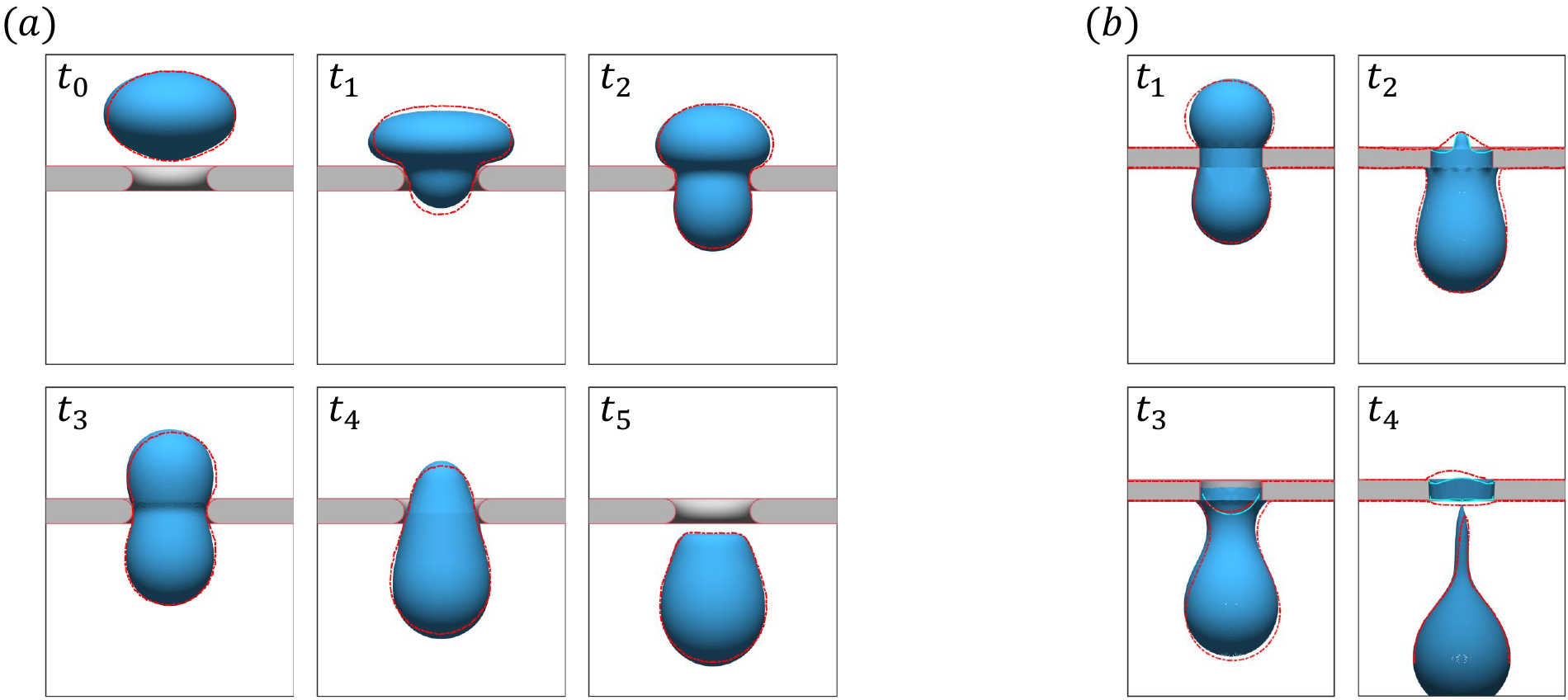}
    \centering
    \caption{
    Snapshots of droplet impact on (a) a round orifice and (b) a sharp orifice. In (a), time instances are normalized as $t_n = t / t_g$ with $t_n = 0$, $0.54$, $1.27$, $1.81$, $2.44$, and $2.96$. In (b), time is normalized as $t_n = t / t_i$ with $t_n = 0.49$, $0.73$, $0.90$, and $1.04$. Blue contours represent simulation results, while red dash-dotted lines correspond to experimental data.
    }
    \label{figr:dropletorifice}
\end{figure}

For the round orifice case, the droplet diameter is $D = 9.315$\,mm, corresponding to a Bond number $Bo = (\rho_d - \rho_s)gD^2/\sigma = 4.9$. Although the plate is mildly hydrophilic in the experiment ($\theta = 55^\circ \pm 1.5^\circ$), the droplet remains non-wetting; accordingly, the contact angle in the simulation is set to $\theta = 180^\circ$. The characteristic gravitational time scale is defined as $t_g = \sqrt{\rho_l D / ((\rho_l - \rho_g)g)}$.
Simulation results are presented in Fig.~\ref{figr:dropletorifice}(a). At $t_0$, just before impact, the droplet decelerates under the influence of the solid surface. By $t_1$, the lower portion of the droplet penetrates the orifice while the upper portion begins to spread laterally. Between $t_2$ and $t_3$, momentum transfer pulls the remaining fluid downward, retracting the lateral spread and facilitating complete passage through the orifice. At $t_5$, the droplet has fully passed through, leaving no residual fluid near the plate, consistent with experimental observations.

For the sharp orifice case, the droplet diameter is increased to $D = 10.307$\,mm, corresponding to $Bo = 6$. Contact angle hysteresis is specified as $[\theta_r, \theta_a] = [42^\circ, 68^\circ]$, based on experimental measurements. To model the observed pinning at the upper edge of the orifice—where minimal spreading occurs—a large advancing angle of $\theta_a = 150^\circ$ is imposed on the top surface. The characteristic time scale for this configuration is defined as $t_i = D^3 / (U_i d^2)$, where $U_i$ is the droplet velocity measured 10\,ms before impact. As shown in Fig.~\ref{figr:dropletorifice}(b), the droplet enters the orifice directly without spreading at the top surface due to the high advancing angle. By $t_2$, the leading edge emerges below the plate and begins to spread, though this is counteracted by mass flow through the narrow orifice. At $t_3$, a crescent-shaped liquid interface forms within the orifice, closely resembling the experimental shape. By $t_4$, the droplet detaches, leaving a small volume of liquid trapped in the orifice, again in agreement with experimental observations.

These 3D simulations are consistent with both our previous axisymmetric results~\cite{huang20252d} and experimental measurements~\cite{bordoloi2014drop}, accurately capturing key interfacial dynamics, including wetting, contact line pinning, and partial entrapment within sharp-edged geometries. The comparison validates the effectiveness of the proposed method for simulating complex droplet impact phenomena involving moving contact lines on intricate solid surfaces.

\subsection{Spontaneous motion of a droplet on a conical surface}\label{sec:conical}

\begin{figure}
    \includegraphics[width=0.36\textwidth]{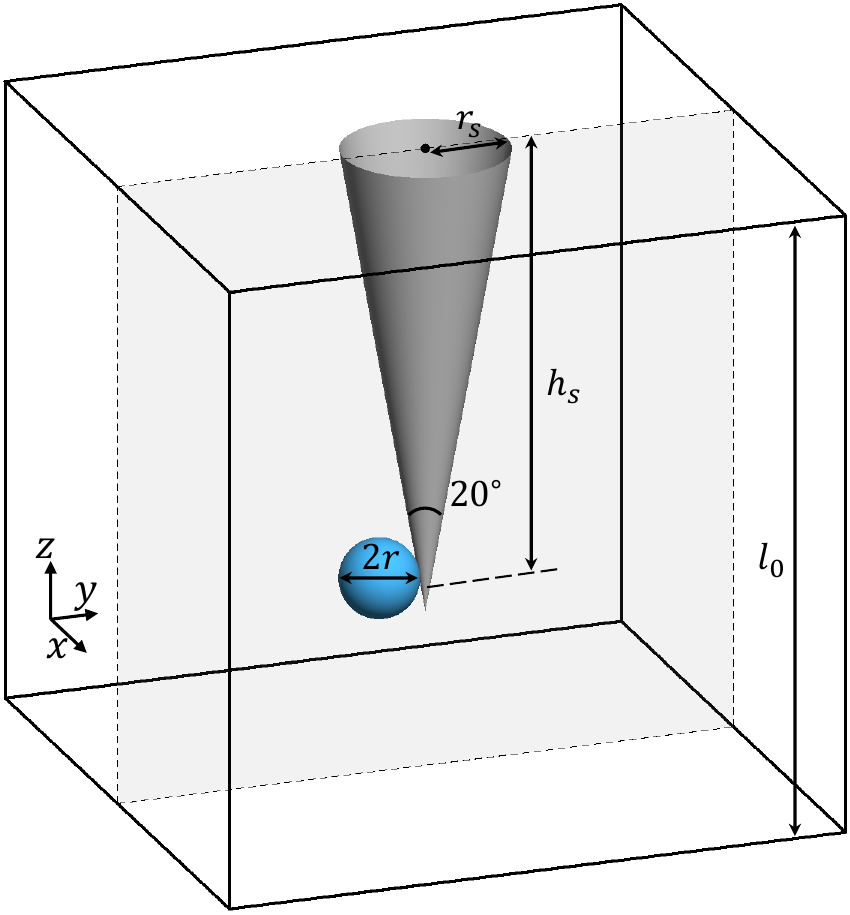}
    \centering
    \caption{
    Sketch of the spontaneous motion of a droplet on a conical surface.
    }
    \label{fig:conical}
\end{figure}

The self-propelled migration of droplets on conical surfaces, driven solely by curvature-induced capillary forces, has attracted significant attention in recent studies. When a droplet is placed on a cone, the geometric gradient in curvature creates an asymmetry between the front and rear contact lines, generating a net capillary force that spontaneously drives the droplet from the apex, where curvature is high, toward the base, where curvature is lower.

Following the setup of Chen et al.~\cite{chen2022simulation}, we simulate this phenomenon using a droplet of initial radius $r = 5~\mu\text{m}$ positioned near the apex of a conical substrate. The droplet properties are $\rho_d = 750~\text{kg}/\text{m}^3$ and $\mu_d = 1.34 \times 10^{-3}~\text{kg}/(\text{m} \cdot \text{s})$, while the surrounding fluid has $\rho_s = 998.2~\text{kg}/\text{m}^3$ and $\mu_s = 1 \times 10^{-3}~\text{kg}/(\text{m} \cdot \text{s})$. Surface tension is set to $\sigma = 0.0528~\text{N}/\text{m}$. Gravity is neglected, as the droplet size is well below the capillary length. The conical substrate has a height $h_s = 60~\mu\text{m}$ and a vertex angle of $20^\circ$. The droplet is initially positioned to intersect a small portion of the cone near $z = l_0 - h_s + r$, where $l_0 = 16r$ defines the domain extent. The computational domain spans $x \in [0, l_0]$, $y \in [-0.5l_0, 0.5l_0]$, and $z \in [0, l_0]$, with a grid resolution of 32 cpr. Due to symmetry, only half of the domain is simulated. The contact angle is prescribed as $\theta = 60^\circ$.

\begin{figure}
    \includegraphics[width=1\textwidth]{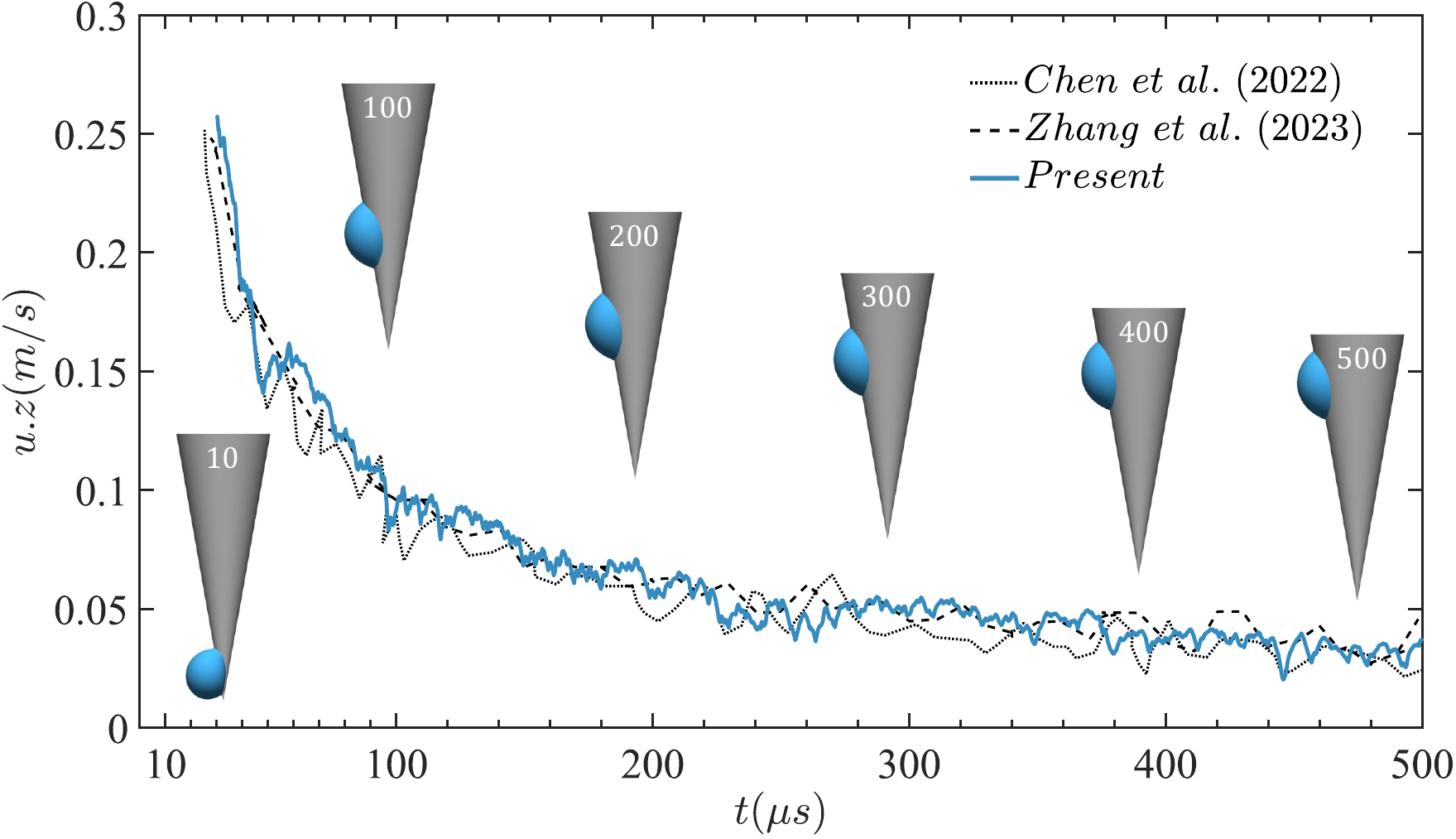}
    \centering
    \caption{
    Time evolution of the droplet climbing velocity on a conical surface. The present results (blue solid line) are compared with those of Chen et al.~\cite{chen2022simulation} (black dotted line) and Zhang et al.~\cite{zhang2023dynamic} (black dashed line). Annotated markers indicate time in $\mu$s.
    }
    \label{figr:conical60}
\end{figure}

The simulation results are presented in Fig.~\ref{figr:conical60}, showing the droplet’s climbing velocity as a function of time. Upon initial contact with the conical surface, the droplet rapidly adheres and attains a high upward velocity due to the strong curvature gradient. As the droplet ascends, the curvature difference between the front and rear contact lines—and consequently the net capillary force—diminishes, while viscous dissipation becomes increasingly significant, gradually reducing the climbing speed. This deceleration behavior is consistent with previous numerical results and theoretical predictions by Chen et al.~\cite{chen2022simulation} and Zhang et al.~\cite{zhang2023dynamic}. The climbing velocities observed in the present simulations exhibit excellent agreement with both studies. Selected droplet shapes during the motion are also visualized in Fig.~\ref{figr:conical60}.

To investigate the influence of geometric perturbations on droplet migration, the conical surface is further modified by superimposing a sinusoidal roughness. The resulting geometry is described by:
\begin{equation}\label{eq:conicaleq}
\frac{r_s^2(z - (l_0 - h_s))^2}{h_s^2} = x^2 + y^2 + \left( \mathcal{A} \cdot r_s \frac{z}{h_s} \sin\left(2\pi \frac{2z}{h_s}\right) \right)^2,
\end{equation}
where $\mathcal{A}$ controls the amplitude of the surface modulation. Simulations are performed for $\mathcal{A} \in [0, 0.35]$ with increments of $0.05$, while maintaining a fixed contact angle of $\theta = 60^\circ$.

\begin{figure}
    \includegraphics[width=1\textwidth]{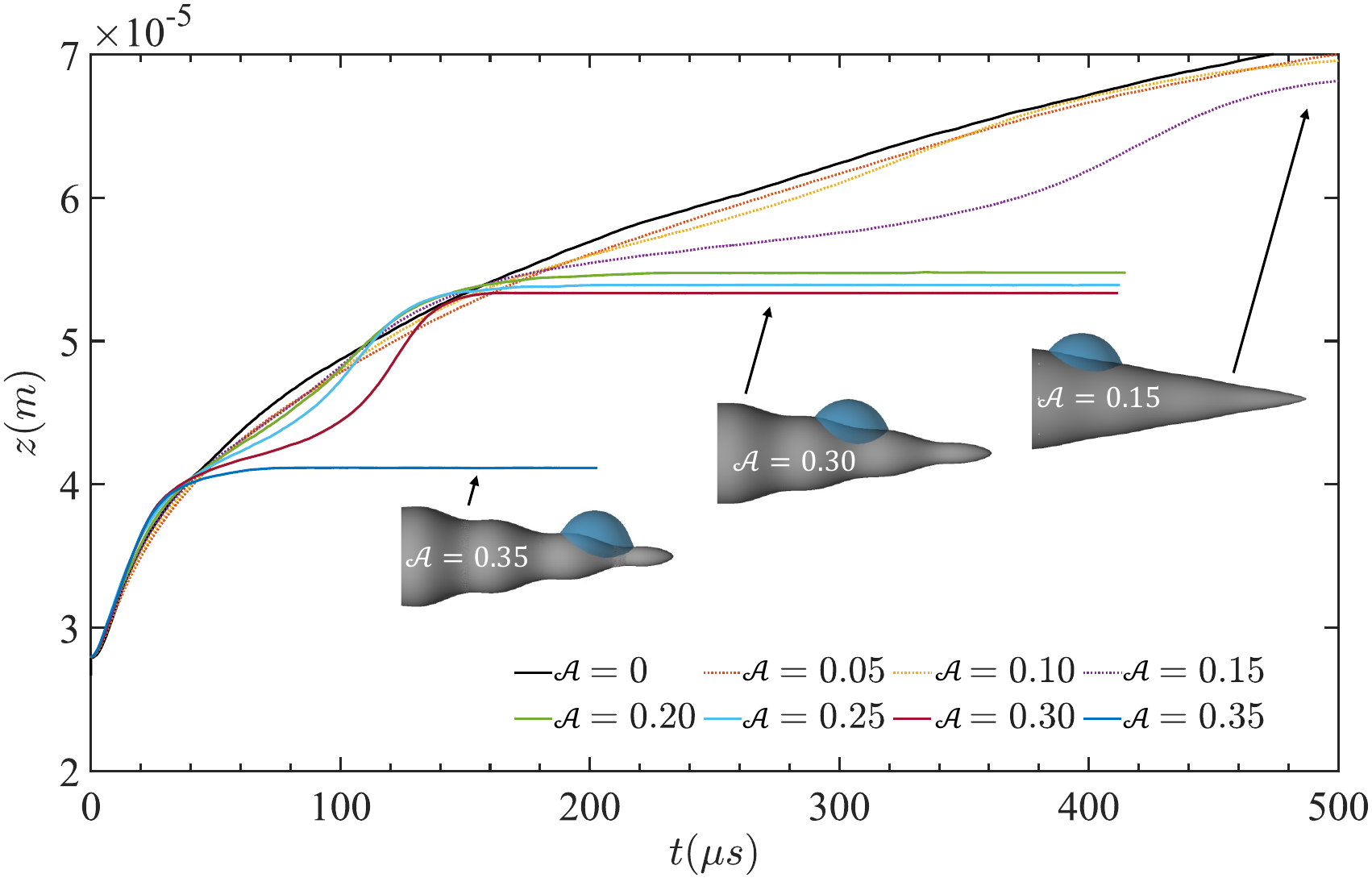}
    \centering
    \caption{
    Time evolution of the droplet position on a conical surface with sinusoidal modulation. Results are presented for different modulation amplitudes, $\mathcal{A}$. Insets show droplet shapes for selected cases.
    }
    \label{figr:conicalsin}
\end{figure}

The time evolution of the droplet's center of mass is shown in Fig.~\ref{figr:conicalsin} for eight representative cases. For $\mathcal{A} = 0.35$, the droplet is unable to overcome the first crest and becomes trapped in the initial trough. As the amplitude decreases to $\mathcal{A} = 0.30$, $0.25$, and $0.20$, the droplet succeeds in crossing the first crest but is arrested in the second trough. When the surface modulation is sufficiently small (e.g., $\mathcal{A} \leq 0.15$), the roughness no longer impedes motion, and the trajectory closely resembles that on a smooth conical surface, indicated by the black solid line. 
The interaction between the droplet and the sinusoidal surface exhibits dynamics similar to those observed in \S~\ref{sec:sinwave}. As the droplet approaches a geometric crest, capillary forces induce fluid accumulation into a ``hill-like'' shape. Once part of the fluid crosses the crest, the remaining mass is rapidly pulled over due to internal capillary-driven motion. These stages are clearly reflected in the displacement curves, which alternate between slower (accumulation) and faster (crossing) phases, highlighting the combined influence of surface geometry and capillarity on droplet migration.

\subsection{Spreading of a droplet on a mesh-shaped surface}\label{sec:mesh}

\begin{figure}
    \includegraphics[width=0.36\textwidth]{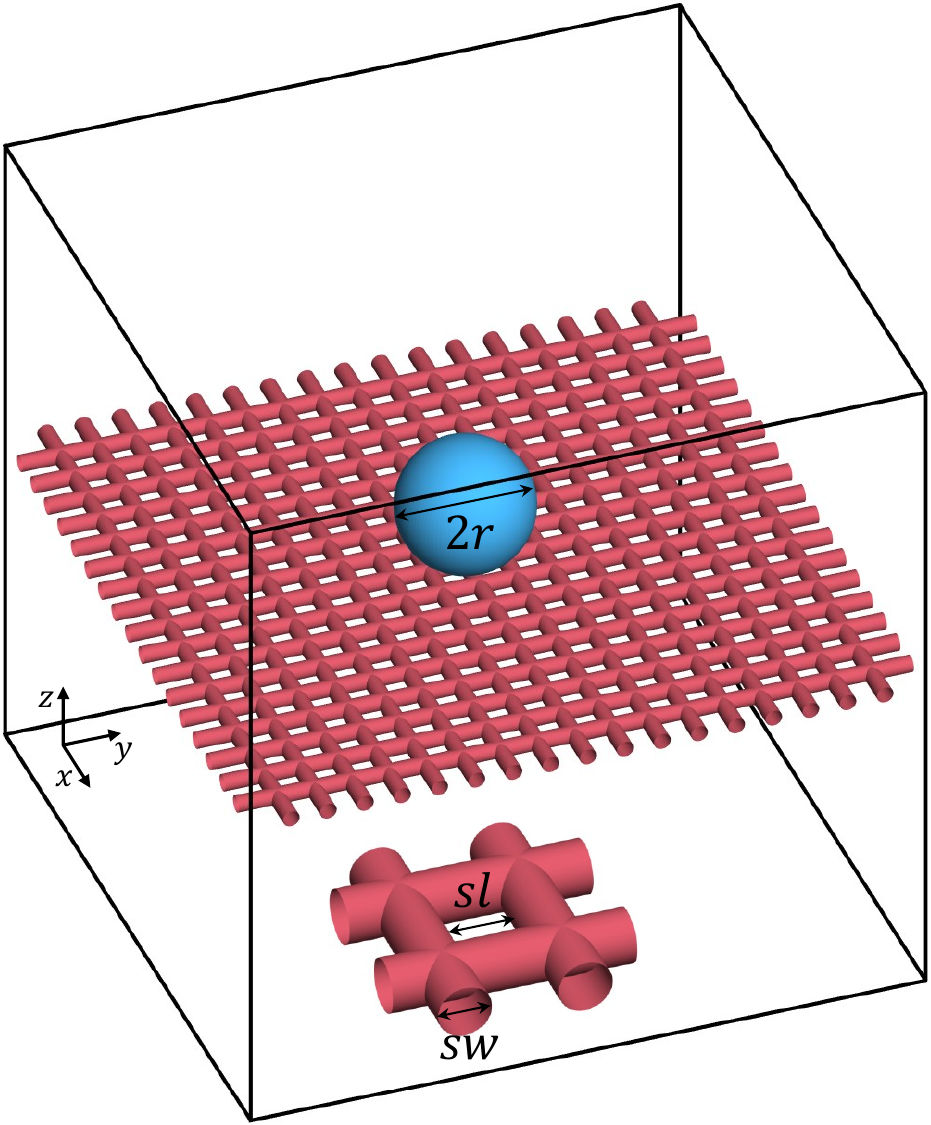}
    \centering
    \caption{Sketch of a droplet spreading on a mesh-shaped surface.}
    \label{fig:mesh}
\end{figure}

The final example investigates droplet spreading on a mesh-shaped surface, representing a geometrically more complex scenario. The mesh consists of intersecting cylindrical structures arranged at uniform intervals along both horizontal and vertical directions in the $xy$-plane. Compared to the previously studied planar and conical surfaces, this configuration presents two primary numerical challenges. First, each cylinder introduces high local curvature, requiring the droplet to spread simultaneously in multiple directions. Second, at the intersections of orthogonal cylinders, abrupt changes in surface geometry induce rapid variations in the curvature of the liquid/gas interface, particularly near corners. These features pose significant difficulties for maintaining numerical stability and accurately capturing the evolution of the contact line.

The initial droplet has a radius of $r = 1 \times 10^{-6}~\text{m}$ and is positioned atop the mesh, as illustrated in Fig.~\ref{fig:mesh}. The mesh itself comprises cylinders with width $sw = r/4$ and spacing $sl = r/3.3$. The physical properties are set as follows: the droplet has a density $\rho_l = 1000~\text{kg}/\text{m}^3$ and viscosity $\mu = 1 \times 10^{-3}~\text{kg}/(\text{m} \cdot \text{s})$, while the surrounding gas has a density $\rho_g = 1~\text{kg}/\text{m}^3$ and viscosity $\mu_g = 2 \times 10^{-5}~\text{kg}/(\text{m} \cdot \text{s})$. The surface tension is $\sigma = 0.072~\text{N}/\text{m}$, and the contact angle is prescribed as $\theta = 30^\circ$. The computational grid resolution is 51.2 cpr.

\begin{figure}
    \includegraphics[width=1\textwidth]{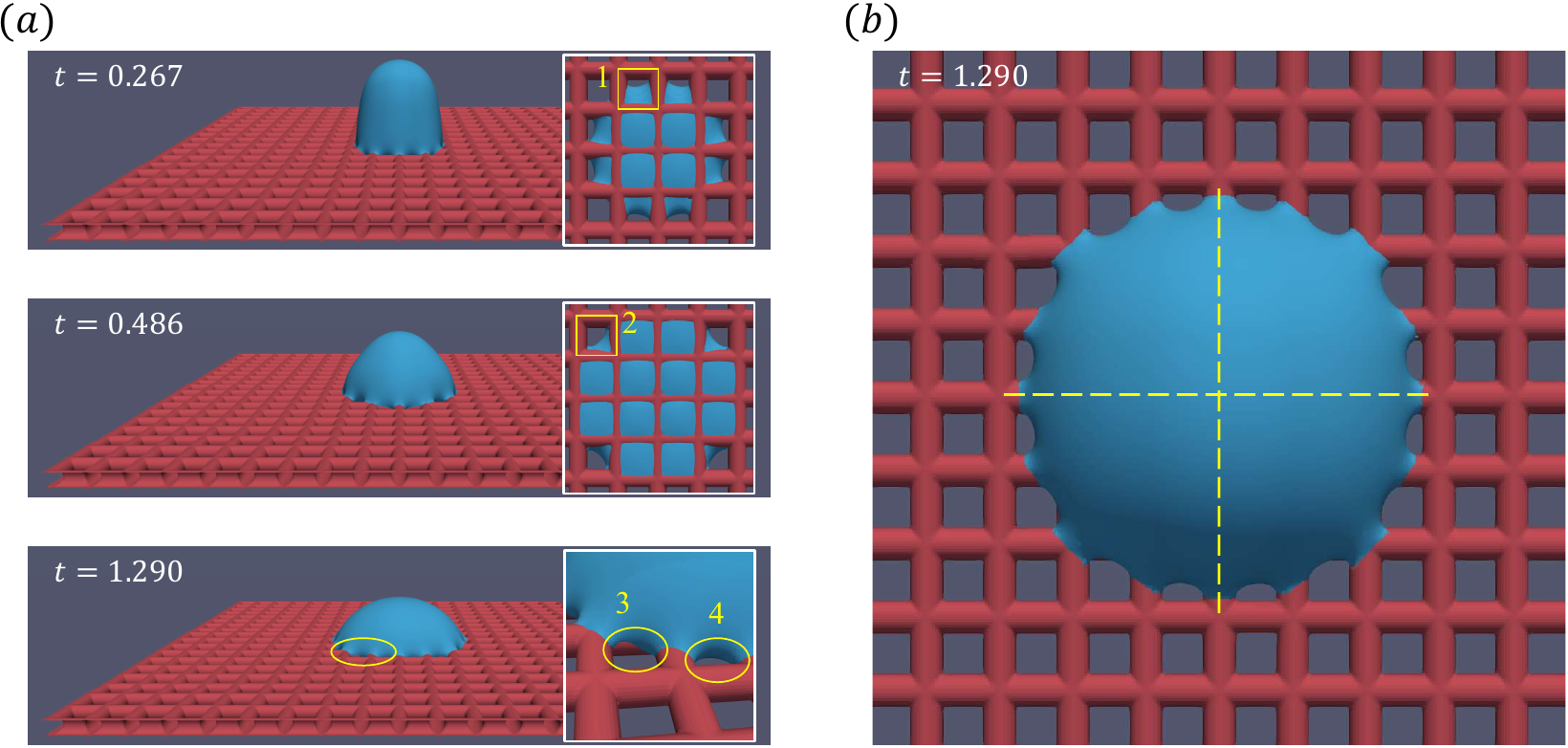}
    \centering
    \caption{Snapshots of droplet spreading on a mesh-shaped surface at various time instances.}
    \label{figr:mesh}
\end{figure}

The simulation results are presented in Fig.~\ref{figr:mesh}, with time normalized using the inertial-capillary time scale $t_i = \sqrt{\rho_l (2r)^3 / \sigma}$. Panel~(a) illustrates the droplet morphology at $t = 0.267$ and $t = 0.486$, highlighting the early-stage spreading behavior. Due to the geometric configuration and prescribed contact angle, the droplet primarily spreads along the tops of the mesh rather than penetrating it. At $t = 0.267$, in grid~1 of panel~(a), the proximity of the right-hand cylinder to the droplet center causes the contact line on that side to advance faster than on the left, resulting in asymmetric spreading. A similar pattern is observed across other grid units. By $t = 0.486$, as shown in grid~2, the contact lines along the right and bottom sides of a unit cell reach the grid corners nearly simultaneously due to symmetry, and then begin extending along the adjacent orthogonal directions. At $t = 1.290$, the droplet attains a quasi-steady state. Grid~3, located diagonally from the center, shows that the droplet partially wets the corner nearest its origin without fully covering the adjacent sidewall. In contrast, in grid~4—closer to the center—the cylindrical wall becomes completely wetted. Panel~(b) presents the top view of the droplet at $t = 1.290$, revealing a symmetric footprint with fourfold patterning aligned with the mesh layout. These results demonstrate that the proposed algorithm accurately captures contact line dynamics on geometrically complex surfaces, maintaining both stability and fidelity even in the presence of abrupt topological transitions.

\section{Conclusion}\label{sec:conclusion}

This work presents a robust and accurate 3D sharp-interface method for simulating moving contact lines on arbitrarily complex solid boundaries. The method is implemented within the open-source \textit{Basilisk} framework and couples the geometric volume-of-fluid (VOF) method for interface tracking with a second-order embedded boundary approach for enforcing boundary conditions in mixed cells. The two main contributions of this study are: (1) the development of a conservative, time-step-unrestricted VOF advection scheme for mixed cells, and (2) the introduction of a novel 3D height-function-based algorithm for imposing contact angle conditions on embedded solid surfaces.

First, we extend the geometric VOF advection scheme from our previous 2D work~\cite{huang20252d} to 3D, adapting it to accommodate arbitrarily complex embedded geometries. A key challenge addressed is the inaccurate estimation of volume fluxes across cell faces in interfacial cells intersected by embedded boundaries, a limitation typically neglected by conventional geometric advection schemes. To correct these flux imbalances, we introduce a volume compensation mechanism that restores consistency in the computed flux of volume fractions. Additionally, a flood-fill algorithm is employed to assist in reconstructing the piecewise-linear liquid/gas interface within irregular mixed cells. To mitigate the severe time-step restrictions imposed by small cut cells—a problem that becomes particularly critical in 3D—we further develop a redistribution strategy. This method transfers excess or deficient volume fractions to neighboring cells in a conservative and physically consistent manner. The resulting algorithm is a fully geometric advection scheme that preserves both local and global volume conservation while entirely eliminating the CFL constraint associated with small cut cells.

Second, we introduce a novel 3D algorithm for imposing static and dynamic contact angle conditions on embedded solid boundaries. This method extends the 2D height-function-based parabola-fitting approach~\cite{huang20252d} to 3D through three key steps: (i) constructing height functions in the fluid portions of mixed cells, even in the presence of irregular geometries; (ii) performing a pre-fitting of a paraboloid surface using all valid HF points surrounding the contact line cell to approximate the local interface profile, thereby enabling the robust identification of well-distributed candidate HF points for the final paraboloid fitting; and (iii) fitting a final objective paraboloid that enforces the prescribed contact angle while supplying missing height function values for curvature estimation near the solid wall. This procedure ensures geometric consistency in the reconstructed interface, improves the stability of curvature estimation, and enables accurate imposition of contact angles on arbitrary 3D embedded surfaces. Furthermore, contact angle hysteresis is incorporated into the methodology based on the principle of volume conservation.

The accuracy and performance of the proposed method are validated through a comprehensive suite of numerical tests. A rotating concentric liquid shell highlights the necessity of conservative advection and the impact of time-step restrictions imposed by cut cells in 3D. A conical cut-out case, extending the classical Zalesak disk problem to include embedded boundaries and three-phase contact lines, demonstrates that the proposed redistribution scheme effectively eliminates cut-cell-induced time-step constraints. Subsequent droplet spreading tests on flat and spherical substrates confirm the accuracy of the paraboloid fitting scheme in enforcing contact angles. Compared to conventional linear fitting methods, the new scheme produces equilibrium contact radii in excellent agreement with theoretical predictions, preserves high contact line symmetry, and exhibits invariance with respect to the position and orientation of embedded solid surfaces—addressing key limitations of existing approaches. Further validation is provided through simulations of shear-driven droplet motion with contact angle hysteresis, where results closely match published benchmarks. Finally, the robustness and versatility of the method are illustrated in complex 3D scenarios, including droplet impact on orifice plates, spontaneous motion on conical and sinusoidally perturbed surfaces, and spreading on mesh-like structures. These cases involve sharp geometric transitions and high interfacial curvatures, yet the proposed method remains stable, accurate, and fully conservative throughout.

In summary, the proposed numerical framework represents a significant advancement for simulating 3D multiphase flows with moving contact lines on complex surfaces. Its geometric accuracy, conservative properties, and capability to handle arbitrary geometries and contact angle models make it a powerful tool for investigating wetting phenomena in both fundamental research and practical engineering applications.

\section*{Acknowledgments}

The authors gratefully acknowledge the support from the National Key R$\&$D Program of China under grants number 2022YFE03130000, 2023YFA1011000 and from the NSFC under grants numbers 12222208, 12588201, and 12472256.

\bibliography{refs}


\end{document}